%% file: xnose.rev.tex
\documentclass[11pt,reqno,oneside]{amsart}
\usepackage{graphics,verbatim,amsmath,amssymb,cite,epsfig}

\numberwithin{equation}{section}

\theoremstyle{plain}
\newtheorem{theorem}{Theorem}[section]

\theoremstyle{definition}

\theoremstyle{remark}

\def\cal{\mathcal}

\def\R{\mathbb{R}}

\def\C{\mathbb{C}}

\def\Z{\mathbb{Z}}

\def\a{\alpha}

\def\e{\epsilon}
\def\g{\gamma}

\def\d{\delta}

\def\s{\sigma}

\def\x{\xi}

\def\inter{\cap}

\def\smallskip{\par\vspace{1mm}}
\def\medskip{\par\vspace{2mm}}
\def\bigskip{\par\vspace{3mm}}

\def\fr#1#2{\frac{#1}{#2}}
\def\smfr#1#2{\tfrac{#1}{#2}}

\newcount\equationcount
\def\thenumber{0}
\def\eq#1{\global\advance\equationcount by 1
   \def\thenumber{\number\equationcount}
                        {$$#1\eqno(\thenumber)$$}}

\newcommand{\Vector}{\mathbb R^n}

\newcommand{\rens}{\mathbb R}

\newcommand{\fracd}[2]{\frac{d #1}{d #2}}
\newcommand{\dis}{\displaystyle}

\begin{document}
\title[Non-ergodicity of the Nos\'e-Hoover thermostatted
harmonic oscillator] {Non-ergodicity of the Nos\'e-Hoover
thermostatted
  harmonic oscillator}


\author{Fr\'ed\'eric Legoll}
\author{Mitchell Luskin}
\author{Richard Moeckel}


\address{Fr\'ed\'eric Legoll \\
Institute for Mathematics and its Applications \\
University of Minnesota \\
207 Church Street SE \\
Minneapolis, MN 55455 \\
U.S.A. and CERMICS and LAMI, ENPC \\
6 et 8 avenue Blaise Pascal \\
77455 Marne la Vall\'ee \\
France}
\email{legoll@lami.enpc.fr}

\address{Mitchell Luskin \\
School of Mathematics \\
University of Minnesota \\
206 Church Street SE \\
Minneapolis, MN 55455 \\
U.S.A.}
\email{luskin@umn.edu}

\address{Richard Moeckel \\
School of Mathematics \\
University of Minnesota \\
206 Church Street SE \\
Minneapolis, MN 55455 \\
U.S.A.}
\email{rick@math.umn.edu}

\thanks{This work was supported in part by
  DMS-0304326, DMS-0500443, the Institute for Mathematics and Its Applications,
  and by the Minnesota Supercomputer Institute. This work is also based on
work supported by the Department of Energy under Award Number
DE-FG02-05ER25706.}

\keywords{Nos\'e-Hoover, Nos\'e, invariant tori, molecular
dynamics, sampling}

\subjclass[2000]{37M25, 65P10, 70F10, 82B80}

\date{March 24, 2006}

\begin{abstract}
The Nos\'e-Hoover thermostat is a deterministic dynamical system
designed for computing phase space integrals for the canonical
Gibbs distribution.  Newton's equations are modified by coupling
an additional reservoir variable to the physical variables. The
correct sampling of the phase space according to the Gibbs measure
is dependent on the Nos\'e-Hoover dynamics being ergodic. Hoover
presented numerical experiments that show the Nos\'e-Hoover
dynamics to be non-ergodic when applied to the harmonic
oscillator.  In this article, we prove that the Nos\'e-Hoover
thermostat does not give an ergodic dynamics for the
one-dimensional harmonic oscillator when the ``mass'' of the
reservoir is large. Our proof of non-ergodicity uses KAM theory to
demonstrate the existence of invariant tori for the Nos\'e-Hoover
dynamical system that separate phase space into invariant regions.

We present numerical experiments motivated by our analysis
that seem to show that the dynamics is not ergodic even for a moderate thermostat mass.
We also give numerical experiments of the Nos\'e-Hoover chain with
two thermostats applied to the one-dimensional harmonic
oscillator. These experiments seem to support the non-ergodicity of the
dynamics if the masses of the reservoirs are large enough and are
consistent with ergodicity for more moderate masses.
\end{abstract}

\maketitle
{
\thispagestyle{empty}

\section{Introduction}
Equilibrium statistical properties of molecular systems
 \cite{mcquarrie,frenkelsmit} are given
 by phase space integrals of the form
\begin{equation}
\label{ps_aver} \langle A \rangle = \int A(q,p) \, d\mu(q,p),
\end{equation}
where $q=(q_1,\ldots,q_M) \in \rens^{nM}$ and $p=(p_1,\ldots,p_M)
\in \rens^{nM}$ denote a set of positions $q_i\in\Vector$ and
momenta $p_i\in\Vector$ of $M$ particles ($n$ denotes the space
dimension), and $A(q,p)$ is an observable, a function defined over
the phase space and related to the macroscopic quantity under
study. If the molecular system is observed at fixed temperature
$\theta,$ then the measure $d\mu$ is the Gibbs measure for the
canonical ensemble \cite{mcquarrie,frenkelsmit}
\begin{equation}
\label{measure_can} d\mu(q,p) =
 \left[ \frac
 {\exp\left({-\beta H(q,p)}\right)}
{\dis{ \int \exp\left({-\beta H(q,p)}\right) \ dq \, dp }} \right]
\, dq \, dp,
\end{equation}
where $H(q,p)$ is the Hamiltonian of the system and is often
simply of the form
$$
H(q,p) = \sum_{i=1}^M \frac{p_i^2}{2m_i} + V(q),
$$
where $p_i^2=p_i\cdot p_i$ and $V(q)$ is the potential energy. The
parameter $\beta$ that appears in (\ref{measure_can}) is related
to the temperature $\theta$
 by $\beta = 1/(k_B \theta)$, where $k_B$ is the Boltzmann constant.
In applications of interest, the number of particles is often very
large ($M \geq 100,000$),
 hence computing integrals such as (\ref{ps_aver}) is
a challenging problem.

Molecular dynamics can be used to compute integrals such as
(\ref{ps_aver}). The method amounts to finding a dynamics on
$(q,p)$ which is ergodic with respect to the measure $d\mu$. As a
consequence, the phase-space average (\ref{ps_aver}) can be
replaced by a time average
\begin{equation}
\label{ergo0} \int A(q,p) \, d\mu(q,p)=\frac{ \dis{ \int A(q,p)
\exp({-\beta H(q,p)) } \, dq \, dp}} {\dis{ \int \exp({-\beta
H(q,p))} \, dq \, dp}}
 = \lim_{T \to +\infty} \frac1T \int_0^T A \left( q(t),p(t) \right) dt
\end{equation}
over a trajectory $(q(t),p(t))_{t \geq 0}$.
 The time average can be approximated by a formula such as
$$
\lim_{T \to +\infty} \frac1T \int_0^T A \left( q(t),p(t) \right)
dt \approx \lim_{N \to \infty} \frac1N \sum_{\ell=1}^N
A(q_\ell,p_\ell),
$$
where $(q_\ell,p_\ell)_{\ell \geq 1}$ is a numerical solution of the chosen
dynamics.

To compute phase space integrals in the canonical ensemble,
several deterministic dynamics have been proposed, such as the
Nos\'e \cite{Nose}, the Nos\'e-Hoover \cite{Hoover}, and the
Nos\'e-Hoover chain dynamics \cite{Martyna92}. More recently, the
Nos\'e-Poincar\'e dynamics \cite{NPoincare99} and the Reversible
Multiple Thermostat method \cite{RMT05} have been proposed.
Stochastic dynamics (such as the Langevin equation) can also be
considered, although we will not discuss them in the following.

The ergodicity condition has not been rigorously proven for any of
the deterministic methods mentioned above. In fact, there is
numerical evidence that shows that the Nos\'e and the
Nos\'e-Hoover methods are not ergodic for some systems
\cite{Hoover,Martyna92,Tuckerman00}, including the one-dimensional
harmonic oscillator. This article further explores non-ergodic
behavior of the  Nos\'e-Hoover dynamics in this simple example.

In Section \ref{Sec_prez}, we present the Nos\'e-Hoover equations
and we recall some of their properties. In Section
\ref{Sec_analysis}, we prove that the Nos\'e-Hoover thermostat
does not give an ergodic dynamics for the Gibbs measure
(\ref{measure_can}) for the one-dimensional harmonic oscillator
when the ``mass'' of the reservoir is large.
Our method is to
apply KAM theory, and more specifically Moser's invariant curve
theorem~\cite{siegelmoser}, to the Poincar\'e return map
and to thus demonstrate the existence of invariant tori that
separate phase space into invariant regions. Finally, in the last
section, we present some numerical experiments with a
Nos\'e-Hoover chain of two thermostats. For large reservoir
masses, results show that the dynamics seems to be non ergodic.
For moderate reservoir masses, they are consistent with
ergodicity.

\section{The Nos\'e-Hoover Thermostat}
\label{Sec_prez} The Nos\'e-Hoover dynamical system \cite{Hoover}
is given by
\begin{equation}
\label{dyn_NH}
\begin{split}
\fracd{q_i}{t} &= \frac{p_i}{m_i} ,\\
\fracd{p_i}{t} &= -\nabla_{q_i} V(q) - \frac{\xi}{Q} \, p_i ,\\
\fracd{\xi}{t} &= \sum_{i=1}^M \frac{p_i^2}{m_i} - nM \beta^{-1},
\end{split}
\end{equation}
where the phase space is described by the physical positions
$q=(q_1,\ldots,q_M) \in \rens^{nM}$ and momenta
$p=(p_1,\ldots,p_M) \in \rens^{nM}$ and an additional variable,
$\xi$, which can be considered as the momentum of the thermostat.
The constant $Q$, which is a parameter of the method, represents
the mass of the reservoir and describes the strength of the
coupling of the reservoir to the physical system. Let us note
that, usually, a second additional variable is introduced
\cite{Hoover}. This variable can be considered as the position of
the thermostat. Since it is decoupled from all the other
variables, we ignore it in the following.

We recall that for the canonical Gibbs measure $d\mu$ given by
\eqref{measure_can} we have
\begin{equation*}
\int \sum_{i=1}^M \frac{p_i^2}{m_i} \, d\mu(q,p)=nM\beta^{-1},
\end{equation*}
so we have that
\begin{equation*}
\lim_{T \to +\infty} \frac1T \int_{0}^{T} \sum_{i=1}^M
\frac{p_i^2(t)}{m_i} \, dt=nM\beta^{-1}
\end{equation*}
for almost all initial conditions for any dynamics on $(q,p)$
which is ergodic with respect to the measure $d\mu.$ Thus, the
right hand side of the Nos\'e-Hoover dynamical equation for
$\fracd{\xi}{t}$ is equal to twice the difference between the
instantaneous kinetic energy of the physical system and the
time-averaged kinetic energy of the physical system at temperature
$\theta=k_B^{-1}\beta^{-1}$ with $nM$ degrees of freedom. Hence,
we see that if the kinetic energy of the physical system is too
high for a sufficiently long time, then the ``thermostat'' added
to the physical momentum equations applies a frictional force to
damp the system.  If the kinetic energy of the physical system is
too low for a sufficiently long time, then the ``thermostat''
added to the physical momentum equations applies an
``anti-frictional'' force to add kinetic energy to the system.

The Nos\'e-Hoover system is not a Hamiltonian system and the
Lebesgue measure $dq\,dp\,d\xi$ is not an invariant measure for
the dynamics.  Instead,  the equations preserve a different
measure which we will now describe. We recall that invariant
measures $\rho(z)dz$ for a general dynamical system
\[
\fracd{z}{t}=f(z)
\]
are determined by the equilibrium equation
\begin{equation}
\label{inv} \operatorname{div}(\rho(z) f(z))=0.
\end{equation}
For the Nos\'e-Hoover system \eqref{dyn_NH} with $z=(q,p,\xi),$ it
is easy to check that
\[
\rho(q,p,\xi)=\exp \left(-\beta \left[ H(q,p) + \frac{\xi^2}{2Q}
\right] \right)
\]
satisfies the condition \eqref{inv}, so a normalized invariant measure for
the Nos\'e-Hoover system \eqref{dyn_NH} is given by
\begin{equation}
\label{measure_NH} d\mu_{\rm NH}(q,p,\xi) =
\left[\frac{\exp \left(-\beta
\left[ H(q,p) + \frac{\xi^2}{2Q} \right] \right) }
{\int \exp \left(-\beta
\left[ H(q,p) + \frac{\xi^2}{2Q} \right] \right)
\, dq \, dp \,
d\xi} \right]\, dq \, dp \,
d\xi .
\end{equation}

We recall that the flow of equations \eqref{dyn_NH} is ergodic (or
metrically indecomposable) with respect to the measure $d\mu_{\rm
NH}$ if the phase space, $\R^3$, cannot be decomposed into two
complementary invariant subsets, each with positive measure
\cite{khinchin}. If the flow of \eqref{dyn_NH} is ergodic
with respect to the measure $d\mu_{\rm
NH}$, then
given any integrable function $A \in L^1(d\mu_{\rm NH})$ the
Birkhoff ergodic theorem shows that for almost all initial
conditions
\begin{equation}
\label{ergo1} \lim_{T \to +\infty} \frac1 T \int_0^{T} A \left(
q(t), p(t),\xi(t)\right) dt = \dis{ \int A \left( q, p ,\xi
\right) d \mu_{\rm NH} },
\end{equation}
where $(q(t),p(t),\xi(t))$ is a solution of the Nos\'e-Hoover
equations~\eqref{dyn_NH}.

If $A=A(q,p)$ is an observable which depends only on the physical
variables $(q,p)$, then
\begin{equation}
\label{ergo5}  \dis{ \int A \left( q, p \right) d \mu_{\rm
NH} }= \dis{ \int A \left( q, p \right) d \mu},
\end{equation}
where $d\mu$ is the Gibbs measure \eqref{measure_can}. To see
this, we observe from Fubini's Theorem that
\begin{equation}
\label{space}
\begin{split}
\int A \left( q, p  \right) d \mu_{\rm NH} &= \frac{\displaystyle\int A(q,p)\exp
\left(-\beta \left[ H(q,p)+
\frac{\xi^2}{2Q} \right] \right) \, dq\,dp\,d\xi
}
{\displaystyle\int \exp
\left(-\beta \left[ H(q,p)+
\frac{\xi^2}{2Q} \right] \right) \, dq\,dp\,d\xi}\\
&=\frac {\dis{ \int A(q,p) \exp({-\beta
H(q,p))} \, dq \, dp}\cdot
\int \exp
\left(-\beta \left[
\frac{\xi^2}{2Q} \right] \right) \,d\xi}
 {\dis{ \int \exp({-\beta H(q,p))} \, dq \,
dp}\cdot
\int \exp
\left(-\beta \left[
\frac{\xi^2}{2Q} \right] \right) \,d\xi}\\
&=\dis{ \int A \left( q, p \right) d \mu}.
\end{split}
\end{equation}
We thus have that if the flow of \eqref{dyn_NH} is ergodic
with respect to the measure $d\mu_{\rm
NH}$ and if the observable $A=A(q,p)$ does not depend on $\xi$, then
for almost all initial
conditions we have from \eqref{ergo1} that
\begin{equation}
\label{ergo22} \lim_{T \to +\infty} \frac1 T \int_0^{T} A \left(
q(t), p(t)\right) dt =\dis{ \int A \left( q, p \right) d \mu}
= \frac {\dis{ \int A(q,p) \exp({-\beta
H(q,p))} \, dq \, dp}} {\dis{ \int \exp({-\beta H(q,p))} \, dq \,
dp}},
\end{equation}
where $(q(t),p(t),\xi(t))$ is a solution of the Nos\'e-Hoover
equations~\eqref{dyn_NH}.

This derivation  has been made under the assumption that the
Nos\'e-Hoover dynamics is ergodic with respect to the measure
$d\mu_{\rm
NH}.$  Numerical experiments show that the equality
(\ref{ergo22}) does not always hold, even for long times $T$ (in
the limit of computationally reachable times). In particular, it
is observed in \cite{Hoover,Martyna92,Tuckerman00} that if the
system under consideration is a one-dimensional harmonic
oscillator, that is, if $n=M=1,$ $m_1=1,$ and $V(q) = \frac12
q^2$, then for lots of initial conditions, there exist $c$ and $C$
with $0 < c < C$ such that the corresponding solution of
(\ref{dyn_NH}) satisfies
\begin{equation}\label{eq_qpbounds}
c \leq q^2(t) + p^2(t) \leq C \qquad\text{for all }t.
\end{equation}
This fact is observed for a wide range of values of $Q$, including
$Q=1$.   This behavior contradicts (\ref{ergo22}), which gives in
this case
\begin{equation}
\label{ave} \lim_{T \to +\infty} \frac1T \int_0^T A \left(
q(t),p(t)\right) \, dt = \frac{\dis{ \int A(q,p) \exp \left(
{-\beta \left( \frac{q^2 +
            p^2}2 \right) } \right) \, dq
\, dp}} {\dis{ \int \exp \left( {-\beta \left( \frac{q^2 + p^2}2
\right)}\right) \, dq \, dp}}.
\end{equation}
For example,  if $A(q,p)$ is a positive function whose support
lies in the disk $q^2+p^2<c$, the left-hand side of (\ref{ave})
will be zero while the right-hand side will be positive.

These numerical experiments give an indication that
 the Nos\'e-Hoover thermostatted harmonic oscillator is
not ergodic. In the next section we will give a rigorous proof of
non-ergodicity if the reservoir mass $Q$ is sufficiently large. We
will apply  KAM theory \cite{siegelmoser} to demonstrate  the
existence of invariant tori that separate the phase space into
invariant regions of positive measure. The projections of these
invariant regions to the $(q,p)$-plane satisfy inequalities of the
form (\ref{eq_qpbounds}) and this explains the numerical
observations. Motivated by our analysis, we show that we can find
initial conditions such that the trajectory does not even sample
the whole ring $\left\{(q,p) \in \R^2: c \leq q^2 + p^2 \leq C
\right\}$ for some $0 < c \leq C$, but only a part of it (see
Figure~\ref{fig_phase_space_NH1} below).

\section{Invariant Tori for the Nos\'e-Hoover Harmonic Oscillator
Dynamics } \label{Sec_analysis}

We now write the Nos\'e-Hoover equations in the case of a
one-dimensional harmonic oscillator. To simplify the notation, let
us assume that the particle mass is $m_1=1$ and the target temperature
is such that $\beta=1/(k_B \theta)=1$. In view
of (\ref{dyn_NH}), the system of differential equations is given
by
\begin{equation}
\label{nh_original}
\begin{split}
\dot q &= p,\\
\dot p &= -q - \e^2 \xi p,\\
\dot \xi &= p^2-1,
\end{split}
\end{equation}
where $\e = 1/\sqrt{Q}$.

We can introduce action-angle variables for the oscillator
\begin{equation}
\label{action_angle}
q = \sqrt{2\tau}\cos\theta \quad\text{and}\quad p =
-\sqrt{2\tau}\sin\theta
\end{equation}
and rescale via $\a = \e \xi$ to get:
\begin{equation}\label{eq_nh}
\begin{split}
\dot \theta &= 1 - \e \a \sin\theta \cos\theta,\\
\dot \tau &= -2\e\tau\a \sin^2\theta,\\
\dot \a &= \e(2\tau\sin^2\theta- 1).
\end{split}
\end{equation}
These equations preserve the volume element
$$
d\Omega = \lambda(\tau,\alpha)\,d\theta\, d\tau \, d\alpha,
$$
where
$$
\lambda(\tau,\alpha) =  e^{-\tau-\alpha^2/2}.
$$

We now make a change of variables as in the averaging method~\cite{siegelmoser}.
Setting
$$
\tau = \hat \tau + \e \hat\tau\hat\a\sin\theta\cos\theta \quad\text{and}\quad
\a = \hat \a - \e \hat\tau \sin\theta\cos\theta
$$
gives a new ODE
of the form:
\begin{equation}\label{eq_nhav}
\begin{split}
\dot \theta &= 1 - \e \hat\a \sin\theta \cos\theta+O(\e^2),\\
\dot {\hat\tau} &= -\e\hat\tau\hat\a +O(\e^2),\\
\dot {\hat\a} &= \e(\hat\tau- 1)+O(\e^2).\end{split}
\end{equation}

The displayed
terms in $\dot {\hat\tau}, \dot {\hat\alpha}$ are the averages
with respect to $\theta$ of the corresponding terms in
(\ref{eq_nh}). These equations preserve the volume element
obtained by transforming $d\Omega$:
$$
d\hat\Omega_\e =
\hat\lambda_\e(\theta,\hat\tau,\hat\alpha)\,d\theta\, d\hat\tau \,
d\hat\alpha,
$$
where
$$
\hat\lambda_\e(\theta,\hat\tau,\hat\alpha) =
 e^{-\hat\tau - \hat\alpha^2/2}+O(\e).
 $$
Indeed,
\begin{eqnarray*}
\hat\lambda_\e(\theta,\hat\tau,\hat\alpha)
&=&
e^{-\tau - \alpha^2/2} \left| \frac{\partial (\tau, \alpha)}{\partial
    (\hat\tau, \hat\alpha)} \right|
\\
&=&
e^{-\tau - \alpha^2/2}\left(1+\frac\epsilon 2\hat\alpha\sin 2\theta
+\frac {\epsilon^2}4\hat\tau\sin^2 2\theta\right)
\\
&=&
e^{-\hat\tau - \hat\alpha^2/2} \ e^{-\frac12 \e^2 \hat\tau^2
  \sin^2 \theta \cos^2 \theta} \left(1+\frac\epsilon 2\hat\alpha\sin 2\theta
+\frac {\epsilon^2}4\hat\tau\sin^2 2\theta\right)
\\
&=&
e^{-\hat\tau - \hat\alpha^2/2}+O(\e).
\end{eqnarray*}
In what follows, the $\hat\,$
will be suppressed and the new variables
 will again be called $(\theta,\tau,\alpha)$.

We will apply the KAM theory \cite{siegelmoser} to the Poincar\'e
return map $P_\e$ of the plane $\Sigma = \{(\theta,\tau,\alpha):
\theta = 0\bmod 2\pi\}$ for the ODE (\ref{eq_nhav}).  This map preserves the area-element:
\begin{equation}\label{eq_areaelement}
d\omega_\e =  \lambda_\e^{0}(\tau,\alpha)\,d\tau\, d\alpha,
\end{equation}
with $\lambda_\e^{0}(\tau,\alpha)=\lambda_\e(0,\tau,\alpha) = e^{-\tau - \alpha^2/2}$.

This can be shown as follows.
Consider a small rectangle $D\subset\Sigma$ with dimensions
$\delta\tau \delta\alpha$ centered at some point $(0,
\tau_0,\alpha_0)\in\Sigma$.  Form a three-dimensional tube ${\cal T}$ by
following the solutions of the ODE (\ref{eq_nhav}) until they
reach the plane $\theta = 2\pi$.  One end of the tube will be the
rectangle $D$ and the other will be its image $P_\e(D)$ under the
Poincar\'e map.

Following ${\cal T}$ forward under the flow of the ODE for a time $\delta
t$ produces a new tube ${\cal T'}$ and the volumes of ${\cal T}$ and ${\cal T'}$ with
respect to the volume element $d\Omega_\e$ are equal.  Now ${\cal T}$ and
${\cal T'}$ differ by two small solid ``cylinders'' with bases $D$ and
$P_\e(D),$ respectively.  Let $\delta = \max(|\delta\tau|, |\delta\alpha|, |\delta t|).$
Then the volume
of the cylinder over $D$ is
$$
\lambda_\e(0,\tau_0,\alpha_0) \, \delta\tau \, \delta\alpha \, \delta\theta
+O(\delta^4) =
\lambda_\e(0,\tau_0,\alpha_0) \ \dot\theta(0,\tau_0,\alpha_0) \,
\delta\tau \, \delta\alpha \, \delta t +O(\delta^4),
$$
where $\dot\theta$ is the first
component of the vectorfield (\ref{eq_nhav}). On the other hand,
the volume of the cylinder over $P_\e(D)$ is
$$
\lambda_\e(2 \pi,\tau_1,\alpha_1) \,
|DP_\e(\tau_0,\alpha_0)| \, \delta\tau \, \delta\alpha \, \delta \theta+O(\delta^4)=
\lambda_\e(2 \pi,\tau_1,\alpha_1) \ \dot\theta(2\pi,\tau_1,\alpha_1) \,
|DP_\e(\tau_0,\alpha_0)| \, \delta\tau \, \delta\alpha \, \delta t+O(\delta^4)
$$
where
$P_\e(\tau_0,\alpha_0) =(\tau_1,\alpha_1)$ and where $
|DP_\e(\tau_0,\alpha_0)|$ is the Jacobian determinant of the
Poincar\'e map.   Note that $\dot\theta(0,\tau,\alpha) =
\dot\theta(2\pi,\tau,\alpha) =1$ since this holds for the ODE
(\ref{eq_nh}) and is preserved by the coordinate change leading to
(\ref{eq_nhav}). It then follows from letting $\delta\to 0$ that
$$
\lambda_\e(0,\tau_0,\alpha_0)=|DP_\e(\tau_0,\alpha_0)| \ \lambda_\e(2 \pi,\tau_1,\alpha_1)
$$
which proves
that the Poincar\'e map $P_\e$ preserves the area element
$\lambda_\e^{0}(\tau,\alpha) \, d\tau \, d\alpha$ as claimed.

We will use the version of Moser's invariant curve theorem from ``Lectures on
Celestial Mechanics" by Siegel and Moser
\cite[sections 32--34] {siegelmoser}.  That
theorem starts with a real-analytic map $P$ of the form
\begin{equation}\label{eq_SM}
\begin{split}
x_1 &= x + \gamma\,y +f(x,y),\\
y_1 &= y + g(x,y),
\end{split}
\end{equation}
where $f$ and $g$ are periodic in $x$ with period $2\pi$.  In the
application they have in mind, $x$ and $y$ are respectively the angle
and radius in a polar coordinate system near a fixed point.   If
$f=g=0,$ then the map reduces to a standard twist map where the radial
variable is preserved while the angular variable is rotated by an
amount which depends on the radius.  They assume a {\it twist
condition}, $\gamma\ne 0$, so the amount of rotation really does
change.    The theorem shows that if $f$ and $g$ are small, then
some of these invariant circles will persist.

In order to prove
this, they also assume that $P$ satisfies the {\it
curve-intersection property}.  This means that for any simple closed
curve, $C$,
of the form
$$
y=\psi(x)\qquad\text{where}\quad \psi(x+2\pi)=\psi(x),
$$
we have that $C\inter P(C)\ne \emptyset$.

 Such a curve represents a simple closed
curve around the fixed point for the map before changing to polar
coordinates.  If this map preserves an area element, then such a
loop cannot map completely inside or completely outside itself and
so the curve intersection property will hold.

To state the theorem precisely, fix an annulus  $a\le y\le b$ with
$b-a=1$ (this condition can always be achieved by rescaling $y$).
Since $f$ and $g$ are real analytic, they extend to complex
analytic functions on some complex neighborhood $D$ of $\R\times
[a,b]$ in $\C^2$.  We also specify a so-called {\it Diophantine
condition} on the rotation numbers  of the unperturbed circles.
Recall that if $f(x)$ satisfies $f(x+2\pi) = f(x) + 2\pi$ and
defines a homeomorphism of the circle $\R\bmod{2\pi}$, then the
rotation number $\omega = \lim_{n\rightarrow\infty}\fr{1}{2\pi
n}(f^n(x)-x)$ exists and is independent of the initial condition,
$x$.  For the unperturbed circles of constant $y$ we have
$f^n(x)-x = n \gamma y$, so $\omega = \frac{\gamma y}{2\pi}$. Fix
any constants $c_0>0$ and $\mu\ge 2$ and consider rotation numbers
satisfying:
\begin{equation}\label{eq_Diophantine}
\left |l\, \omega - k \right| = \left |\frac{l \gamma
y}{2 \pi} - k \right| \ge \fr{\gamma c_0}{l^\mu} \qquad
\text{for all $k, l \in \Z, l>0$}.
\end{equation}
 It is shown in Siegel-Moser that for $c_0$ sufficiently small, the set
of $y$ which satisfy this condition forms a positive measure
Cantor set in the interval $[a,b]$ whose measure tends to  $b-a=1$
as $c_0\rightarrow 0$.  Moser's theorem states that if $f$ and $g$
are sufficiently small, then for each such $y$ there is a nearby
invariant curve of $P$ with the same rotation number.  More
precisely,

\begin{theorem} Let $P$ be a real analytic map of the form (\ref{eq_SM})
which extends to a complex domain $D$.   Assume that $P$ is
periodic in $x$ of period $2\pi$ and has the curve intersection
property.  Fix a Diophantine condition (\ref{eq_Diophantine}) and
an annulus $a\le y \le b$ where $b-a=1$.

Given any $\tilde\e>0$ there is a $\d>0$ such that if
\begin{equation}\label{eq_delta}
|f| + |g| < \d \g,
\end{equation}
then for every $y_0 \in [a,b]$ which satisfies
(\ref{eq_Diophantine}), $P$ has an invariant curve of the form
$y=\psi(x)$ where $ \psi(x+2\pi)=\psi(x)$ with $|\psi(x) - y_0| <
\tilde\epsilon$  for all $x$.  Moreover, the restriction of $P$ to
this curve has rotation number $\omega = \frac{\gamma y_0}{2
\pi}$.  Here $|\cdot|$ denotes the sup norm in the complex domain
$D$.  The constant $\d$ depends on $c_0, \mu, \tilde\e, D$ but not
on $\gamma$.
\end{theorem}

To apply this result to $P_\e$, some further coordinate changes
will be needed.   First, note that
$$P_\e(\tau,\alpha) = Q_\e(\tau,\alpha) + O(\e^2),$$
where $Q_\e$ is the time-$2\pi$ advance map of the differential equation:
\begin{equation*}
\begin{split}
\dot \tau &= -\e\tau \a,\\
\dot \a &= \e(\tau- 1),\end{split}
\end{equation*}
or equivalently, the time-$2\pi\e$ advance map of
\begin{equation}\label{eq_unpert}
\begin{split}
\tau^\prime &= -\tau \a,\\
\a^\prime &= \tau- 1.
\end{split}
\end{equation}
This follows because the return time to the section $\Sigma$ is
$$T(\tau,\alpha) = 2\pi+O(\e).$$

The ODE (\ref{eq_unpert}) has an integral
\begin{equation}
\label{defG}
G(\tau,\alpha) = \tau - \ln\tau + \smfr12 \alpha^2 -1.
\end{equation}
Furthermore there is an equilibrium point at $(\tau_0,\alpha_0) =
(1,0)$ on the level set $G=0$.  All of the other level sets in the
half-plane $\tau>0$ are simple closed curves around the
equilibrium point (see Figure~\ref{fig_Gcontour}).  These are all
invariant curves for the map $Q_\e,$ and the goal is to show that
the actual Poincar\'e maps $P_\e$ with $\e$ sufficiently small
also have such invariant curves.  We note that the existence of
invariant curves for $P_\e$ corresponds to the existence of
invariant tori for the 3D flow of (\ref{eq_nhav}).   These tori
separate the phase space into invariant 3D regions, showing that
the flow is not ergodic.   Note that for an invariant torus close to a level set
$G=g_0 >0$ we have
$$\tau - \ln\tau  \le 1+ g_0 + O(\e).$$
Since $\tau = \fr12(p^2+q^2)$, it follows easily that  the
projection of the torus to the $(q,p)$-plane satisfies bounds of
the form (\ref{eq_qpbounds}).

It is worth noting that even without using the KAM theory, the
existence of the integral $G$ for the averaged system
(\ref{eq_unpert}) shows that the convergence of ergodic averages
in (\ref{ave}) would be slow for $\epsilon$ small. Indeed, the
standard averaging theory shows that the trajectories of the
actual system remain $\epsilon$ close to those of the averaged
system on a time scale of order $1/\epsilon$.   The KAM theory
shows that the trajectories which lie on the invariant tori are
close to the averaged ones for {\it all} time.  Moreover, other
trajectories are trapped between the tori and so their $G$
coordinates are prevented from wandering very far.

To get the Poincar\'e map into the form (\ref{eq_SM}), we introduce
action-angle variables around the equilibrium point. The integral
$G$ provides a natural radial coordinate or action variable.   To
construct the corresponding angle variable, let $T(g)$ denote the
period of the periodic solutions of (\ref{eq_unpert}) which
correspond to the level curve $G=g$. We define an angular coordinate
$\phi$ to be the time along this orbit multiplied by $2\pi/T(g)$,
taking the initial point $\phi=0$ along the $\tau$-axis to the
right of the equilibrium point.  By definition, the ODE
(\ref{eq_unpert}) will become:
\begin{equation*}
\begin{split}
\phi^\prime &= 2\pi/T(G),\\
G^\prime &= 0,
\end{split}
\end{equation*}
in these coordinates and so the time-$2\pi\e$ advance map
$Q_\e(\phi,G)=(\phi_1,G_1)$ is
\begin{equation}\label{eq_unpertmap}
\begin{split}
\phi_1 &= \phi + 2\pi\e/T(G),\\
G_1 &= G.
\end{split}
\end{equation}

\begin{figure}
\includegraphics[width=7cm]{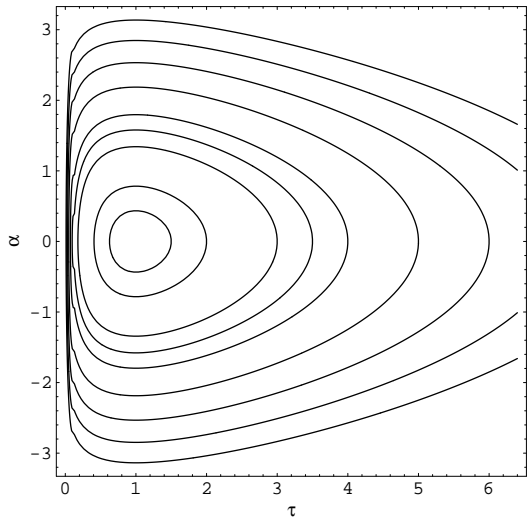}
\caption{Level
curves of $G$ (see (\ref{defG})), the integral for the averaged
equations (\ref{eq_unpert}).}\label{fig_Gcontour}
\end{figure}

The Poincar\'e maps $P_\e$ can be viewed as $O(\e^2)$
perturbations of (\ref{eq_unpertmap}). To see this, let us first
show that $P_\e$ has a fixed point $(\tau_\e,\alpha_\e) = (1,0)
+O(\e).$ Since $Q_\e$ is the time-$2\pi\e$ advance map of
(\ref{eq_unpert}), we compute that
$$
Q_\e(\tau,\a) = (\tau-2 \pi \e \tau \a,\a+2 \pi \e (\tau-1)) + O(\e^2),
$$
hence
\begin{equation}\label{implicit}
P_\e(\tau,\a) - (\tau,\a) = 2 \pi \e R(\e,\tau,\a),
\end{equation}
with
$$
R(\e,\tau,\a) = (-\tau \a + \e u(\e,\tau,\a),\tau - 1 + \e v(\e,\tau,\a))
$$
for some smooth functions $u$ and $v$. We then apply the
implicit function theorem to $R(\e,\tau,\a)$ to continue the
solution in $\e$ from $R(0,1,0) =
0$. We can compute
$$
\frac{\partial R}{\partial (\tau,\a)}(0,\tau,\a) =
\begin{pmatrix}
-\a & -\tau \\
1 & 0
\end{pmatrix},
$$
so the matrix $\dis{ \frac{\partial R}{\partial (\tau,\a)}(0,1,0) }$ is
invertible. As a consequence, the implicit equation $R(\e,\tau_\e,\a_\e)
= 0$ defines a function $\e \mapsto (\tau_\e,\a_\e)$ from a neighborhood of
0 to a neighborhood of $(1,0)$. In view of (\ref{implicit}), we see that
$(\tau_\e,\a_\e)$ is a fixed point of $P_\e$. From the equation
$R(\e,\tau_\e,\a_\e) = 0$, we obtain that $(\tau_\e,\alpha_\e) = (1,0)
+O(\e).$

After a translation of the coordinates, one can assume that all
maps $P_\e$ for $\e$ sufficiently small fix the point $(1,0)$.
Now use the same coordinates $(\phi, G)$ as for the unperturbed
map $Q_\e$.  $G$ is not an integral but one has that $P_\e (\phi,G)=
(\phi_1,G_1)$ with
\begin{equation}\label{eq_poincmap}
\begin{split}
\phi_1 &= \phi + 2\pi\e/T(G)+O(\e^2),\\
G_1 &= G+O(\e^2).
\end{split}
\end{equation}
The fact that $P_\e$ preserves the area element $d\omega_\e$ (see
(\ref{eq_areaelement}))
implies that (\ref{eq_poincmap}) satisfies the  curve-intersection
property.   We also need  a  twist condition on the period $T(G).$
It will be shown below that $T'(G) >0$, that is,
the period of the periodic solutions
increases as we move out from the
equilibrium point.  This implies that one can replace the action
coordinate $G$ by the period $T(G)$ or its reciprocal $1/T(G)$.
To match the notation in Siegel-Moser, let $y = 1/T(G)$ and $x =
\phi$. Then one finds that the Poincar\'e maps take the form
\begin{equation}\label{eq_poincmap2}
\begin{split}
x_1 &= x + 2\pi\e\,y +\e^2 \tilde f(x,y,\e),\\
y_1 &= y + \e^2 \tilde g(x,y,\e).
\end{split}
\end{equation}
This is of the form (\ref{eq_SM}) with $\gamma = 2\pi\e$ and with
the perturbing functions $f = \e^2\tilde f$ and $g = \e^2\tilde g$ of
order $O(\e^2)$.

We are now in a position to apply Moser's theorem to obtain:

\begin{theorem}
Fix an annulus $\cal A$ of the form $c \le G \le d$ with $0<c<d$.
Then for $\epsilon$ sufficiently small, the Poincar\'e map
$P_\epsilon$ of the Nos\'e-Hoover system has infinitely many
invariant curves close to the level curves of $G$ in $\cal A$.
In particular, the corresponding flow is not ergodic.
\end{theorem}

The proof consists in checking the remaining  hypotheses of Moser's
theorem.  Note
that the differential equation (\ref{eq_nhav}) is real-analytic in
$(\theta,\tau,\alpha,\e)$.  It follows that the Poincar\'e map
$P_\e(\tau,\alpha)$ is real-analytic in $(\tau,\alpha,\e)$.   The
location of the fixed point $(\tau_\e,\alpha_\e)$ is a
real-analytic function of $\e$ (for $\e$ sufficiently small) and
so composing the Poincar\'e map with a translation to move this
point to $(1,0)$ preserves analyticity.  The action angle
variables $(\phi, G)$ are independent of $\e$ and are
real-analytic in $(\tau,\alpha)$ away from the fixed point itself.
Finally, the period function $T(G)$ is analytic in $G$.  Thus the
map (\ref{eq_poincmap2}) is real-analytic with respect to
$(x,y,\e)$ for $0<y<Y$, $|\e|<\e_0$ where $Y>0, \e_0>0$ are
constants. Since $x$ is an angular variable, the map is also periodic
with respect to $x$ with period $2\pi$.

The annulus $\cal A$ contains an annulus of the form $0<a\le y\le
b < Y$.  If $b-a = k \ne 1,$ then we can set $y= y'/k$ where $y'$
is a new variable. The map will retain the same form except that
$\gamma=2\pi\e$ is replaced by $\gamma' = 2\pi\e/k$. Once $a$ and
$b$ are fixed, the Poincar\'e map admits a complex-analytic
extension to some domain of the form $E = D\times \{|\e|\le
\e_1\}\subset \C^3$ where $D$ is a closed $\delta$-neighborhood of
$\R \times [a,b]$ in $\C^2$, and $\delta>0, \epsilon_1>0$ are
constants. For any fixed choice of $\e$ with $|\e|\le \e_1$, the
perturbing functions $\tilde f$ and $\tilde g$ in
(\ref{eq_poincmap2}) will have complex analytic extensions to $D$.
We now fix the Diophantine constants $c_0$ and $\mu$ in
(\ref{eq_Diophantine}) and choose any $\tilde \e>0$.  Letting
$\delta>0$ be the constant in (\ref{eq_delta}) guaranteed by
Moser's theorem, the theorem shows that for any $y_0\in [a,b]$
satisfying (\ref{eq_Diophantine}), and for any real analytic
functions $f, g$ with complex extensions to $D$ which satisfy
(\ref{eq_delta}), the corresponding map admits an invariant curve
$\tilde \e$-close to $y=y_0$.

We now let $f = \e^2\tilde f$ and $g = \e^2\tilde g$ where we think
of $\e$ as fixed, and we let $K = \sup_E |\tilde f|+|\tilde g|$.   Then
$|f|+|g|\le \e^2\,K$ on $D$.  Since $\gamma = 2\pi\e$ (or $\gamma'
= 2\pi\e/k$), it follows that (\ref{eq_delta}) holds for all $\e$
sufficiently small.

To complete the proof, it only remains to verify that the period
satisfies $T'(G) >0$.  We first introduce a new variable $\sigma =
\ln\tau$. Then (\ref{eq_unpert}) becomes:
\begin{equation}
\label{eq_unpert_ham}
\begin{split}
\s^\prime &= -\a,\\
\a^\prime &= e^\s -1,
\end{split}
\end{equation}
which is a planar Hamiltonian system with Hamiltonian function
$$G(\s,\a) = \smfr12 \a^2 + e^\s -1-\s.$$
In fact, it is of the classical kinetic plus potential form with potential
$$V(\s) = e^\s-1-\s$$
(except that in the ODE (\ref{eq_unpert_ham}), $\s$ plays the role of the momentum
variable and $\a$ the role of the position). The equilibrium point
is now at $(\s,\a)= (0,0)$.

The behavior of the period as a function of energy for such
systems is a well-studied problem. A result of Chicone
\cite{chicone} shows that $T(G)$ will be a strictly increasing
function of energy $G$ provided that $\dis{
\frac{V(\s)}{(V'(\s))^2} }$ is a strictly convex function (except
at $\s = 0$). This condition reads
$$6 V V''^2 - 3V'^2V'' - 2V V' V'' > 0,$$
except at $\s=0$.  It is not hard to check that this is true for
the potential $V(\s)$ above.

\begin{figure}
\includegraphics[width=14cm]{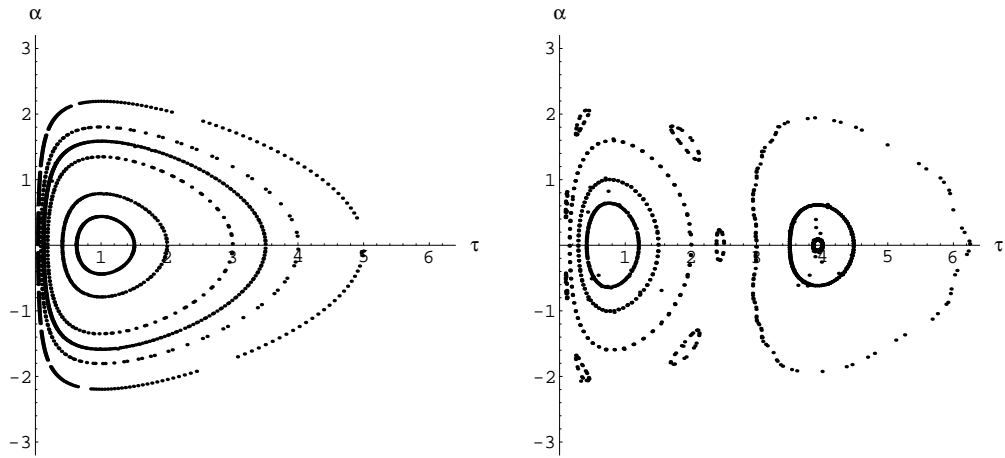}
\caption{Numerically computed orbits of the Poincar\'e map of the
plane $\theta = 0\bmod 2\pi$ for the ODE (\ref{eq_nh}); $\e=0.1$ (left)
and $\e=1.0$ (right).}\label{fig_Poinc}
\end{figure}

Figure~\ref{fig_Poinc} shows numerically computed Poincar\'e maps of the
ODE (\ref{eq_nh})
for two values of $\e$.  The Poincar\'e map for $\e=0.1$  has
invariant curves close to the level curves of $G$ in
Figure~\ref{fig_Gcontour}. Apparently the Poincar\'e map for
$\e=1$ still has many invariant curves although they are not
particularly close to the level curves of $G$.  If these really do
exist, then the Nos\'e-Hoover system is non-ergodic even for
$\e=1$.

It is interesting to remark that, for $\e= 1$, the Poincar\'e map
invariant curves are sometimes composed of a set of islands
(instead of being simple closed curves). This is the case for
instance with the initial condition $(\theta,\tau,\alpha) =
(0,2.42,0)$ which corresponds to $(q,p,\xi) = (2.2,0,0)$ (7
islands can be seen on the right hand side of
Figure~\ref{fig_Poinc}). Starting from this initial condition, we
numerically integrate (\ref{nh_original}) with a time step $\Delta
t = 0.001$ for $5.10^7$ time steps with the algorithm proposed in
\cite{molphys96}. This second-order algorithm is based on an
operator splitting technique. It preserves the measure
(\ref{measure_NH}) as well as the time-reversibility of
(\ref{nh_original}). Figure~\ref{fig_phase_space_NH1} shows the
projection of the time trajectory of the $(q,p)$-plane. With this
initial condition, the trajectory seems not even to sample a ring
of the form $\left\{(q,p) \in \R^2: c \leq q^2 + p^2 \leq C
\right\}$ for some $0 < c \leq C$, but only a part of it.

\begin{figure}
\includegraphics[width=8cm]{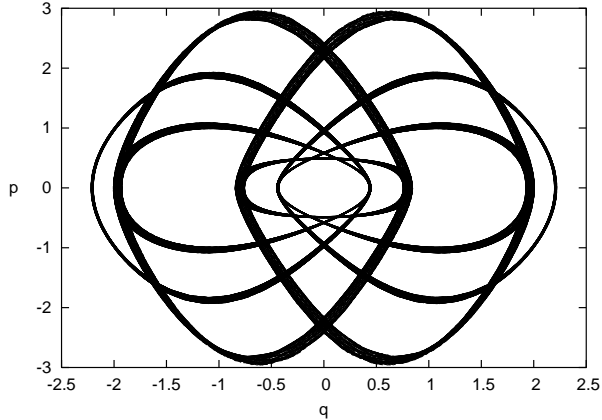}
\caption{Projection on the $(q,p)$-plane of the numerically computed
  trajectory of (\ref{nh_original}) with $\e=1$,
  starting from the initial condition $(q,p,\xi) = (2.2,0,0).$}
\label{fig_phase_space_NH1}
\end{figure}

\section{Harmonic oscillator coupled with a
  Nos\'e-Hoover chain}
\label{sec_2therm}

In the previous section, we have proven that a one-dimensional
harmonic oscillator coupled to a Nos\'e-Hoover thermostat is not
an ergodic system with respect to the Gibbs measure, when the mass
$Q$ of the reservoir is large. As mentioned above, numerical
observations of this fact have already been reported in
\cite{Hoover,Martyna92,Tuckerman00}, even for moderate values of
$Q$ such as $Q=1$. As a consequence, it is known that one should
be cautious when making use of the Nos\'e-Hoover equations
(\ref{dyn_NH}) to compute phase space integrals such as
(\ref{ps_aver}).

One way that has been proposed to circumvent this difficulty is to
generalize the Nos\'e-Hoover equations to the so-called
Nos\'e-Hoover chain equations \cite{Martyna92}. The idea consists
in coupling the physical variables $(q,p)$ with a first thermostat
as in (\ref{dyn_NH}) and to then couple this thermostat with a
second one, which can be coupled to a third one, and so on. The
variables now include $M_{\rm ext}$ additional scalar variables,
$\xi_j$ for $j=1, \ldots, M_{\rm ext}$, where the number $M_{\rm
ext}$ of thermostats can be freely specified. The Nos\'e-Hoover
chain dynamics is given by
\begin{equation}
\label{dyn_NHC}
\begin{split}
\frac{d q_i}{d t} &= \frac{p_i}{m_i},
\\
\frac{d p_i}{d t} &= -\nabla_{q_i} V - \frac{\xi_1}{Q_1} \, p_i,
\\
\frac{d \xi_1}{d t} &= \left( \sum_{i=1}^M \frac{p_i^2}{m_i}  - n
M
  \beta^{-1} \right) - \frac{\xi_2}{Q_2} \, \xi_1,
\\
\frac{d \xi_j}{d t} &= \left( \frac{\xi_{j-1}^2}{Q_{j-1}}  -
\beta^{-1}
  \right) - \frac{\xi_{j+1}}{Q_{j+1}} \, \xi_{j}
   \quad \text{for }j=2, \ldots, M_{\rm ext}-1,
\\
\frac{d \xi_{M_{\rm ext}}}{d t} &= \frac{\xi_{M_{\rm
      ext}-1}^2}{Q_{M_{\rm ext}-1}}  - \beta^{-1},
\end{split}
\end{equation}
where the masses $Q_1, \ldots, Q_{M_{\rm ext}}$ are free
parameters that can be arbitrarily specified.

For the Nos\'e-Hoover system \eqref{dyn_NHC} with
$z=(q,p,\xi_1,\dots,\xi_{M_{\rm ext}}),$ it is easy to check that
\[
\rho(q,p,\xi_1,\dots,\xi_{M_{\rm ext}})=\exp \left(-\beta \left[
H(q,p) + \sum_{j=1}^{M_{\rm ext}} \frac{\xi_j^2}{2Q_j} \right]
\right)
\]
satisfies the condition \eqref{inv}, so an invariant measure for
the Nos\'e-Hoover system \eqref{dyn_NHC} is given by
\begin{equation}
\label{measure_NHC} d\mu_{\rm NHC}(q,p,\xi_1,\dots,\xi_{M_{\rm
ext}}) = \frac{\exp \left(-\beta \left[ H(q,p) + \sum_{j=1}^{M_{\rm
ext}}\frac{\xi_j^2}{2Q_j}
 \right] \right) \, dq \, dp \, d \xi_1 \ldots d \xi_{M_{\rm ext}}}
 {\displaystyle\int\exp \left(-\beta \left[ H(q,p) + \sum_{j=1}^{M_{\rm
ext}}\frac{\xi_j^2}{2Q_j}
 \right] \right) \, dq \, dp \, d \xi_1 \ldots d \xi_{M_{\rm ext}}}.
\end{equation}

We can now calculate that
\[
{\dis \int \sum_{i=1}^M \frac{p_i^2}{m_i} \,d \mu_{\rm
NHC}}= nM\beta^{-1},
\]
and
\[
{\dis \int \frac{\xi_j^2}{Q_j} \,d \mu_{\rm NHC}} = \beta^{-1}\quad\text{ for } j=1,
\ldots, M_{\rm ext}.
\]
Hence, we can see that the evolution of $\xi_1$ in (\ref{dyn_NHC})
is controlled by the difference between the instantaneous value of
twice the kinetic energy $\sum_{i=1}^M \frac{p_i^2}{m_i}$
 and its average value with respect to the
invariant measure $d\mu_{\rm NHC}$ and the evolution of $\xi_j$
for $j=2, \ldots, M_{\rm ext}$ in (\ref{dyn_NHC}) is controlled by
the difference between the instantaneous value of twice the
``kinetic energy'' $\xi_{j-1}^2/Q_{j-1}$ and its average value
with respect to the invariant measure $d\mu_{\rm NHC}.$

If (\ref{dyn_NHC}) is ergodic with respect to $d \mu_{\rm
  NHC}$, computations similar to the ones performed in Section
\ref{Sec_prez} show that
\begin{equation}
\label{ergo_NHC} \lim_{T \to +\infty} \frac1 T \int_0^{T} A \left(
q(t), p(t)\right) dt = \int A(q,p) d \mu(q,p),
\end{equation}
where $(q(t),p(t),\xi_1(t),\dots,\xi_{M_{\rm ext}}(t))$ is a
solution of the Nos\'e-Hoover chain dynamics and $d\mu(q,p)$ is
the Gibbs measure (\ref{measure_can}). Hence, the Nos\'e-Hoover
chain dynamics can be used to compute phase space integrals in the
canonical ensemble if (\ref{dyn_NHC}) is ergodic with respect to
$d \mu_{\rm
  NHC}.$

In this section, we numerically study the case of a
one-dimensional harmonic oscillator coupled to a Nos\'e-Hoover
chain of length $M_{\rm
  ext} = 2$. We consider the case $Q_1 = Q_2$.  In
view of (\ref{dyn_NHC}), the dynamics reads
\begin{equation}
\label{nhc_original}
\begin{split}
\dot q &= p,\\
\dot p &= -q - \e^2 p \xi_1,\\
\dot \xi_1 &= p^2 - 1 - \e^2 \xi_1 \xi_2, \\
\dot \xi_2 &= \e^2 \xi_1^2 - 1,
\end{split}
\end{equation}
where $\e = 1/\sqrt{Q_1} = 1/\sqrt{Q_2}$. We again introduce the
action-angle variables $(\tau,\theta)$ for the oscillator (see
(\ref{action_angle})) and rescale via $\a_j = \e \xi_j$ to get
\begin{equation}
\label{eq_nhc}
\begin{split}
\dot \theta &= 1 - \e \a_1 \sin\theta \cos\theta,\\
\dot \tau &= -2\e\tau\a_1 \sin^2\theta,\\
\dot \a_1 &= \e (2\tau\sin^2\theta - 1 - \a_1 \a_2), \\
\dot \a_2 &= \e (\a_1^2 - 1).
\end{split}
\end{equation}
For small $\e$, the corresponding averaged system reads after
rescaling time,
\begin{equation}
\label{eq_nhc_unpert}
\begin{split}
\dot \tau &= -\tau\a_1,\\
\dot \a_1 &= \tau - 1 - \a_1 \a_2, \\
\dot \a_2 &= \a_1^2 - 1.
\end{split}
\end{equation}

In the case of the Nos\'e-Hoover equation, the averaged system is
(\ref{eq_unpert}), and the analysis conducted in
Section~\ref{Sec_analysis} is based on the knowledge of a first
integral for (\ref{eq_unpert}), namely the function $G(\tau,\a)$
defined by (\ref{defG}). In the case of (\ref{eq_nhc_unpert}), we
were not able to find such a first integral. Let us follow a
different route and study the system by numerically computing the
Poincar\'e return map of the plane $\Sigma_{\rm NHC}^{\rm av} =
\{(\tau,\a_1,\a_2): \a_2 = 0 \}$ for the ODE
(\ref{eq_nhc_unpert}).  Two trajectories of the return map with
different initial conditions are shown in
Figure~\ref{fig_Poinc_NHC_av}.  Some initial conditions lead to
trajectories which lie on invariant curves (see right-hand side of
Figure~\ref{fig_Poinc_NHC_av}) and the corresponding values of
$\tau$ that are sampled are bounded from below and isolated from
0: there exists $\tau_c^{\e = 0} > 0$ such that for all $t$ we
have $\tau(t) \geq \tau_c^{\e = 0}.$  For other initial
conditions, the trajectory does not seem to be confined to a curve
(see left-hand side of Figure~\ref{fig_Poinc_NHC_av}), but it
still does not sample the entire plane.

\begin{figure}[h]
\centerline{\input{fig_P_return_map1.tex}
\input{fig_P_return_map2.tex}}
\caption{Numerically computed orbits of the Poincar\'e return map
of the plane $\a_2 = 0$ for the ODE (\ref{eq_nhc_unpert}). Left:
the initial condition is $(\tau,\a_1,\a_2) = (q_0^2/2,0,0)$ with
$q_0 = 0.5$; a similar trajectory is obtained for $q_0 = 1.3$ and
$q_0 = 1.5$. Right: the initial condition is $(\tau,\a_1,\a_2) =
(q_0^2/2,0,0)$ with $q_0 = 1.1$; a similar trajectory is obtained
for $q_0$ = 0.65, 0.7, 0.75, 0.8, 0.85, 0.9, 0.95 and 1.0.}
\label{fig_Poinc_NHC_av}
\end{figure}
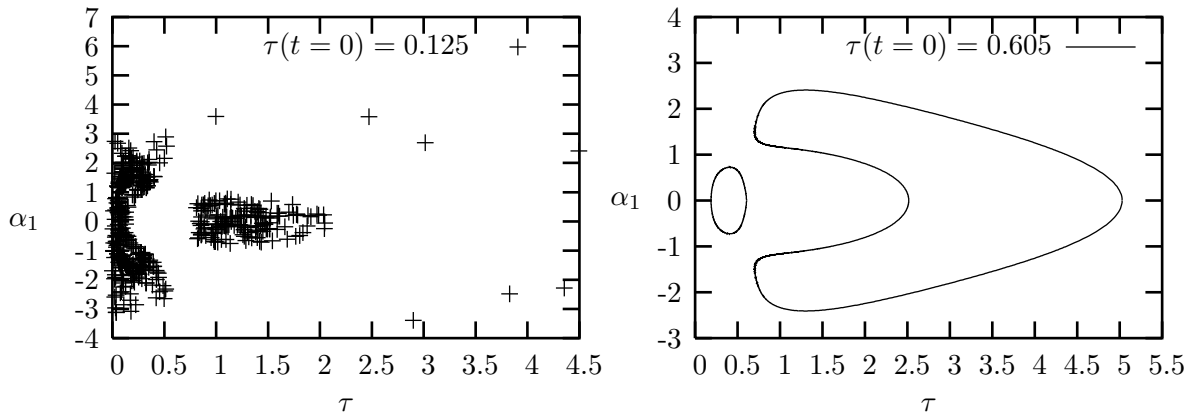

We now consider the ODE (\ref{eq_nhc}) and study whether the
behaviour we observe in the case $\e \to 0$ still persists in the
case of a small but positive $\e$. We choose $Q_1 = Q_2 = 10$,
that is $\e =1/\sqrt{10} = 0.316$, and we numerically integrate
(\ref{nhc_original}) with the second-order algorithm proposed in
\cite{molphys96} from the initial condition $(q,p,\xi_1,\x_2) =
(1.1,0,0,0)$. This condition corresponds to an initial condition
for which the Poincar\'e return map of (\ref{eq_nhc_unpert}) seems
to have invariant curves (see Figure~\ref{fig_Poinc_NHC_av},
right-hand side). Figure~\ref{fig_Poinc_NHC} shows the trace of
the time trajectory of (\ref{eq_nhc}) on the plane $\Sigma_{\rm
NHC} = \{(\theta,\tau,\a_1,\a_2): \a_2 = 0 \}$. We can see the
topological similarity between curves on
Figures~\ref{fig_Poinc_NHC_av} (right-hand side) and
\ref{fig_Poinc_NHC}, although they are quantitatively quite
different. Figure~\ref{fig_phase_space_NH2} shows the projection
of the same trajectory on the $(q,p)$-plane. We see that the
values of $(q,p)$ that are sampled still satisfy $\tau(t) =
(q^2(t)+p^2(t))/2 \geq \tau_c^{\e}$ for some $\tau_c^{\e} > 0$ and
all $t \geq 0$. Thus, the lower bound on the values of $\tau$ that
we observe for the averaged dynamics persists in the case $\e =
1/\sqrt{10}$. For the initial condition we consider here,
$\tau_c^{\e = 0} = 0.188$ whereas $\tau_c^{\e} = 0.194$. Because
of this lower bound, we obtain a contradiction with
(\ref{ergo_NHC}), and hence this numerical experiment indicates
that a one-dimensional harmonic oscillator coupled to a
Nos\'e-Hoover chain of two thermostats is not always an ergodic
system. We observe this non-ergodicity for many different initial
conditions (including all those listed in
Figure~\ref{fig_Poinc_NHC_av}), and for many different values of
$Q \geq 10$.

\begin{figure}[h]
\centerline{\input{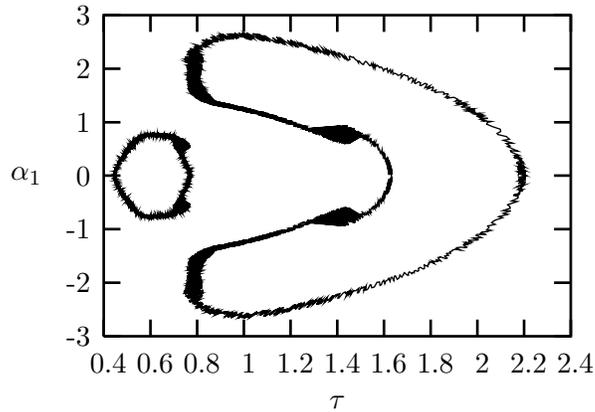}} \caption{Trace of the
trajectory of (\ref{eq_nhc}) with $\e =
  1/\sqrt{10}$ on the plane $\a_2 = 0$. The initial condition is
  $(\theta,\tau,\a_1,\a_2) = (0,q_0^2/2,0,0)$ with $q_0 = 1.1$.}
\label{fig_Poinc_NHC}
\end{figure}

\begin{figure}
\includegraphics[width=6cm]{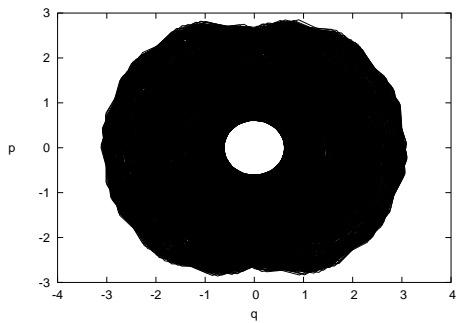}
\caption{Projection on the $(q,p)$-plane of the
  trajectory of (\ref{eq_nhc}) with $\e=1/\sqrt{10}$,
  starting from the initial condition $(\theta,\tau,\a_1,\a_2) =
  (0,q_0^2/2,0,0)$ with $q_0 = 1.1$.}
\label{fig_phase_space_NH2}
\end{figure}

In the case $Q = 1$, we did not find any initial condition such
that the values of $\tau(t)$ are isolated from 0. Results obtained
with $Q=1$ and the same initial condition as previously, namely
$(q,p,\xi_1,\x_2) = (1.1,0,0,0)$, are shown on
Figure~\ref{fig_2therm_Q1_amplitude}. We compare the theoretical
distributions of the angular variable $\theta$ and of the
amplitude variable $r = \sqrt{q^2+p^2}$ (as given by the Gibbs
measure (\ref{measure_can})) with the empirical distributions
obtained from the time trajectory. Since we work with $\beta = 1$,
the theoretical distributions respectively read
$$
f_{\rm theo}^{\rm ang}(\theta) = \frac{1}{2 \pi} 1_{[0,2\pi]} \
d\theta, \quad f_{\rm theo}^{\rm amp}(r) = r \exp \left( -
\frac{r^2}2 \right) 1_{[0,+\infty)} \ dr,
$$
where $1_{[0,2\pi]}$ is the characteristic function of $[0,2\pi]$.
The numerical distributions are denoted $f_{\rm num}^{\rm
ang}(\theta)$ and $f_{\rm num}^{\rm amp}(r)$. They have been
computed from a trajectory of length $T = 2.5 \ 10^6$, with a time
step of $2.5 \ 10^{-3}$ (the trajectory is thus composed of $10^9$
time steps), by partitioning the sampled interval into 100 bins (a
cutoff of $r_c=4$ has been used). Note that all distributions have
been normalized so that their integrals are equal to 1.

For the angular variable, we plot on
Figure~\ref{fig_2therm_Q1_amplitude} both distributions $f_{\rm
  theo}^{\rm ang}(\theta)$ and $f_{\rm num}^{\rm ang}(\theta)$. We can
see some small oscillations of the empirical distribution around
the theoretical uniform distribution. Regarding the amplitude
variable, it is more convenient to directly plot the difference
$\left| f_{\rm theo}^{\rm amp}(r) - f_{\rm num}^{\rm amp}(r)
\right|$. Again, we can see some oscillations, which are small
(let us recall that the maximum value of $f_{\rm theo}^{\rm
amp}(r)$ is $f_{\rm theo}^{\rm amp}(1) = 0.6$). We have checked
that the error does not change when the time step is reduced, the
total length $T$ of the trajectory being kept fixed.

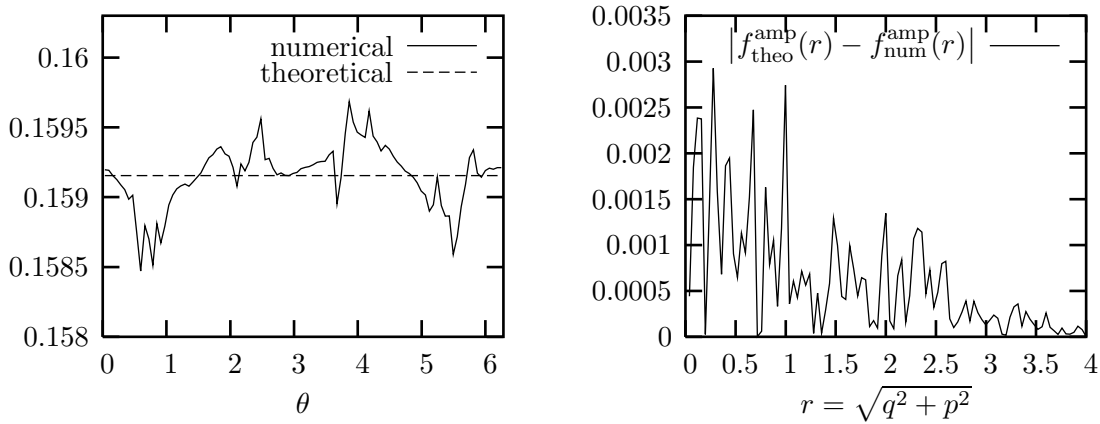
\begin{figure}[h]
\centerline{\input{fig_2therm_Q1_angle.tex}
\input{fig_2therm_Q1_amplitude.tex}}
\caption{Left: Numerical and theoretical distributions of
$\theta$;
  right: difference between the numerical and theoretical distributions
 of $r = \sqrt{q^2+p^2}$. The numerical distribution is obtained from the
  simulation of (\ref{eq_nhc}) with $\e=1$, starting from the initial
  condition $(q,p,\xi_1,\x_2) = (1.1,0,0,0)$.}
\label{fig_2therm_Q1_amplitude}
\end{figure}

To study the evolution of the error with $T$ (the time step being
now kept fixed), the indicator we consider is the star
discrepancy. Recall that the star discrepancy $D_N$ of a sequence
$x = \{ x_n \}_{1 \leq n \leq N}$ with values in $[0,1]^2$ is
defined as \cite{LPS}
\begin{equation}
  \label{discrepancy_definition}
  D_N^*(x) = \sup_{y \in [0,1]^2 } \left | \frac{1}{N} \sum_{n=1}^N
  1_{ [0,y] }(x_n)
    - \textrm{Volume}([0,y]) \right |,
\end{equation}
where, for $2$-dimensional vectors $y$ and $z,$ we write $y \leq
z$ when $y_i \leq z_i$ for all $1 \leq i \leq 2$, and note that
$[0,y] = \{ z \in [0,1]^2, \ z \leq y \}$. The fact that $D_N^*(x)
\to 0$ when $N \to \infty$ is equivalent \cite[p.15]{LPS} to the
fact that, for any Riemann integrable function $A$ defined on
$[0,1]^2$,
\[
\lim_{N \to \infty} \frac{1}{N} \sum_{n=1}^N A(x_n) = \int_{
[0,1]^d } A(x) \, d x.
\]
In addition, for functions $A$ which have bounded variations
$V_{\rm
  HK}(A)$ in the sense of Hardy and Krause\cite{Nei92}, the following
error estimate holds true:
\begin{equation}
\label{estim_discr} \left| \frac1N \sum_{n=1}^N A(x_n) - \int_{
[0,1]^2 } A(x) \, d x \right| \leq V_{\rm HK}(A) D_N^*(x).
\end{equation}
If $A \in C^2([0,1]^2)$, then its variation $V_{\rm HK}(A)$ has a
simple expression:
\begin{equation*}
V_{\rm HK}(A) =  \int_{ [0,1]^2 } \left| \frac{\partial^2 A}{\partial x_{1}
    \partial x_{2}} \right| \ dx
+\int_0^1 \left| \frac{\partial A}{\partial x_1}(x_1,1)\right|\,dx_1
+\int_0^1 \left| \frac{\partial A}{\partial x_2}(1,x_2)\right|\,dx_2.
\end{equation*}

In view of (\ref{estim_discr}), we can see that the convergence of
$D_N^*(x)$ toward 0 implies the convergence of the empirical average
of $A$ towards its statistical average, and the rate of
convergence of $D_N^*(x)$ gives information about the convergence
rate of the observable average.

Here, we work with the discrepancy criterion
\begin{equation}
\label{disc_continuous} D_N^* \left( \{ \theta_n,r_n \} \right) =
\sup_{(\theta,r) \in [0,2 \pi]
  \times [0,r_c] } \left | \frac{1}{N} \sum_{n=1}^N
  1_{ [0,\theta] } (\theta_n) \, 1_{[0,r]} (r_n)
    - \int_{\bar \theta = 0}^\theta \int_{\bar r = 0}^r \frac{\bar r}{2 \pi}
  \exp \left( - \frac{\bar r^2}2 \right) d \bar \theta \, d \bar r \right|,
\end{equation}
where the sample $(\theta_n,r_n)$ is obtained by the numerical
integration of the ODE with a time step of $\Delta t = 2.5 \
10^{-3}$. We again set the cutoff amplitude at $r_c = 4$. Note
that the integral that appears in (\ref{disc_continuous}) can be
exactly computed. In practice, $D_N^*$ is approximated by
considering the supremum only over $(\theta,r)$ of the form
$(\theta_k,r_l)$ with $\theta_k = 2 k \pi / 100$, $r_l = l
r_c/100$, and $1 \leq k, l \leq 100$.

We have considered 8 different initial conditions, giving rise to
8 different samples $(\theta_n,r_n)$, and for each of them, we
have computed the star discrepancy $D_N^*$ for several values of
$N$. Results are shown on Figure~\ref{fig_2therm_Q1_disc}, where
we plot the mean of the computed discrepancies as a function of
$N$. A least mean-square fit shows that the mean discrepancy
decreases as $\dis{ \frac{11.1}{N^{0.483}} }$.

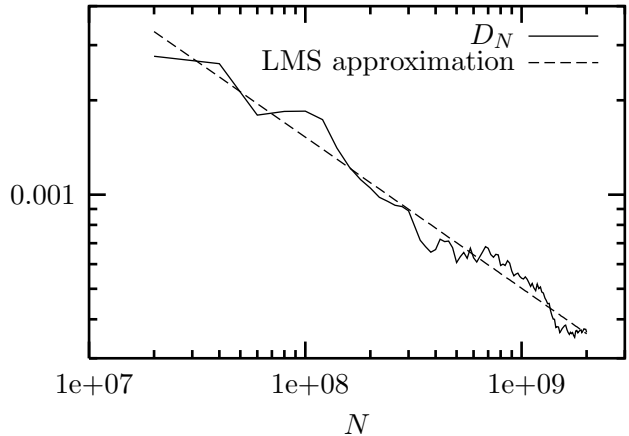
\begin{figure}[h]
\centerline{\input{fig_2therm_Q1_disc.tex}} \caption{Star
discrepancy (\ref{disc_continuous}) as a function of the sample
  size $N$, and comparison to the least mean-square fit $D_N^* \approx
  11.1/N^{0.483}$.}
\label{fig_2therm_Q1_disc}
\end{figure}

We thus see from our numerical experiments with these thermostat
masses and these initial conditions that the numerical
distribution converges to the theoretical distribution, and the
numerical results are thus consistent with ergodicity. However,
more extensive numerical testing would be needed to assess whether
this would be the case for any initial condition or any smaller
thermostat masses.

\section{Acknowledgments}
We thank Alain Chenciner and Jacques Fejoz for helpful
suggestions on this research.

}
\end{document}

%% file: fig_P_return_map1.tex
\begingroup
  \catcode`\@=11\relax
  \def\GNUPLOTspecial{%
    \def\do##1{\catcode`##1=12\relax}\dospecials
    \catcode`\{=1\catcode`\}=2\catcode\%=14\relax\special}%
\expandafter\ifx\csname GNUPLOTpicture\endcsname\relax
  \csname newdimen\endcsname\GNUPLOTunit
  \gdef\GNUPLOTpicture(#1,#2){\vbox to#2\GNUPLOTunit\bgroup
    \def\put(##1,##2)##3{\unskip\raise##2\GNUPLOTunit
      \hbox to0pt{\kern##1\GNUPLOTunit ##3\hss}\ignorespaces}%
    \def\ljust##1{\vbox to0pt{\vss\hbox to0pt{##1\hss}\vss}}%
    \def\cjust##1{\vbox to0pt{\vss\hbox to0pt{\hss ##1\hss}\vss}}%
    \def\rjust##1{\vbox to0pt{\vss\hbox to0pt{\hss ##1}\vss}}%
    \def\stack##1{\let\\=\cr\tabskip=0pt\halign{\hfil ####\hfil\cr ##1\crcr}}%
    \def\lstack##1{\hbox to0pt{\vbox to0pt{\vss\stack{##1}}\hss}}%
    \def\cstack##1{\hbox to0pt{\hss\vbox to0pt{\vss\stack{##1}}\hss}}%
    \def\rstack##1{\hbox to0pt{\vbox to0pt{\stack{##1}\vss}\hss}}%
    \vss\hbox to#1\GNUPLOTunit\bgroup\ignorespaces}%
  \gdef\endGNUPLOTpicture{\hss\egroup\egroup}%
\fi
\GNUPLOTunit=0.1bp
{\GNUPLOTspecial{!
/gnudict 256 dict def
gnudict begin
/Color false def
/Solid false def
/gnulinewidth 5.000 def
/userlinewidth gnulinewidth def
/vshift -33 def
/dl {10 mul} def
/hpt_ 31.5 def
/vpt_ 31.5 def
/hpt hpt_ def
/vpt vpt_ def
/M {moveto} bind def
/L {lineto} bind def
/R {rmoveto} bind def
/V {rlineto} bind def
/vpt2 vpt 2 mul def
/hpt2 hpt 2 mul def
/Lshow { currentpoint stroke M
  0 vshift R show } def
/Rshow { currentpoint stroke M
  dup stringwidth pop neg vshift R show } def
/Cshow { currentpoint stroke M
  dup stringwidth pop -2 div vshift R show } def
/UP { dup vpt_ mul /vpt exch def hpt_ mul /hpt exch def
  /hpt2 hpt 2 mul def /vpt2 vpt 2 mul def } def
/DL { Color {setrgbcolor Solid {pop []} if 0 setdash }
 {pop pop pop Solid {pop []} if 0 setdash} ifelse } def
/BL { stroke userlinewidth 2 mul setlinewidth } def
/AL { stroke userlinewidth 2 div setlinewidth } def
/UL { dup gnulinewidth mul /userlinewidth exch def
      dup 1 lt {pop 1} if 10 mul /udl exch def } def
/PL { stroke userlinewidth setlinewidth } def
/LTb { BL [] 0 0 0 DL } def
/LTa { AL [1 udl mul 2 udl mul] 0 setdash 0 0 0 setrgbcolor } def
/LT0 { PL [] 1 0 0 DL } def
/LT1 { PL [4 dl 2 dl] 0 1 0 DL } def
/LT2 { PL [2 dl 3 dl] 0 0 1 DL } def
/LT3 { PL [1 dl 1.5 dl] 1 0 1 DL } def
/LT4 { PL [5 dl 2 dl 1 dl 2 dl] 0 1 1 DL } def
/LT5 { PL [4 dl 3 dl 1 dl 3 dl] 1 1 0 DL } def
/LT6 { PL [2 dl 2 dl 2 dl 4 dl] 0 0 0 DL } def
/LT7 { PL [2 dl 2 dl 2 dl 2 dl 2 dl 4 dl] 1 0.3 0 DL } def
/LT8 { PL [2 dl 2 dl 2 dl 2 dl 2 dl 2 dl 2 dl 4 dl] 0.5 0.5 0.5 DL } def
/Pnt { stroke [] 0 setdash
   gsave 1 setlinecap M 0 0 V stroke grestore } def
/Dia { stroke [] 0 setdash 2 copy vpt add M
  hpt neg vpt neg V hpt vpt neg V
  hpt vpt V hpt neg vpt V closepath stroke
  Pnt } def
/Pls { stroke [] 0 setdash vpt sub M 0 vpt2 V
  currentpoint stroke M
  hpt neg vpt neg R hpt2 0 V stroke
  } def
/Box { stroke [] 0 setdash 2 copy exch hpt sub exch vpt add M
  0 vpt2 neg V hpt2 0 V 0 vpt2 V
  hpt2 neg 0 V closepath stroke
  Pnt } def
/Crs { stroke [] 0 setdash exch hpt sub exch vpt add M
  hpt2 vpt2 neg V currentpoint stroke M
  hpt2 neg 0 R hpt2 vpt2 V stroke } def
/TriU { stroke [] 0 setdash 2 copy vpt 1.12 mul add M
  hpt neg vpt -1.62 mul V
  hpt 2 mul 0 V
  hpt neg vpt 1.62 mul V closepath stroke
  Pnt  } def
/Star { 2 copy Pls Crs } def
/BoxF { stroke [] 0 setdash exch hpt sub exch vpt add M
  0 vpt2 neg V  hpt2 0 V  0 vpt2 V
  hpt2 neg 0 V  closepath fill } def
/TriUF { stroke [] 0 setdash vpt 1.12 mul add M
  hpt neg vpt -1.62 mul V
  hpt 2 mul 0 V
  hpt neg vpt 1.62 mul V closepath fill } def
/TriD { stroke [] 0 setdash 2 copy vpt 1.12 mul sub M
  hpt neg vpt 1.62 mul V
  hpt 2 mul 0 V
  hpt neg vpt -1.62 mul V closepath stroke
  Pnt  } def
/TriDF { stroke [] 0 setdash vpt 1.12 mul sub M
  hpt neg vpt 1.62 mul V
  hpt 2 mul 0 V
  hpt neg vpt -1.62 mul V closepath fill} def
/DiaF { stroke [] 0 setdash vpt add M
  hpt neg vpt neg V hpt vpt neg V
  hpt vpt V hpt neg vpt V closepath fill } def
/Pent { stroke [] 0 setdash 2 copy gsave
  translate 0 hpt M 4 {72 rotate 0 hpt L} repeat
  closepath stroke grestore Pnt } def
/PentF { stroke [] 0 setdash gsave
  translate 0 hpt M 4 {72 rotate 0 hpt L} repeat
  closepath fill grestore } def
/Circle { stroke [] 0 setdash 2 copy
  hpt 0 360 arc stroke Pnt } def
/CircleF { stroke [] 0 setdash hpt 0 360 arc fill } def
/C0 { BL [] 0 setdash 2 copy moveto vpt 90 450  arc } bind def
/C1 { BL [] 0 setdash 2 copy        moveto
       2 copy  vpt 0 90 arc closepath fill
               vpt 0 360 arc closepath } bind def
/C2 { BL [] 0 setdash 2 copy moveto
       2 copy  vpt 90 180 arc closepath fill
               vpt 0 360 arc closepath } bind def
/C3 { BL [] 0 setdash 2 copy moveto
       2 copy  vpt 0 180 arc closepath fill
               vpt 0 360 arc closepath } bind def
/C4 { BL [] 0 setdash 2 copy moveto
       2 copy  vpt 180 270 arc closepath fill
               vpt 0 360 arc closepath } bind def
/C5 { BL [] 0 setdash 2 copy moveto
       2 copy  vpt 0 90 arc
       2 copy moveto
       2 copy  vpt 180 270 arc closepath fill
               vpt 0 360 arc } bind def
/C6 { BL [] 0 setdash 2 copy moveto
      2 copy  vpt 90 270 arc closepath fill
              vpt 0 360 arc closepath } bind def
/C7 { BL [] 0 setdash 2 copy moveto
      2 copy  vpt 0 270 arc closepath fill
              vpt 0 360 arc closepath } bind def
/C8 { BL [] 0 setdash 2 copy moveto
      2 copy vpt 270 360 arc closepath fill
              vpt 0 360 arc closepath } bind def
/C9 { BL [] 0 setdash 2 copy moveto
      2 copy  vpt 270 450 arc closepath fill
              vpt 0 360 arc closepath } bind def
/C10 { BL [] 0 setdash 2 copy 2 copy moveto vpt 270 360 arc closepath fill
       2 copy moveto
       2 copy vpt 90 180 arc closepath fill
               vpt 0 360 arc closepath } bind def
/C11 { BL [] 0 setdash 2 copy moveto
       2 copy  vpt 0 180 arc closepath fill
       2 copy moveto
       2 copy  vpt 270 360 arc closepath fill
               vpt 0 360 arc closepath } bind def
/C12 { BL [] 0 setdash 2 copy moveto
       2 copy  vpt 180 360 arc closepath fill
               vpt 0 360 arc closepath } bind def
/C13 { BL [] 0 setdash  2 copy moveto
       2 copy  vpt 0 90 arc closepath fill
       2 copy moveto
       2 copy  vpt 180 360 arc closepath fill
               vpt 0 360 arc closepath } bind def
/C14 { BL [] 0 setdash 2 copy moveto
       2 copy  vpt 90 360 arc closepath fill
               vpt 0 360 arc } bind def
/C15 { BL [] 0 setdash 2 copy vpt 0 360 arc closepath fill
               vpt 0 360 arc closepath } bind def
/Rec   { newpath 4 2 roll moveto 1 index 0 rlineto 0 exch rlineto
       neg 0 rlineto closepath } bind def
/Square { dup Rec } bind def
/Bsquare { vpt sub exch vpt sub exch vpt2 Square } bind def
/S0 { BL [] 0 setdash 2 copy moveto 0 vpt rlineto BL Bsquare } bind def
/S1 { BL [] 0 setdash 2 copy vpt Square fill Bsquare } bind def
/S2 { BL [] 0 setdash 2 copy exch vpt sub exch vpt Square fill Bsquare } bind def
/S3 { BL [] 0 setdash 2 copy exch vpt sub exch vpt2 vpt Rec fill Bsquare } bind def
/S4 { BL [] 0 setdash 2 copy exch vpt sub exch vpt sub vpt Square fill Bsquare } bind def
/S5 { BL [] 0 setdash 2 copy 2 copy vpt Square fill
       exch vpt sub exch vpt sub vpt Square fill Bsquare } bind def
/S6 { BL [] 0 setdash 2 copy exch vpt sub exch vpt sub vpt vpt2 Rec fill Bsquare } bind def
/S7 { BL [] 0 setdash 2 copy exch vpt sub exch vpt sub vpt vpt2 Rec fill
       2 copy vpt Square fill
       Bsquare } bind def
/S8 { BL [] 0 setdash 2 copy vpt sub vpt Square fill Bsquare } bind def
/S9 { BL [] 0 setdash 2 copy vpt sub vpt vpt2 Rec fill Bsquare } bind def
/S10 { BL [] 0 setdash 2 copy vpt sub vpt Square fill 2 copy exch vpt sub exch vpt Square fill
       Bsquare } bind def
/S11 { BL [] 0 setdash 2 copy vpt sub vpt Square fill 2 copy exch vpt sub exch vpt2 vpt Rec fill
       Bsquare } bind def
/S12 { BL [] 0 setdash 2 copy exch vpt sub exch vpt sub vpt2 vpt Rec fill Bsquare } bind def
/S13 { BL [] 0 setdash 2 copy exch vpt sub exch vpt sub vpt2 vpt Rec fill
       2 copy vpt Square fill Bsquare } bind def
/S14 { BL [] 0 setdash 2 copy exch vpt sub exch vpt sub vpt2 vpt Rec fill
       2 copy exch vpt sub exch vpt Square fill Bsquare } bind def
/S15 { BL [] 0 setdash 2 copy Bsquare fill Bsquare } bind def
/D0 { gsave translate 45 rotate 0 0 S0 stroke grestore } bind def
/D1 { gsave translate 45 rotate 0 0 S1 stroke grestore } bind def
/D2 { gsave translate 45 rotate 0 0 S2 stroke grestore } bind def
/D3 { gsave translate 45 rotate 0 0 S3 stroke grestore } bind def
/D4 { gsave translate 45 rotate 0 0 S4 stroke grestore } bind def
/D5 { gsave translate 45 rotate 0 0 S5 stroke grestore } bind def
/D6 { gsave translate 45 rotate 0 0 S6 stroke grestore } bind def
/D7 { gsave translate 45 rotate 0 0 S7 stroke grestore } bind def
/D8 { gsave translate 45 rotate 0 0 S8 stroke grestore } bind def
/D9 { gsave translate 45 rotate 0 0 S9 stroke grestore } bind def
/D10 { gsave translate 45 rotate 0 0 S10 stroke grestore } bind def
/D11 { gsave translate 45 rotate 0 0 S11 stroke grestore } bind def
/D12 { gsave translate 45 rotate 0 0 S12 stroke grestore } bind def
/D13 { gsave translate 45 rotate 0 0 S13 stroke grestore } bind def
/D14 { gsave translate 45 rotate 0 0 S14 stroke grestore } bind def
/D15 { gsave translate 45 rotate 0 0 S15 stroke grestore } bind def
/DiaE { stroke [] 0 setdash vpt add M
  hpt neg vpt neg V hpt vpt neg V
  hpt vpt V hpt neg vpt V closepath stroke } def
/BoxE { stroke [] 0 setdash exch hpt sub exch vpt add M
  0 vpt2 neg V hpt2 0 V 0 vpt2 V
  hpt2 neg 0 V closepath stroke } def
/TriUE { stroke [] 0 setdash vpt 1.12 mul add M
  hpt neg vpt -1.62 mul V
  hpt 2 mul 0 V
  hpt neg vpt 1.62 mul V closepath stroke } def
/TriDE { stroke [] 0 setdash vpt 1.12 mul sub M
  hpt neg vpt 1.62 mul V
  hpt 2 mul 0 V
  hpt neg vpt -1.62 mul V closepath stroke } def
/PentE { stroke [] 0 setdash gsave
  translate 0 hpt M 4 {72 rotate 0 hpt L} repeat
  closepath stroke grestore } def
/CircE { stroke [] 0 setdash 
  hpt 0 360 arc stroke } def
/Opaque { gsave closepath 1 setgray fill grestore 0 setgray closepath } def
/DiaW { stroke [] 0 setdash vpt add M
  hpt neg vpt neg V hpt vpt neg V
  hpt vpt V hpt neg vpt V Opaque stroke } def
/BoxW { stroke [] 0 setdash exch hpt sub exch vpt add M
  0 vpt2 neg V hpt2 0 V 0 vpt2 V
  hpt2 neg 0 V Opaque stroke } def
/TriUW { stroke [] 0 setdash vpt 1.12 mul add M
  hpt neg vpt -1.62 mul V
  hpt 2 mul 0 V
  hpt neg vpt 1.62 mul V Opaque stroke } def
/TriDW { stroke [] 0 setdash vpt 1.12 mul sub M
  hpt neg vpt 1.62 mul V
  hpt 2 mul 0 V
  hpt neg vpt -1.62 mul V Opaque stroke } def
/PentW { stroke [] 0 setdash gsave
  translate 0 hpt M 4 {72 rotate 0 hpt L} repeat
  Opaque stroke grestore } def
/CircW { stroke [] 0 setdash 
  hpt 0 360 arc Opaque stroke } def
/BoxFill { gsave Rec 1 setgray fill grestore } def
/Symbol-Oblique /Symbol findfont [1 0 .167 1 0 0] makefont
dup length dict begin {1 index /FID eq {pop pop} {def} ifelse} forall
currentdict end definefont
end
}}%
\GNUPLOTpicture(2160,1511)
{\GNUPLOTspecial{"
gnudict begin
gsave
0 0 translate
0.100 0.100 scale
0 setgray
newpath
1.000 UL
LTb
250 300 M
63 0 V
1697 0 R
-63 0 V
250 410 M
63 0 V
1697 0 R
-63 0 V
250 520 M
63 0 V
1697 0 R
-63 0 V
250 630 M
63 0 V
1697 0 R
-63 0 V
250 740 M
63 0 V
1697 0 R
-63 0 V
250 850 M
63 0 V
1697 0 R
-63 0 V
250 961 M
63 0 V
1697 0 R
-63 0 V
250 1071 M
63 0 V
1697 0 R
-63 0 V
250 1181 M
63 0 V
1697 0 R
-63 0 V
250 1291 M
63 0 V
1697 0 R
-63 0 V
250 1401 M
63 0 V
1697 0 R
-63 0 V
250 1511 M
63 0 V
1697 0 R
-63 0 V
250 300 M
0 63 V
0 1148 R
0 -63 V
446 300 M
0 63 V
0 1148 R
0 -63 V
641 300 M
0 63 V
0 1148 R
0 -63 V
837 300 M
0 63 V
0 1148 R
0 -63 V
1032 300 M
0 63 V
0 1148 R
0 -63 V
1228 300 M
0 63 V
0 1148 R
0 -63 V
1423 300 M
0 63 V
0 1148 R
0 -63 V
1619 300 M
0 63 V
0 1148 R
0 -63 V
1814 300 M
0 63 V
0 1148 R
0 -63 V
2010 300 M
0 63 V
0 1148 R
0 -63 V
1.000 UL
LTb
250 300 M
1760 0 V
0 1211 V
-1760 0 V
250 300 L
1.000 UP
1.000 UL
LT0
423 545 Pls
1953 489 Pls
312 673 Pls
391 594 Pls
575 741 Pls
454 1024 Pls
393 889 Pls
441 498 Pls
597 716 Pls
300 730 Pls
929 804 Pls
584 725 Pls
446 978 Pls
404 586 Pls
567 792 Pls
571 675 Pls
337 645 Pls
373 617 Pls
374 616 Pls
374 616 Pls
378 868 Pls
378 867 Pls
333 832 Pls
572 807 Pls
1049 713 Pls
309 794 Pls
615 821 Pls
586 757 Pls
909 673 Pls
299 746 Pls
849 817 Pls
640 1136 Pls
309 684 Pls
1046 765 Pls
574 678 Pls
382 875 Pls
382 874 Pls
451 485 Pls
577 698 Pls
408 909 Pls
683 825 Pls
1217 1135 Pls
381 907 Pls
724 813 Pls
258 524 Pls
298 698 Pls
602 781 Pls
968 698 Pls
416 521 Pls
252 712 Pls
753 668 Pls
387 567 Pls
1384 367 Pls
811 662 Pls
415 555 Pls
1747 467 Pls
324 652 Pls
368 611 Pls
369 611 Pls
587 802 Pls
451 1059 Pls
300 775 Pls
775 678 Pls
373 564 Pls
817 703 Pls
349 545 Pls
368 918 Pls
780 800 Pls
694 657 Pls
406 571 Pls
581 732 Pls
338 639 Pls
339 638 Pls
357 623 Pls
357 622 Pls
358 622 Pls
358 622 Pls
359 622 Pls
361 861 Pls
360 860 Pls
360 860 Pls
359 860 Pls
337 840 Pls
448 478 Pls
1429 1037 Pls
406 910 Pls
697 824 Pls
288 715 Pls
643 754 Pls
351 541 Pls
374 917 Pls
774 804 Pls
428 961 Pls
404 570 Pls
586 780 Pls
1050 734 Pls
369 871 Pls
369 871 Pls
598 808 Pls
595 712 Pls
298 747 Pls
607 764 Pls
602 687 Pls
370 883 Pls
369 883 Pls
255 630 Pls
255 630 Pls
255 629 Pls
255 629 Pls
256 628 Pls
256 628 Pls
292 692 Pls
945 768 Pls
759 698 Pls
336 570 Pls
666 756 Pls
325 562 Pls
267 650 Pls
275 612 Pls
275 612 Pls
772 717 Pls
809 757 Pls
267 664 Pls
280 599 Pls
272 729 Pls
282 871 Pls
282 871 Pls
271 829 Pls
299 986 Pls
282 595 Pls
321 943 Pls
789 755 Pls
284 881 Pls
300 986 Pls
284 594 Pls
319 940 Pls
284 886 Pls
284 885 Pls
698 740 Pls
284 595 Pls
318 937 Pls
319 541 Pls
284 887 Pls
300 495 Pls
269 665 Pls
284 599 Pls
318 934 Pls
283 886 Pls
267 810 Pls
271 654 Pls
282 609 Pls
282 608 Pls
319 931 Pls
280 882 Pls
280 882 Pls
808 724 Pls
772 763 Pls
275 870 Pls
267 829 Pls
312 1017 Pls
275 765 Pls
667 725 Pls
280 776 Pls
639 709 Pls
360 903 Pls
731 802 Pls
632 671 Pls
368 601 Pls
368 600 Pls
2008 1006 Pls
853 667 Pls
418 549 Pls
596 728 Pls
329 644 Pls
358 618 Pls
358 617 Pls
359 617 Pls
593 803 Pls
956 687 Pls
424 529 Pls
386 953 Pls
642 750 Pls
282 753 Pls
806 786 Pls
362 877 Pls
362 877 Pls
634 808 Pls
388 974 Pls
944 766 Pls
353 968 Pls
874 760 Pls
286 997 Pls
265 643 Pls
269 619 Pls
269 618 Pls
269 618 Pls
276 705 Pls
318 467 Pls
279 585 Pls
327 963 Pls
315 568 Pls
287 868 Pls
286 867 Pls
276 832 Pls
701 753 Pls
284 568 Pls
309 570 Pls
295 889 Pls
272 807 Pls
290 950 Pls
732 724 Pls
304 575 Pls
271 792 Pls
288 937 Pls
721 726 Pls
298 584 Pls
306 907 Pls
736 758 Pls
290 525 Pls
292 598 Pls
815 733 Pls
323 516 Pls
262 783 Pls
825 729 Pls
329 518 Pls
276 888 Pls
681 726 Pls
325 910 Pls
284 475 Pls
342 905 Pls
747 789 Pls
984 723 Pls
358 886 Pls
672 806 Pls
261 539 Pls
302 670 Pls
358 596 Pls
270 1035 Pls
284 697 Pls
338 426 Pls
341 964 Pls
846 757 Pls
264 665 Pls
276 599 Pls
274 710 Pls
694 751 Pls
281 580 Pls
313 564 Pls
290 880 Pls
272 818 Pls
292 963 Pls
745 723 Pls
309 566 Pls
295 894 Pls
291 951 Pls
737 726 Pls
304 570 Pls
299 903 Pls
290 940 Pls
727 728 Pls
300 577 Pls
270 780 Pls
288 928 Pls
270 683 Pls
296 585 Pls
271 774 Pls
286 916 Pls
312 988 Pls
291 598 Pls
312 916 Pls
282 903 Pls
292 508 Pls
317 916 Pls
756 765 Pls
277 886 Pls
682 727 Pls
324 913 Pls
663 721 Pls
283 782 Pls
629 699 Pls
360 893 Pls
696 805 Pls
254 780 Pls
305 668 Pls
363 598 Pls
415 452 Pls
363 957 Pls
876 771 Pls
278 742 Pls
651 729 Pls
284 769 Pls
265 397 Pls
346 860 Pls
345 860 Pls
345 860 Pls
345 859 Pls
326 841 Pls
326 840 Pls
747 684 Pls
354 574 Pls
648 765 Pls
278 710 Pls
318 443 Pls
332 951 Pls
818 760 Pls
327 943 Pls
800 760 Pls
794 721 Pls
325 540 Pls
324 929 Pls
676 732 Pls
276 765 Pls
660 723 Pls
283 778 Pls
624 698 Pls
363 892 Pls
401 500 Pls
259 922 Pls
661 672 Pls
365 595 Pls
998 758 Pls
283 714 Pls
320 401 Pls
342 936 Pls
293 426 Pls
354 916 Pls
264 397 Pls
428 475 Pls
369 959 Pls
354 565 Pls
343 940 Pls
342 545 Pls
871 710 Pls
361 524 Pls
407 1041 Pls
367 891 Pls
677 809 Pls
678 672 Pls
300 679 Pls
987 762 Pls
361 957 Pls
871 771 Pls
341 936 Pls
810 773 Pls
351 916 Pls
775 790 Pls
358 875 Pls
358 875 Pls
636 807 Pls
255 659 Pls
258 604 Pls
293 686 Pls
948 764 Pls
282 695 Pls
880 760 Pls
263 658 Pls
271 605 Pls
746 709 Pls
845 753 Pls
279 575 Pls
274 694 Pls
292 873 Pls
291 873 Pls
276 826 Pls
261 723 Pls
322 978 Pls
304 584 Pls
301 893 Pls
721 759 Pls
286 536 Pls
275 665 Pls
295 599 Pls
309 902 Pls
734 763 Pls
287 517 Pls
317 906 Pls
333 511 Pls
275 889 Pls
676 722 Pls
328 905 Pls
280 471 Pls
350 902 Pls
376 508 Pls
417 1010 Pls
626 802 Pls
290 688 Pls
932 763 Pls
281 691 Pls
328 579 Pls
261 681 Pls
274 586 Pls
315 581 Pls
284 956 Pls
328 983 Pls
305 590 Pls
300 882 Pls
277 815 Pls
260 748 Pls
280 649 Pls
292 613 Pls
293 613 Pls
842 735 Pls
276 797 Pls
278 909 Pls
686 720 Pls
321 901 Pls
269 874 Pls
269 873 Pls
262 826 Pls
900 720 Pls
285 790 Pls
629 695 Pls
355 886 Pls
678 804 Pls
262 529 Pls
354 601 Pls
355 600 Pls
742 690 Pls
343 577 Pls
342 968 Pls
325 574 Pls
681 757 Pls
276 587 Pls
274 699 Pls
698 756 Pls
282 565 Pls
308 577 Pls
296 886 Pls
274 811 Pls
287 945 Pls
320 979 Pls
301 586 Pls
303 898 Pls
727 759 Pls
709 725 Pls
293 602 Pls
825 734 Pls
273 785 Pls
281 909 Pls
319 1007 Pls
275 781 Pls
274 884 Pls
286 493 Pls
329 906 Pls
749 779 Pls
352 902 Pls
378 508 Pls
622 679 Pls
619 800 Pls
351 577 Pls
281 1015 Pls
329 571 Pls
675 756 Pls
273 605 Pls
273 605 Pls
274 707 Pls
829 754 Pls
263 687 Pls
751 718 Pls
314 566 Pls
289 875 Pls
274 824 Pls
703 752 Pls
284 567 Pls
309 568 Pls
295 892 Pls
715 752 Pls
287 554 Pls
304 573 Pls
300 902 Pls
270 790 Pls
289 938 Pls
315 975 Pls
299 580 Pls
271 781 Pls
287 925 Pls
313 983 Pls
294 591 Pls
309 913 Pls
745 759 Pls
291 516 Pls
275 653 Pls
287 609 Pls
287 609 Pls
273 772 Pls
280 897 Pls
689 728 Pls
321 913 Pls
272 870 Pls
272 870 Pls
265 830 Pls
672 724 Pls
279 780 Pls
645 711 Pls
354 903 Pls
264 462 Pls
321 646 Pls
348 617 Pls
349 616 Pls
349 616 Pls
611 799 Pls
352 567 Pls
799 707 Pls
338 549 Pls
842 713 Pls
787 778 Pls
364 911 Pls
757 801 Pls
638 664 Pls
387 588 Pls
1027 764 Pls
587 693 Pls
392 898 Pls
654 818 Pls
389 963 Pls
368 570 Pls
790 701 Pls
342 556 Pls
834 711 Pls
796 777 Pls
363 913 Pls
634 664 Pls
389 588 Pls
445 449 Pls
585 691 Pls
391 894 Pls
640 818 Pls
383 961 Pls
915 775 Pls
278 726 Pls
658 740 Pls
279 756 Pls
640 717 Pls
289 774 Pls
647 662 Pls
392 584 Pls
585 789 Pls
585 689 Pls
388 891 Pls
764 682 Pls
362 568 Pls
644 760 Pls
341 552 Pls
843 711 Pls
349 530 Pls
638 716 Pls
367 911 Pls
259 455 Pls
304 688 Pls
1021 768 Pls
585 739 Pls
388 892 Pls
307 800 Pls
260 1042 Pls
363 569 Pls
790 703 Pls
340 555 Pls
830 712 Pls
796 775 Pls
357 913 Pls
262 414 Pls
316 659 Pls
368 605 Pls
596 798 Pls
941 684 Pls
298 757 Pls
796 679 Pls
376 558 Pls
358 926 Pls
802 789 Pls
250 554 Pls
359 555 Pls
640 731 Pls
290 758 Pls
805 665 Pls
1778 1398 Pls
stroke
grestore
end
showpage
}}%
\put(1597,1398){\rjust{$\tau (t=0) = 0.125$}}%
\put(-140,740){\ljust{$\alpha_1$}}%
\put(1130,50){\cjust{$\tau$}}%
\put(2010,200){\cjust{ 4.5}}%
\put(1814,200){\cjust{ 4}}%
\put(1619,200){\cjust{ 3.5}}%
\put(1423,200){\cjust{ 3}}%
\put(1228,200){\cjust{ 2.5}}%
\put(1032,200){\cjust{ 2}}%
\put(837,200){\cjust{ 1.5}}%
\put(641,200){\cjust{ 1}}%
\put(446,200){\cjust{ 0.5}}%
\put(250,200){\cjust{ 0}}%
\put(200,1511){\rjust{ 7}}%
\put(200,1401){\rjust{ 6}}%
\put(200,1291){\rjust{ 5}}%
\put(200,1181){\rjust{ 4}}%
\put(200,1071){\rjust{ 3}}%
\put(200,961){\rjust{ 2}}%
\put(200,850){\rjust{ 1}}%
\put(200,740){\rjust{ 0}}%
\put(200,630){\rjust{-1}}%
\put(200,520){\rjust{-2}}%
\put(200,410){\rjust{-3}}%
\put(200,300){\rjust{-4}}%
\endGNUPLOTpicture
\endgroup
 

%% file: fig_P_return_map2.tex
\begingroup
  \catcode`\@=11\relax
  \def\GNUPLOTspecial{%
    \def\do##1{\catcode`##1=12\relax}\dospecials
    \catcode`\{=1\catcode`\}=2\catcode\%=14\relax\special}%
\expandafter\ifx\csname GNUPLOTpicture\endcsname\relax
  \csname newdimen\endcsname\GNUPLOTunit
  \gdef\GNUPLOTpicture(#1,#2){\vbox to#2\GNUPLOTunit\bgroup
    \def\put(##1,##2)##3{\unskip\raise##2\GNUPLOTunit
      \hbox to0pt{\kern##1\GNUPLOTunit ##3\hss}\ignorespaces}%
    \def\ljust##1{\vbox to0pt{\vss\hbox to0pt{##1\hss}\vss}}%
    \def\cjust##1{\vbox to0pt{\vss\hbox to0pt{\hss ##1\hss}\vss}}%
    \def\rjust##1{\vbox to0pt{\vss\hbox to0pt{\hss ##1}\vss}}%
    \def\stack##1{\let\\=\cr\tabskip=0pt\halign{\hfil ####\hfil\cr ##1\crcr}}%
    \def\lstack##1{\hbox to0pt{\vbox to0pt{\vss\stack{##1}}\hss}}%
    \def\cstack##1{\hbox to0pt{\hss\vbox to0pt{\vss\stack{##1}}\hss}}%
    \def\rstack##1{\hbox to0pt{\vbox to0pt{\stack{##1}\vss}\hss}}%
    \vss\hbox to#1\GNUPLOTunit\bgroup\ignorespaces}%
  \gdef\endGNUPLOTpicture{\hss\egroup\egroup}%
\fi
\GNUPLOTunit=0.1bp
{\GNUPLOTspecial{!
/gnudict 256 dict def
gnudict begin
/Color false def
/Solid false def
/gnulinewidth 5.000 def
/userlinewidth gnulinewidth def
/vshift -33 def
/dl {10 mul} def
/hpt_ 31.5 def
/vpt_ 31.5 def
/hpt hpt_ def
/vpt vpt_ def
/M {moveto} bind def
/L {lineto} bind def
/R {rmoveto} bind def
/V {rlineto} bind def
/vpt2 vpt 2 mul def
/hpt2 hpt 2 mul def
/Lshow { currentpoint stroke M
  0 vshift R show } def
/Rshow { currentpoint stroke M
  dup stringwidth pop neg vshift R show } def
/Cshow { currentpoint stroke M
  dup stringwidth pop -2 div vshift R show } def
/UP { dup vpt_ mul /vpt exch def hpt_ mul /hpt exch def
  /hpt2 hpt 2 mul def /vpt2 vpt 2 mul def } def
/DL { Color {setrgbcolor Solid {pop []} if 0 setdash }
 {pop pop pop Solid {pop []} if 0 setdash} ifelse } def
/BL { stroke userlinewidth 2 mul setlinewidth } def
/AL { stroke userlinewidth 2 div setlinewidth } def
/UL { dup gnulinewidth mul /userlinewidth exch def
      dup 1 lt {pop 1} if 10 mul /udl exch def } def
/PL { stroke userlinewidth setlinewidth } def
/LTb { BL [] 0 0 0 DL } def
/LTa { AL [1 udl mul 2 udl mul] 0 setdash 0 0 0 setrgbcolor } def
/LT0 { PL [] 1 0 0 DL } def
/LT1 { PL [4 dl 2 dl] 0 1 0 DL } def
/LT2 { PL [2 dl 3 dl] 0 0 1 DL } def
/LT3 { PL [1 dl 1.5 dl] 1 0 1 DL } def
/LT4 { PL [5 dl 2 dl 1 dl 2 dl] 0 1 1 DL } def
/LT5 { PL [4 dl 3 dl 1 dl 3 dl] 1 1 0 DL } def
/LT6 { PL [2 dl 2 dl 2 dl 4 dl] 0 0 0 DL } def
/LT7 { PL [2 dl 2 dl 2 dl 2 dl 2 dl 4 dl] 1 0.3 0 DL } def
/LT8 { PL [2 dl 2 dl 2 dl 2 dl 2 dl 2 dl 2 dl 4 dl] 0.5 0.5 0.5 DL } def
/Pnt { stroke [] 0 setdash
   gsave 1 setlinecap M 0 0 V stroke grestore } def
/Dia { stroke [] 0 setdash 2 copy vpt add M
  hpt neg vpt neg V hpt vpt neg V
  hpt vpt V hpt neg vpt V closepath stroke
  Pnt } def
/Pls { stroke [] 0 setdash vpt sub M 0 vpt2 V
  currentpoint stroke M
  hpt neg vpt neg R hpt2 0 V stroke
  } def
/Box { stroke [] 0 setdash 2 copy exch hpt sub exch vpt add M
  0 vpt2 neg V hpt2 0 V 0 vpt2 V
  hpt2 neg 0 V closepath stroke
  Pnt } def
/Crs { stroke [] 0 setdash exch hpt sub exch vpt add M
  hpt2 vpt2 neg V currentpoint stroke M
  hpt2 neg 0 R hpt2 vpt2 V stroke } def
/TriU { stroke [] 0 setdash 2 copy vpt 1.12 mul add M
  hpt neg vpt -1.62 mul V
  hpt 2 mul 0 V
  hpt neg vpt 1.62 mul V closepath stroke
  Pnt  } def
/Star { 2 copy Pls Crs } def
/BoxF { stroke [] 0 setdash exch hpt sub exch vpt add M
  0 vpt2 neg V  hpt2 0 V  0 vpt2 V
  hpt2 neg 0 V  closepath fill } def
/TriUF { stroke [] 0 setdash vpt 1.12 mul add M
  hpt neg vpt -1.62 mul V
  hpt 2 mul 0 V
  hpt neg vpt 1.62 mul V closepath fill } def
/TriD { stroke [] 0 setdash 2 copy vpt 1.12 mul sub M
  hpt neg vpt 1.62 mul V
  hpt 2 mul 0 V
  hpt neg vpt -1.62 mul V closepath stroke
  Pnt  } def
/TriDF { stroke [] 0 setdash vpt 1.12 mul sub M
  hpt neg vpt 1.62 mul V
  hpt 2 mul 0 V
  hpt neg vpt -1.62 mul V closepath fill} def
/DiaF { stroke [] 0 setdash vpt add M
  hpt neg vpt neg V hpt vpt neg V
  hpt vpt V hpt neg vpt V closepath fill } def
/Pent { stroke [] 0 setdash 2 copy gsave
  translate 0 hpt M 4 {72 rotate 0 hpt L} repeat
  closepath stroke grestore Pnt } def
/PentF { stroke [] 0 setdash gsave
  translate 0 hpt M 4 {72 rotate 0 hpt L} repeat
  closepath fill grestore } def
/Circle { stroke [] 0 setdash 2 copy
  hpt 0 360 arc stroke Pnt } def
/CircleF { stroke [] 0 setdash hpt 0 360 arc fill } def
/C0 { BL [] 0 setdash 2 copy moveto vpt 90 450  arc } bind def
/C1 { BL [] 0 setdash 2 copy        moveto
       2 copy  vpt 0 90 arc closepath fill
               vpt 0 360 arc closepath } bind def
/C2 { BL [] 0 setdash 2 copy moveto
       2 copy  vpt 90 180 arc closepath fill
               vpt 0 360 arc closepath } bind def
/C3 { BL [] 0 setdash 2 copy moveto
       2 copy  vpt 0 180 arc closepath fill
               vpt 0 360 arc closepath } bind def
/C4 { BL [] 0 setdash 2 copy moveto
       2 copy  vpt 180 270 arc closepath fill
               vpt 0 360 arc closepath } bind def
/C5 { BL [] 0 setdash 2 copy moveto
       2 copy  vpt 0 90 arc
       2 copy moveto
       2 copy  vpt 180 270 arc closepath fill
               vpt 0 360 arc } bind def
/C6 { BL [] 0 setdash 2 copy moveto
      2 copy  vpt 90 270 arc closepath fill
              vpt 0 360 arc closepath } bind def
/C7 { BL [] 0 setdash 2 copy moveto
      2 copy  vpt 0 270 arc closepath fill
              vpt 0 360 arc closepath } bind def
/C8 { BL [] 0 setdash 2 copy moveto
      2 copy vpt 270 360 arc closepath fill
              vpt 0 360 arc closepath } bind def
/C9 { BL [] 0 setdash 2 copy moveto
      2 copy  vpt 270 450 arc closepath fill
              vpt 0 360 arc closepath } bind def
/C10 { BL [] 0 setdash 2 copy 2 copy moveto vpt 270 360 arc closepath fill
       2 copy moveto
       2 copy vpt 90 180 arc closepath fill
               vpt 0 360 arc closepath } bind def
/C11 { BL [] 0 setdash 2 copy moveto
       2 copy  vpt 0 180 arc closepath fill
       2 copy moveto
       2 copy  vpt 270 360 arc closepath fill
               vpt 0 360 arc closepath } bind def
/C12 { BL [] 0 setdash 2 copy moveto
       2 copy  vpt 180 360 arc closepath fill
               vpt 0 360 arc closepath } bind def
/C13 { BL [] 0 setdash  2 copy moveto
       2 copy  vpt 0 90 arc closepath fill
       2 copy moveto
       2 copy  vpt 180 360 arc closepath fill
               vpt 0 360 arc closepath } bind def
/C14 { BL [] 0 setdash 2 copy moveto
       2 copy  vpt 90 360 arc closepath fill
               vpt 0 360 arc } bind def
/C15 { BL [] 0 setdash 2 copy vpt 0 360 arc closepath fill
               vpt 0 360 arc closepath } bind def
/Rec   { newpath 4 2 roll moveto 1 index 0 rlineto 0 exch rlineto
       neg 0 rlineto closepath } bind def
/Square { dup Rec } bind def
/Bsquare { vpt sub exch vpt sub exch vpt2 Square } bind def
/S0 { BL [] 0 setdash 2 copy moveto 0 vpt rlineto BL Bsquare } bind def
/S1 { BL [] 0 setdash 2 copy vpt Square fill Bsquare } bind def
/S2 { BL [] 0 setdash 2 copy exch vpt sub exch vpt Square fill Bsquare } bind def
/S3 { BL [] 0 setdash 2 copy exch vpt sub exch vpt2 vpt Rec fill Bsquare } bind def
/S4 { BL [] 0 setdash 2 copy exch vpt sub exch vpt sub vpt Square fill Bsquare } bind def
/S5 { BL [] 0 setdash 2 copy 2 copy vpt Square fill
       exch vpt sub exch vpt sub vpt Square fill Bsquare } bind def
/S6 { BL [] 0 setdash 2 copy exch vpt sub exch vpt sub vpt vpt2 Rec fill Bsquare } bind def
/S7 { BL [] 0 setdash 2 copy exch vpt sub exch vpt sub vpt vpt2 Rec fill
       2 copy vpt Square fill
       Bsquare } bind def
/S8 { BL [] 0 setdash 2 copy vpt sub vpt Square fill Bsquare } bind def
/S9 { BL [] 0 setdash 2 copy vpt sub vpt vpt2 Rec fill Bsquare } bind def
/S10 { BL [] 0 setdash 2 copy vpt sub vpt Square fill 2 copy exch vpt sub exch vpt Square fill
       Bsquare } bind def
/S11 { BL [] 0 setdash 2 copy vpt sub vpt Square fill 2 copy exch vpt sub exch vpt2 vpt Rec fill
       Bsquare } bind def
/S12 { BL [] 0 setdash 2 copy exch vpt sub exch vpt sub vpt2 vpt Rec fill Bsquare } bind def
/S13 { BL [] 0 setdash 2 copy exch vpt sub exch vpt sub vpt2 vpt Rec fill
       2 copy vpt Square fill Bsquare } bind def
/S14 { BL [] 0 setdash 2 copy exch vpt sub exch vpt sub vpt2 vpt Rec fill
       2 copy exch vpt sub exch vpt Square fill Bsquare } bind def
/S15 { BL [] 0 setdash 2 copy Bsquare fill Bsquare } bind def
/D0 { gsave translate 45 rotate 0 0 S0 stroke grestore } bind def
/D1 { gsave translate 45 rotate 0 0 S1 stroke grestore } bind def
/D2 { gsave translate 45 rotate 0 0 S2 stroke grestore } bind def
/D3 { gsave translate 45 rotate 0 0 S3 stroke grestore } bind def
/D4 { gsave translate 45 rotate 0 0 S4 stroke grestore } bind def
/D5 { gsave translate 45 rotate 0 0 S5 stroke grestore } bind def
/D6 { gsave translate 45 rotate 0 0 S6 stroke grestore } bind def
/D7 { gsave translate 45 rotate 0 0 S7 stroke grestore } bind def
/D8 { gsave translate 45 rotate 0 0 S8 stroke grestore } bind def
/D9 { gsave translate 45 rotate 0 0 S9 stroke grestore } bind def
/D10 { gsave translate 45 rotate 0 0 S10 stroke grestore } bind def
/D11 { gsave translate 45 rotate 0 0 S11 stroke grestore } bind def
/D12 { gsave translate 45 rotate 0 0 S12 stroke grestore } bind def
/D13 { gsave translate 45 rotate 0 0 S13 stroke grestore } bind def
/D14 { gsave translate 45 rotate 0 0 S14 stroke grestore } bind def
/D15 { gsave translate 45 rotate 0 0 S15 stroke grestore } bind def
/DiaE { stroke [] 0 setdash vpt add M
  hpt neg vpt neg V hpt vpt neg V
  hpt vpt V hpt neg vpt V closepath stroke } def
/BoxE { stroke [] 0 setdash exch hpt sub exch vpt add M
  0 vpt2 neg V hpt2 0 V 0 vpt2 V
  hpt2 neg 0 V closepath stroke } def
/TriUE { stroke [] 0 setdash vpt 1.12 mul add M
  hpt neg vpt -1.62 mul V
  hpt 2 mul 0 V
  hpt neg vpt 1.62 mul V closepath stroke } def
/TriDE { stroke [] 0 setdash vpt 1.12 mul sub M
  hpt neg vpt 1.62 mul V
  hpt 2 mul 0 V
  hpt neg vpt -1.62 mul V closepath stroke } def
/PentE { stroke [] 0 setdash gsave
  translate 0 hpt M 4 {72 rotate 0 hpt L} repeat
  closepath stroke grestore } def
/CircE { stroke [] 0 setdash 
  hpt 0 360 arc stroke } def
/Opaque { gsave closepath 1 setgray fill grestore 0 setgray closepath } def
/DiaW { stroke [] 0 setdash vpt add M
  hpt neg vpt neg V hpt vpt neg V
  hpt vpt V hpt neg vpt V Opaque stroke } def
/BoxW { stroke [] 0 setdash exch hpt sub exch vpt add M
  0 vpt2 neg V hpt2 0 V 0 vpt2 V
  hpt2 neg 0 V Opaque stroke } def
/TriUW { stroke [] 0 setdash vpt 1.12 mul add M
  hpt neg vpt -1.62 mul V
  hpt 2 mul 0 V
  hpt neg vpt 1.62 mul V Opaque stroke } def
/TriDW { stroke [] 0 setdash vpt 1.12 mul sub M
  hpt neg vpt 1.62 mul V
  hpt 2 mul 0 V
  hpt neg vpt -1.62 mul V Opaque stroke } def
/PentW { stroke [] 0 setdash gsave
  translate 0 hpt M 4 {72 rotate 0 hpt L} repeat
  Opaque stroke grestore } def
/CircW { stroke [] 0 setdash 
  hpt 0 360 arc Opaque stroke } def
/BoxFill { gsave Rec 1 setgray fill grestore } def
/Symbol-Oblique /Symbol findfont [1 0 .167 1 0 0] makefont
dup length dict begin {1 index /FID eq {pop pop} {def} ifelse} forall
currentdict end definefont
end
}}%
\GNUPLOTpicture(2160,1511)
{\GNUPLOTspecial{"
gnudict begin
gsave
0 0 translate
0.100 0.100 scale
0 setgray
newpath
1.000 UL
LTb
250 300 M
63 0 V
1697 0 R
-63 0 V
250 473 M
63 0 V
1697 0 R
-63 0 V
250 646 M
63 0 V
1697 0 R
-63 0 V
250 819 M
63 0 V
1697 0 R
-63 0 V
250 992 M
63 0 V
1697 0 R
-63 0 V
250 1165 M
63 0 V
1697 0 R
-63 0 V
250 1338 M
63 0 V
1697 0 R
-63 0 V
250 1511 M
63 0 V
1697 0 R
-63 0 V
250 300 M
0 63 V
0 1148 R
0 -63 V
410 300 M
0 63 V
0 1148 R
0 -63 V
570 300 M
0 63 V
0 1148 R
0 -63 V
730 300 M
0 63 V
0 1148 R
0 -63 V
890 300 M
0 63 V
0 1148 R
0 -63 V
1050 300 M
0 63 V
0 1148 R
0 -63 V
1210 300 M
0 63 V
0 1148 R
0 -63 V
1370 300 M
0 63 V
0 1148 R
0 -63 V
1530 300 M
0 63 V
0 1148 R
0 -63 V
1690 300 M
0 63 V
0 1148 R
0 -63 V
1850 300 M
0 63 V
0 1148 R
0 -63 V
2010 300 M
0 63 V
0 1148 R
0 -63 V
1.000 UL
LTb
250 300 M
1760 0 V
0 1211 V
-1760 0 V
250 300 L
1.000 UL
LT0
1647 1398 M
263 0 V
444 819 M
0 1 V
0 1 V
0 1 V
0 1 V
0 1 V
-1 1 V
0 1 V
0 1 V
0 1 V
0 1 V
0 1 V
0 1 V
0 1 V
0 1 V
0 1 V
0 1 V
0 1 V
0 1 V
-1 1 V
0 1 V
0 1 V
0 1 V
0 1 V
0 1 V
0 1 V
0 1 V
-1 0 V
0 1 V
0 1 V
0 1 V
0 1 V
0 1 V
0 1 V
-1 1 V
0 1 V
0 1 V
0 1 V
0 1 V
-1 1 V
0 1 V
0 1 V
0 1 V
0 1 V
-1 1 V
0 1 V
0 1 V
0 1 V
0 1 V
-1 0 V
0 1 V
0 1 V
0 1 V
0 1 V
-1 1 V
0 1 V
0 1 V
0 1 V
-1 1 V
0 1 V
0 1 V
0 1 V
-1 1 V
0 1 V
0 1 V
0 1 V
0 1 V
-1 0 V
0 1 V
0 1 V
0 1 V
-1 1 V
0 1 V
0 1 V
0 1 V
-1 0 V
0 1 V
0 1 V
0 1 V
-1 1 V
0 1 V
0 1 V
-1 1 V
0 1 V
0 1 V
-1 1 V
0 1 V
0 1 V
-1 1 V
0 1 V
0 1 V
-1 0 V
0 1 V
0 1 V
-1 1 V
0 1 V
0 1 V
-1 0 V
0 1 V
-1 1 V
0 1 V
-1 1 V
0 1 V
0 1 V
-1 0 V
0 1 V
-1 1 V
0 1 V
-1 1 V
0 1 V
-1 0 V
0 1 V
-1 0 V
0 1 V
0 1 V
-1 0 V
0 1 V
0 1 V
-1 0 V
0 1 V
-1 0 V
0 1 V
-1 1 V
-1 0 V
0 1 V
-1 1 V
-1 1 V
-1 0 V
0 1 V
-1 0 V
0 1 V
-1 0 V
0 1 V
-1 0 V
0 1 V
-1 0 V
0 1 V
-1 0 V
-1 1 V
-1 0 V
0 1 V
-1 0 V
-1 1 V
-1 0 V
-1 0 V
0 1 V
-1 0 V
-1 0 V
0 1 V
-1 0 V
-1 0 V
-1 1 V
-1 0 V
-1 0 V
-1 0 V
0 1 V
-1 0 V
-1 0 V
-1 0 V
-1 0 V
-1 0 V
-1 0 V
-1 0 V
0 1 V
0 -1 V
0 1 V
-1 0 V
0 -1 V
-1 1 V
0 -1 V
-1 0 V
-1 1 V
0 -1 V
-1 0 V
-1 0 V
-1 0 V
-1 0 V
-1 0 V
-1 0 V
-1 0 V
0 -1 V
-1 0 V
-1 0 V
-1 0 V
-1 0 V
-1 -1 V
-1 0 V
-1 0 V
-1 -1 V
-1 0 V
-1 0 V
0 -1 V
-1 0 V
-1 0 V
-1 -1 V
-1 0 V
-1 -1 V
-1 0 V
0 -1 V
-1 0 V
-1 0 V
0 -1 V
-1 0 V
-1 -1 V
-1 0 V
0 -1 V
-1 0 V
0 -1 V
-1 0 V
-1 0 V
0 -1 V
-1 0 V
0 -1 V
-1 0 V
0 -1 V
-1 0 V
0 -1 V
-1 0 V
0 -1 V
-1 0 V
-1 -1 V
-1 -1 V
0 -1 V
-1 0 V
-1 -1 V
0 -1 V
-1 0 V
0 -1 V
-1 0 V
0 -1 V
-1 -1 V
0 -1 V
-1 0 V
-1 -1 V
0 -1 V
-1 0 V
0 -1 V
-1 -1 V
0 -1 V
-1 0 V
0 -1 V
-1 -1 V
0 -1 V
-1 -1 V
0 -1 V
-1 -1 V
0 -1 V
-1 0 V
0 -1 V
0 -1 V
-1 0 V
0 -1 V
0 -1 V
-1 -1 V
0 -1 V
-1 -1 V
0 -1 V
-1 -1 V
0 -1 V
0 -1 V
-1 -1 V
0 -1 V
-1 -1 V
0 -1 V
0 -1 V
-1 0 V
0 -1 V
0 -1 V
0 -1 V
-1 -1 V
0 -1 V
0 -1 V
-1 -1 V
0 -1 V
0 -1 V
0 -1 V
-1 -1 V
0 -1 V
0 -1 V
0 -1 V
-1 0 V
0 -1 V
0 -1 V
0 -1 V
0 -1 V
-1 -1 V
0 -1 V
0 -1 V
0 -1 V
0 -1 V
-1 -1 V
0 -1 V
0 -1 V
0 -1 V
0 -1 V
-1 -1 V
0 -1 V
0 -1 V
0 -1 V
0 -1 V
0 -1 V
0 -1 V
-1 -1 V
0 -1 V
0 -1 V
0 -1 V
0 -1 V
0 -1 V
0 -1 V
0 -1 V
0 -1 V
-1 -1 V
0 -1 V
0 -1 V
0 -1 V
0 -1 V
0 -1 V
0 -1 V
0 -1 V
0 -1 V
0 -1 V
0 -1 V
0 -1 V
0 -1 V
0 -1 V
-1 0 V
0 -1 V
0 -1 V
0 -1 V
0 -1 V
0 -1 V
0 -1 V
0 -1 V
0 -1 V
0 -1 V
0 -1 V
0 -1 V
0 -1 V
0 -1 V
0 -1 V
0 -1 V
0 -1 V
0 -1 V
0 -1 V
0 -1 V
0 -1 V
0 -1 V
0 -1 V
0 -1 V
0 -1 V
0 -1 V
0 -1 V
0 -1 V
0 -1 V
1 -1 V
0 -1 V
0 -1 V
0 -1 V
0 -1 V
0 -1 V
0 -1 V
0 -1 V
0 -1 V
0 -1 V
0 -1 V
0 -1 V
0 -1 V
1 -1 V
0 -1 V
0 -1 V
0 -1 V
0 -1 V
0 -1 V
0 -1 V
0 -1 V
0 -1 V
1 -1 V
0 -1 V
0 -1 V
0 -1 V
0 -1 V
0 -1 V
1 -1 V
0 -1 V
0 -1 V
0 -1 V
0 -1 V
0 -1 V
0 -1 V
1 0 V
0 -1 V
0 -1 V
0 -1 V
0 -1 V
0 -1 V
1 0 V
0 -1 V
0 -1 V
currentpoint stroke M
0 -1 V
1 -1 V
0 -1 V
0 -1 V
0 -1 V
0 -1 V
1 0 V
0 -1 V
0 -1 V
0 -1 V
1 -1 V
0 -1 V
0 -1 V
0 -1 V
1 0 V
0 -1 V
0 -1 V
0 -1 V
1 0 V
0 -1 V
0 -1 V
1 -1 V
0 -1 V
0 -1 V
1 -1 V
0 -1 V
1 -1 V
0 -1 V
0 -1 V
1 0 V
0 -1 V
0 -1 V
1 0 V
0 -1 V
0 -1 V
1 -1 V
0 -1 V
1 0 V
0 -1 V
0 -1 V
1 0 V
0 -1 V
1 -1 V
0 -1 V
1 0 V
0 -1 V
1 -1 V
0 -1 V
1 0 V
0 -1 V
1 -1 V
0 -1 V
1 0 V
0 -1 V
1 -1 V
0 -1 V
1 0 V
0 -1 V
1 0 V
0 -1 V
1 0 V
0 -1 V
1 -1 V
1 -1 V
1 -1 V
1 0 V
0 -1 V
1 -1 V
1 0 V
0 -1 V
1 0 V
0 -1 V
1 0 V
0 -1 V
1 0 V
0 -1 V
1 0 V
1 0 V
0 -1 V
1 0 V
0 -1 V
1 0 V
1 -1 V
1 0 V
0 -1 V
1 0 V
1 0 V
0 -1 V
1 0 V
1 0 V
0 -1 V
1 0 V
1 -1 V
1 0 V
1 0 V
0 -1 V
1 0 V
1 0 V
1 -1 V
1 0 V
1 0 V
1 -1 V
1 0 V
1 0 V
1 0 V
0 -1 V
1 1 V
0 -1 V
1 0 V
1 0 V
1 0 V
1 0 V
1 0 V
1 0 V
1 -1 V
0 1 V
1 0 V
1 0 V
0 -1 V
0 1 V
1 0 V
1 0 V
1 0 V
1 0 V
1 0 V
1 0 V
1 0 V
1 0 V
1 1 V
1 0 V
1 0 V
1 0 V
1 1 V
1 0 V
1 0 V
0 1 V
1 0 V
1 0 V
0 1 V
1 0 V
1 0 V
0 1 V
1 0 V
1 0 V
0 1 V
1 0 V
0 1 V
1 0 V
1 0 V
0 1 V
1 0 V
0 1 V
1 0 V
1 1 V
0 1 V
1 0 V
1 1 V
0 1 V
1 0 V
1 1 V
1 1 V
0 1 V
1 0 V
0 1 V
1 0 V
0 1 V
0 1 V
1 0 V
1 1 V
0 1 V
1 0 V
0 1 V
0 1 V
1 0 V
0 1 V
1 1 V
0 1 V
1 0 V
0 1 V
0 1 V
1 0 V
0 1 V
0 1 V
1 0 V
0 1 V
0 1 V
1 1 V
0 1 V
1 0 V
0 1 V
0 1 V
1 1 V
0 1 V
0 1 V
1 1 V
0 1 V
1 1 V
0 1 V
0 1 V
1 1 V
0 1 V
0 1 V
1 1 V
0 1 V
0 1 V
1 1 V
0 1 V
0 1 V
0 1 V
1 0 V
0 1 V
0 1 V
0 1 V
1 1 V
0 1 V
0 1 V
0 1 V
1 1 V
0 1 V
0 1 V
0 1 V
1 0 V
0 1 V
0 1 V
0 1 V
0 1 V
1 0 V
0 1 V
0 1 V
0 1 V
0 1 V
1 0 V
0 1 V
0 1 V
0 1 V
0 1 V
1 1 V
0 1 V
0 1 V
0 1 V
0 1 V
1 1 V
0 1 V
0 1 V
0 1 V
0 1 V
1 1 V
0 1 V
0 1 V
0 1 V
0 1 V
0 1 V
1 1 V
0 1 V
0 1 V
0 1 V
0 1 V
0 1 V
1 0 V
0 1 V
0 1 V
0 1 V
0 1 V
0 1 V
0 1 V
0 1 V
1 1 V
0 1 V
0 1 V
0 1 V
0 1 V
0 1 V
0 1 V
0 1 V
0 1 V
0 1 V
0 1 V
0 1 V
0 1 V
1 1 V
0 1 V
0 1 V
0 1 V
0 1 V
0 1 V
1.000 UL
LT0
1859 819 M
0 8 V
-1 7 V
-2 8 V
-3 6 V
-4 9 V
-5 8 V
-4 7 V
-6 6 V
-6 8 V
-6 6 V
-11 9 V
-7 6 V
-11 9 V
-8 6 V
-9 5 V
-12 8 V
-11 7 V
-10 5 V
-14 7 V
-16 9 V
-12 6 V
-9 4 V
-13 6 V
-12 5 V
-11 5 V
-9 4 V
-4 1 V
-3 2 V
-4 1 V
-4 2 V
-5 2 V
-4 1 V
-4 2 V
-5 2 V
-4 1 V
-3 2 V
-5 1 V
-4 2 V
-4 2 V
-5 1 V
-4 2 V
-7 2 V
-6 3 V
-6 2 V
-7 2 V
-6 2 V
-4 2 V
-5 1 V
-4 2 V
-4 1 V
-5 2 V
-4 1 V
-6 3 V
-7 2 V
-6 2 V
-5 2 V
-5 1 V
-5 2 V
-4 1 V
-5 2 V
-6 1 V
-7 3 V
-6 2 V
-4 1 V
-5 2 V
-5 1 V
-6 2 V
-7 2 V
-5 2 V
-5 1 V
-5 2 V
-5 1 V
-7 2 V
-5 2 V
-5 1 V
-5 2 V
-5 2 V
-8 2 V
-4 1 V
-5 1 V
-7 2 V
-6 2 V
-5 2 V
-5 1 V
-6 2 V
-7 2 V
-5 1 V
-6 2 V
-6 2 V
-5 1 V
-4 1 V
-6 2 V
-5 1 V
-5 2 V
-6 1 V
-6 2 V
-4 1 V
-6 2 V
-6 1 V
-5 2 V
-7 2 V
-6 1 V
-5 2 V
-7 1 V
-4 2 V
-6 1 V
-6 2 V
-5 1 V
-7 2 V
-5 1 V
-6 2 V
-6 2 V
-5 1 V
-6 2 V
-5 1 V
-6 1 V
-5 2 V
-5 1 V
-6 2 V
-5 1 V
-7 2 V
-5 1 V
-6 1 V
-5 2 V
-6 1 V
-5 1 V
-6 2 V
-5 1 V
-7 2 V
-5 1 V
-6 2 V
-6 1 V
-6 2 V
-5 1 V
-6 1 V
-6 2 V
-5 1 V
-6 2 V
-5 1 V
-6 1 V
-6 1 V
-6 2 V
-5 1 V
-6 1 V
-6 2 V
-5 1 V
-6 1 V
-5 2 V
-6 1 V
-6 1 V
-6 2 V
-6 1 V
-5 1 V
-6 1 V
-5 2 V
-6 1 V
-5 1 V
-6 1 V
-6 2 V
-6 1 V
-4 1 V
-7 1 V
-5 1 V
-5 1 V
-6 1 V
-5 2 V
-6 1 V
-4 0 V
-6 2 V
-5 1 V
-5 1 V
-5 1 V
-6 1 V
-6 1 V
-5 1 V
-5 0 V
-5 1 V
-5 1 V
-4 1 V
-6 1 V
-5 1 V
-5 1 V
-4 0 V
-5 1 V
-6 1 V
-4 0 V
-6 1 V
-4 1 V
-5 0 V
-6 1 V
-5 1 V
-4 0 V
-5 1 V
-5 1 V
-4 0 V
-5 0 V
-4 1 V
-5 0 V
-5 1 V
-4 0 V
-4 0 V
-4 1 V
-5 0 V
-5 0 V
-4 1 V
-4 0 V
-4 0 V
-3 0 V
-4 0 V
-5 0 V
-4 1 V
-3 0 V
-4 0 V
-4 0 V
-3 0 V
-4 -1 V
-4 0 V
-3 0 V
-3 0 V
-4 0 V
-3 0 V
-4 0 V
-3 -1 V
-3 0 V
-3 0 V
-4 0 V
-2 -1 V
-3 0 V
-3 0 V
-3 -1 V
-3 0 V
-2 -1 V
-4 0 V
-2 -1 V
-2 0 V
-3 0 V
-3 -1 V
-2 -1 V
-3 0 V
-2 -1 V
-2 0 V
-3 -1 V
-2 -1 V
-3 0 V
-2 -1 V
-2 -1 V
-2 0 V
-2 -1 V
-1 0 V
-2 -1 V
-2 -1 V
-1 0 V
-2 -1 V
-2 -1 V
-2 0 V
-1 -1 V
-2 -1 V
-2 -1 V
-1 0 V
-2 -1 V
-1 -1 V
-2 -1 V
-1 0 V
-1 -1 V
-2 -1 V
-1 -1 V
-1 -1 V
-1 0 V
-2 -1 V
-1 -1 V
-1 -1 V
-1 0 V
-1 -1 V
-1 -1 V
-1 -1 V
-1 0 V
-1 -1 V
-2 -1 V
0 -1 V
-1 -1 V
-1 -1 V
-1 0 V
-1 -1 V
-1 -1 V
-1 -1 V
-1 0 V
0 -1 V
-1 -1 V
-1 -1 V
-1 0 V
-1 -1 V
0 -1 V
-1 -1 V
-1 -1 V
0 -1 V
-1 0 V
-1 -1 V
-1 -1 V
0 -1 V
-1 0 V
0 -1 V
-1 -1 V
-1 -1 V
0 -1 V
-1 0 V
0 -1 V
-1 -1 V
0 -1 V
-1 -1 V
-1 0 V
0 -1 V
-1 -1 V
0 -1 V
-1 -1 V
-1 -1 V
0 -1 V
-1 -1 V
0 -1 V
-1 -1 V
0 -1 V
-1 -1 V
0 -1 V
-1 -1 V
0 -1 V
-1 -1 V
0 -1 V
-1 -1 V
0 -1 V
-1 -1 V
0 -1 V
-1 -1 V
0 -1 V
0 -1 V
-1 -1 V
-1 -1 V
0 -1 V
0 -1 V
0 -1 V
-1 -1 V
0 -1 V
-1 0 V
0 -1 V
0 -1 V
0 -1 V
-1 -1 V
0 -1 V
-1 -1 V
1 0 V
-1 -1 V
0 -1 V
0 -1 V
-1 -1 V
0 -1 V
0 -1 V
-1 0 V
0 -1 V
0 -1 V
-1 -1 V
0 -1 V
0 -1 V
0 -1 V
0 -1 V
-1 -1 V
0 -1 V
0 -1 V
0 -1 V
-1 0 V
0 -1 V
0 -1 V
0 -1 V
0 -1 V
-1 0 V
0 -1 V
0 -1 V
0 -1 V
-1 0 V
1 -1 V
-1 -1 V
0 -1 V
0 -1 V
0 -1 V
-1 0 V
0 -1 V
0 -1 V
0 -1 V
0 -1 V
0 -1 V
-1 -1 V
0 -1 V
0 -1 V
0 -1 V
0 -1 V
-1 0 V
1 -1 V
-1 -1 V
1 0 V
-1 -1 V
0 -1 V
0 -1 V
currentpoint stroke M
-1 -1 V
1 0 V
-1 -1 V
1 0 V
-1 -1 V
1 -1 V
-1 0 V
0 -1 V
0 -1 V
0 -1 V
0 -1 V
-1 0 V
1 -1 V
0 -1 V
-1 0 V
1 0 V
0 -1 V
-1 0 V
0 -1 V
1 0 V
0 -1 V
-1 -1 V
1 0 V
-1 -1 V
1 0 V
-1 0 V
0 -1 V
1 0 V
-1 0 V
1 -1 V
-1 0 V
1 -1 V
0 -1 V
-1 0 V
1 0 V
-1 0 V
1 -1 V
-1 0 V
1 0 V
0 -1 V
-1 0 V
1 0 V
-1 -1 V
1 0 V
-1 0 V
1 -1 V
0 -1 V
0 -1 V
0 -1 V
0 -1 V
0 -1 V
1 0 V
-1 0 V
1 0 V
-1 -1 V
1 0 V
0 -1 V
0 -1 V
0 -1 V
1 -1 V
-1 0 V
1 0 V
0 -1 V
-1 0 V
1 -1 V
1 0 V
-1 0 V
0 -1 V
1 0 V
-1 0 V
1 -1 V
0 -1 V
1 -1 V
-1 0 V
1 0 V
-1 0 V
1 -1 V
1 0 V
-1 0 V
1 -1 V
-1 0 V
1 0 V
-1 0 V
1 0 V
0 -1 V
1 0 V
-1 0 V
1 0 V
-1 -1 V
1 0 V
-1 0 V
1 0 V
-1 0 V
1 -1 V
1 0 V
-1 0 V
1 0 V
-1 0 V
1 0 V
-1 -1 V
2 0 V
-1 0 V
1 0 V
0 -1 V
-1 0 V
1 0 V
-1 0 V
1 0 V
-1 0 V
1 -1 V
1 0 V
-1 0 V
1 -1 V
-1 0 V
1 0 V
1 0 V
0 -1 V
-1 0 V
1 0 V
1 -1 V
1 0 V
-1 0 V
1 0 V
0 -1 V
-1 0 V
1 0 V
1 0 V
-1 0 V
1 -1 V
1 0 V
-2 0 V
2 0 V
-1 0 V
1 0 V
1 -1 V
-2 0 V
2 0 V
-1 0 V
1 0 V
1 0 V
-1 -1 V
1 0 V
-1 0 V
1 0 V
1 0 V
-1 0 V
1 -1 V
-1 1 V
1 -1 V
1 0 V
-1 0 V
1 0 V
0 -1 V
-1 0 V
1 0 V
1 0 V
-1 0 V
2 0 V
-1 0 V
1 0 V
-1 0 V
1 -1 V
-1 0 V
1 0 V
1 0 V
1 0 V
-2 0 V
2 0 V
-1 0 V
1 -1 V
1 0 V
-2 0 V
1 0 V
1 0 V
-1 0 V
1 0 V
1 -1 V
-1 0 V
1 0 V
1 0 V
1 0 V
-2 0 V
1 0 V
-1 0 V
2 -1 V
-1 0 V
1 0 V
2 0 V
-2 0 V
1 0 V
1 0 V
-2 0 V
2 -1 V
1 0 V
-1 0 V
1 0 V
1 0 V
-2 0 V
1 0 V
1 0 V
-1 0 V
1 0 V
1 -1 V
1 0 V
-2 0 V
1 0 V
1 0 V
1 0 V
1 0 V
-1 0 V
0 -1 V
1 0 V
-1 0 V
1 0 V
1 0 V
-2 0 V
2 0 V
1 0 V
-1 0 V
-1 0 V
1 0 V
1 -1 V
1 0 V
1 0 V
1 -1 V
1 0 V
-1 0 V
2 0 V
2 0 V
-2 0 V
1 0 V
0 -1 V
2 0 V
-1 0 V
2 0 V
-2 0 V
1 0 V
3 0 V
-1 0 V
1 0 V
0 -1 V
-1 0 V
1 0 V
2 0 V
2 0 V
-4 0 V
2 0 V
3 0 V
-1 0 V
-2 0 V
3 -1 V
1 0 V
-2 0 V
3 0 V
-2 0 V
1 0 V
-2 0 V
2 0 V
-1 0 V
4 0 V
-2 0 V
-1 0 V
2 0 V
-1 0 V
1 0 V
0 -1 V
3 0 V
1 0 V
-2 0 V
1 0 V
-1 0 V
1 0 V
1 0 V
2 0 V
-2 0 V
2 0 V
-2 0 V
4 -1 V
-1 0 V
1 0 V
3 0 V
1 0 V
-1 0 V
1 0 V
-1 0 V
2 -1 V
-1 0 V
2 0 V
-2 0 V
4 0 V
-1 0 V
2 0 V
1 0 V
-3 0 V
2 0 V
-1 0 V
3 -1 V
3 0 V
-3 0 V
4 0 V
-2 0 V
1 0 V
1 0 V
-1 0 V
1 0 V
2 0 V
-3 0 V
4 0 V
-3 0 V
4 -1 V
-2 0 V
2 0 V
-2 0 V
3 0 V
-2 0 V
2 0 V
-2 0 V
2 0 V
-1 0 V
2 0 V
3 0 V
1 0 V
-1 0 V
-2 0 V
1 0 V
3 -1 V
-4 1 V
3 -1 V
-1 0 V
1 0 V
2 0 V
-3 0 V
2 0 V
-1 0 V
2 0 V
1 0 V
3 0 V
-1 0 V
1 0 V
-1 0 V
-1 0 V
1 0 V
1 0 V
-1 0 V
4 -1 V
1 0 V
-3 0 V
1 0 V
3 0 V
-2 0 V
2 0 V
-1 0 V
1 0 V
3 0 V
2 0 V
-4 0 V
3 -1 V
1 0 V
-2 0 V
2 0 V
1 0 V
-1 0 V
3 0 V
-1 0 V
2 0 V
2 0 V
-1 0 V
2 0 V
-2 0 V
5 -1 V
-3 0 V
4 0 V
-3 0 V
1 0 V
4 0 V
-3 0 V
3 0 V
1 0 V
0 -1 V
1 0 V
5 0 V
1 0 V
1 0 V
1 0 V
0 -1 V
1 0 V
1 0 V
1 0 V
1 0 V
1 0 V
1 0 V
1 0 V
3 -1 V
1 0 V
-2 0 V
2 0 V
1 0 V
-1 0 V
1 0 V
3 0 V
currentpoint stroke M
-1 0 V
2 0 V
1 0 V
0 -1 V
2 0 V
-1 0 V
2 0 V
3 0 V
-1 0 V
2 0 V
-2 0 V
1 0 V
1 0 V
2 -1 V
-1 0 V
2 0 V
2 0 V
1 0 V
1 0 V
-1 0 V
2 0 V
2 0 V
-2 0 V
2 -1 V
1 0 V
1 0 V
-1 0 V
2 0 V
1 0 V
-1 0 V
2 0 V
-1 0 V
2 0 V
2 -1 V
1 0 V
1 0 V
1 0 V
1 0 V
2 0 V
-1 0 V
1 -1 V
1 0 V
1 0 V
1 0 V
-1 0 V
2 0 V
-1 0 V
2 0 V
1 0 V
1 0 V
2 -1 V
1 0 V
2 0 V
2 0 V
1 0 V
1 -1 V
1 0 V
1 0 V
1 0 V
1 0 V
1 0 V
1 0 V
1 -1 V
1 0 V
2 0 V
1 0 V
1 0 V
1 -1 V
2 0 V
1 0 V
1 0 V
1 0 V
1 0 V
2 -1 V
1 0 V
1 0 V
1 0 V
1 0 V
1 0 V
0 -1 V
2 0 V
1 0 V
1 0 V
1 0 V
1 0 V
1 -1 V
1 0 V
1 0 V
1 0 V
1 0 V
2 0 V
0 -1 V
1 0 V
1 0 V
1 0 V
1 0 V
1 0 V
1 0 V
1 -1 V
1 0 V
1 0 V
1 0 V
1 0 V
1 -1 V
1 0 V
1 0 V
1 0 V
1 0 V
2 0 V
0 -1 V
1 0 V
1 0 V
1 0 V
1 0 V
2 -1 V
1 0 V
1 0 V
1 0 V
1 0 V
1 0 V
1 -1 V
1 0 V
1 0 V
1 0 V
1 0 V
1 -1 V
1 0 V
1 0 V
1 0 V
1 0 V
1 -1 V
1 0 V
1 0 V
1 0 V
1 0 V
1 -1 V
1 0 V
1 0 V
1 0 V
1 0 V
1 -1 V
1 0 V
1 0 V
1 0 V
1 0 V
1 -1 V
1 0 V
1 0 V
2 0 V
3 -1 V
3 0 V
2 -1 V
1 0 V
3 -1 V
2 0 V
2 -1 V
3 0 V
2 -1 V
2 0 V
3 -1 V
2 0 V
2 -1 V
3 0 V
2 -1 V
2 0 V
3 -1 V
4 -1 V
3 -1 V
4 -1 V
4 -1 V
2 0 V
3 -1 V
2 -1 V
3 0 V
2 -1 V
2 -1 V
2 0 V
4 -1 V
3 -1 V
4 -1 V
3 -1 V
3 -1 V
2 0 V
2 -1 V
3 -1 V
3 -1 V
3 -1 V
3 -1 V
2 0 V
3 -1 V
1 -1 V
3 -1 V
3 -1 V
4 -1 V
2 -1 V
3 0 V
2 -1 V
3 -1 V
3 -1 V
2 -1 V
2 -1 V
2 -1 V
3 0 V
3 -2 V
3 -1 V
2 -1 V
2 0 V
3 -2 V
2 -1 V
2 0 V
3 -1 V
2 -1 V
3 -2 V
3 -1 V
1 0 V
4 -2 V
2 -1 V
2 -1 V
3 -1 V
2 -1 V
2 -1 V
2 -1 V
3 -2 V
2 -1 V
2 -1 V
3 -1 V
2 -1 V
2 -1 V
2 -1 V
2 -1 V
2 -1 V
3 -2 V
1 -1 V
2 -1 V
3 -1 V
1 -2 V
2 -1 V
2 -1 V
2 -1 V
2 -1 V
2 -2 V
2 -1 V
2 -1 V
1 -1 V
2 -2 V
2 -1 V
1 -1 V
3 -2 V
1 -1 V
1 -1 V
2 -1 V
1 -1 V
2 -2 V
2 -1 V
1 -1 V
2 -1 V
1 -2 V
1 -1 V
2 -1 V
1 -2 V
1 -1 V
2 -1 V
1 -1 V
1 -2 V
1 -1 V
2 -2 V
1 -1 V
1 -1 V
1 -1 V
1 -1 V
1 -2 V
1 -1 V
1 -2 V
1 -1 V
1 -1 V
1 -1 V
1 -2 V
0 -1 V
1 -2 V
1 -1 V
1 -1 V
1 -1 V
0 -2 V
1 -1 V
1 -2 V
1 -1 V
0 -1 V
1 -2 V
0 -1 V
1 -2 V
0 -1 V
1 -1 V
0 -2 V
1 -1 V
0 -1 V
1 -2 V
0 -1 V
0 -2 V
0 -1 V
1 -1 V
0 -2 V
0 -1 V
0 -1 V
1 -2 V
0 -1 V
0 -2 V
0 -1 V
0 -1 V
0 -2 V
0 -1 V
0 -2 V
0 -2 V
0 -1 V
0 -2 V
0 -1 V
-1 -1 V
0 -2 V
0 -1 V
0 -2 V
-1 -1 V
0 -1 V
0 -2 V
-1 -1 V
0 -1 V
0 -2 V
-1 -1 V
0 -2 V
-1 -1 V
0 -1 V
-1 -2 V
0 -1 V
-1 -1 V
-1 -2 V
0 -1 V
-1 -2 V
-1 -1 V
0 -2 V
-1 -1 V
-1 -1 V
-1 -2 V
-1 -1 V
-1 -2 V
0 -1 V
-1 -1 V
-1 -1 V
-1 -2 V
-1 -1 V
-1 -1 V
-1 -2 V
-2 -1 V
-1 -1 V
-1 -2 V
-1 -1 V
-1 -1 V
-1 -2 V
-2 -1 V
-1 -1 V
-1 -2 V
-2 -1 V
-1 -1 V
-1 -2 V
-1 -1 V
-2 -1 V
-1 -1 V
-2 -1 V
-2 -2 V
-1 -1 V
-3 -2 V
-1 -1 V
-2 -2 V
-1 -1 V
-2 -1 V
-2 -2 V
-2 -1 V
-2 -1 V
-1 -1 V
-2 -2 V
-2 -1 V
-2 -1 V
-3 -2 V
-1 -1 V
-2 -1 V
-2 -1 V
-2 -1 V
-2 -1 V
-2 -2 V
-2 0 V
-2 -2 V
-2 -1 V
-2 -1 V
-3 -1 V
-2 -1 V
-2 -1 V
-3 -2 V
-2 -1 V
-2 -1 V
-2 -1 V
-3 -1 V
-2 -1 V
-2 -1 V
currentpoint stroke M
-3 -1 V
-2 -1 V
-2 -1 V
-3 -1 V
-2 -2 V
-3 -1 V
-1 0 V
-2 -1 V
-4 -2 V
-2 -1 V
-2 0 V
-2 -1 V
-3 -1 V
-3 -1 V
-2 -1 V
-2 -1 V
-2 -1 V
-5 -1 V
-2 -1 V
-2 -1 V
-3 -1 V
-3 -1 V
-4 -1 V
-2 -1 V
-2 -1 V
-3 0 V
-2 -1 V
-3 -1 V
-4 -1 V
-2 -1 V
-3 -1 V
-2 -1 V
-2 0 V
-3 -1 V
-4 -1 V
-3 -1 V
-2 -1 V
-2 0 V
-3 -1 V
-2 -1 V
-2 0 V
-3 -1 V
-2 -1 V
-4 -1 V
-4 -1 V
-3 0 V
-3 -1 V
-2 -1 V
-2 0 V
-3 -1 V
-2 0 V
-2 -1 V
-3 -1 V
-2 0 V
-2 -1 V
-3 0 V
-3 -1 V
-4 -1 V
-2 0 V
-5 -1 V
-3 -1 V
-4 -1 V
-2 0 V
-3 -1 V
-2 -1 V
-5 0 V
-3 -1 V
-1 0 V
-1 -1 V
-1 0 V
-1 0 V
-1 0 V
-1 0 V
-1 -1 V
-1 0 V
-1 0 V
-1 0 V
-1 0 V
-1 -1 V
-1 0 V
-1 0 V
-1 0 V
-1 0 V
-1 0 V
-1 -1 V
-2 0 V
-1 0 V
-1 0 V
-1 -1 V
-1 0 V
-1 0 V
-1 0 V
-2 0 V
-1 -1 V
-1 0 V
-1 0 V
-1 0 V
-1 0 V
-1 -1 V
-1 0 V
-1 0 V
-1 0 V
-1 0 V
-1 0 V
-1 -1 V
-2 0 V
-1 0 V
-1 0 V
-1 0 V
0 -1 V
-1 0 V
-1 0 V
-1 0 V
-2 0 V
-1 0 V
0 -1 V
-2 0 V
-1 0 V
-1 0 V
-1 0 V
-1 0 V
-1 -1 V
-1 0 V
-1 0 V
-1 0 V
-1 0 V
-1 -1 V
-2 0 V
-1 0 V
-1 0 V
-1 0 V
-1 0 V
0 -1 V
-2 0 V
-1 0 V
-1 0 V
-1 0 V
-1 0 V
-1 0 V
-1 -1 V
-1 0 V
-1 0 V
-1 0 V
-1 0 V
-1 0 V
-1 0 V
-1 -1 V
-1 0 V
-1 0 V
-1 0 V
-1 0 V
-2 0 V
-1 -1 V
-1 0 V
-1 0 V
-2 0 V
-1 0 V
-2 -1 V
-1 0 V
-2 0 V
-1 0 V
1 0 V
-2 0 V
-1 0 V
0 -1 V
-3 0 V
-1 0 V
-2 0 V
-1 0 V
0 -1 V
-2 0 V
-1 0 V
1 0 V
-1 0 V
-2 0 V
-1 0 V
1 0 V
-1 0 V
-1 0 V
-1 0 V
-1 -1 V
1 0 V
-2 0 V
-2 0 V
-2 0 V
-1 0 V
-2 0 V
1 -1 V
-1 0 V
-3 0 V
1 0 V
-2 0 V
-1 0 V
-2 0 V
1 0 V
-2 -1 V
-2 0 V
1 0 V
-2 0 V
1 0 V
-2 0 V
-1 0 V
-1 0 V
-1 -1 V
1 0 V
-3 0 V
-1 0 V
1 0 V
-1 0 V
-3 0 V
2 0 V
-3 0 V
2 0 V
-2 0 V
-1 -1 V
-1 0 V
-1 0 V
-1 0 V
-1 0 V
-1 0 V
-1 0 V
-1 0 V
-1 0 V
0 -1 V
-1 0 V
-3 0 V
2 0 V
-2 0 V
-1 0 V
-1 0 V
-1 0 V
-1 0 V
0 -1 V
-4 0 V
-1 0 V
3 0 V
-1 0 V
-3 0 V
-3 0 V
2 0 V
-2 0 V
2 0 V
-2 0 V
2 0 V
-5 -1 V
2 0 V
-3 0 V
-1 0 V
2 0 V
-3 0 V
1 0 V
2 0 V
-3 0 V
1 0 V
-1 0 V
2 0 V
-3 0 V
-2 0 V
-1 -1 V
2 0 V
-2 0 V
2 0 V
-2 0 V
-1 0 V
-3 0 V
1 0 V
1 0 V
1 0 V
-4 0 V
1 0 V
-2 0 V
2 0 V
1 0 V
-1 0 V
-3 0 V
1 0 V
-2 -1 V
1 0 V
-2 0 V
-1 0 V
1 0 V
-1 0 V
2 0 V
-1 0 V
-1 0 V
1 0 V
-1 0 V
-2 0 V
1 0 V
-5 -1 V
2 0 V
-1 0 V
2 0 V
-1 0 V
-3 0 V
-2 0 V
1 0 V
1 0 V
-1 0 V
-3 0 V
1 0 V
-2 -1 V
-1 0 V
1 0 V
1 0 V
-2 0 V
2 0 V
-1 0 V
-2 0 V
-3 0 V
1 0 V
1 0 V
-2 0 V
-2 0 V
1 0 V
2 0 V
-1 0 V
1 0 V
-2 -1 V
-1 0 V
-1 0 V
-2 0 V
-1 0 V
2 0 V
-2 0 V
2 0 V
-1 0 V
1 0 V
-3 0 V
-1 0 V
0 -1 V
1 0 V
-1 0 V
-1 0 V
1 0 V
-1 0 V
-3 0 V
1 0 V
-2 0 V
4 0 V
-2 0 V
-2 0 V
2 0 V
-1 0 V
-2 -1 V
1 0 V
1 0 V
-2 0 V
1 0 V
-1 0 V
-1 0 V
-2 0 V
1 0 V
1 0 V
-4 0 V
3 0 V
-1 0 V
-2 0 V
2 0 V
1 0 V
-2 0 V
1 0 V
-4 -1 V
3 0 V
-2 0 V
-2 0 V
2 0 V
-2 0 V
2 0 V
-3 0 V
2 0 V
-3 0 V
0 -1 V
-1 0 V
3 0 V
-3 0 V
-1 0 V
1 0 V
-1 0 V
-1 0 V
2 0 V
-2 0 V
1 0 V
-1 -1 V
-2 0 V
1 0 V
-1 0 V
-3 0 V
2 0 V
-2 0 V
0 -1 V
1 0 V
-1 0 V
-1 0 V
-1 0 V
2 0 V
-3 0 V
1 0 V
0 -1 V
-3 0 V
currentpoint stroke M
2 0 V
-1 0 V
-1 0 V
1 0 V
-1 0 V
-2 -1 V
-1 0 V
-1 0 V
1 0 V
-1 0 V
-1 0 V
1 0 V
-1 -1 V
-1 0 V
-1 0 V
2 0 V
-1 0 V
-1 0 V
0 -1 V
-1 0 V
-1 0 V
2 0 V
-2 0 V
1 0 V
-1 0 V
-1 0 V
1 0 V
-1 0 V
0 -1 V
2 0 V
-1 0 V
-1 0 V
-2 0 V
1 0 V
-1 0 V
0 -1 V
-1 0 V
1 0 V
-2 0 V
1 0 V
-1 -1 V
1 0 V
-2 0 V
1 0 V
-2 0 V
2 0 V
-2 0 V
2 0 V
-1 0 V
-1 0 V
1 0 V
-2 -1 V
1 0 V
-1 0 V
1 0 V
-1 0 V
1 0 V
-2 0 V
2 0 V
-1 0 V
0 -1 V
-1 0 V
-1 0 V
1 0 V
-1 0 V
-1 -1 V
2 0 V
-2 0 V
1 0 V
-1 0 V
-1 0 V
1 0 V
0 -1 V
-1 0 V
-1 0 V
1 0 V
-1 0 V
1 0 V
-1 0 V
1 -1 V
-1 0 V
-1 0 V
1 0 V
-1 0 V
0 -1 V
-1 0 V
1 0 V
-1 0 V
1 0 V
-2 0 V
1 -1 V
-1 0 V
1 0 V
-1 0 V
-1 -1 V
1 0 V
-1 0 V
-1 -1 V
1 0 V
-1 0 V
-1 -1 V
1 0 V
-1 0 V
-1 0 V
1 -1 V
-1 0 V
1 0 V
-1 0 V
0 -1 V
-1 0 V
1 0 V
0 -1 V
-1 0 V
1 0 V
-1 0 V
-1 -1 V
1 0 V
-1 0 V
0 -1 V
-1 0 V
0 -1 V
1 0 V
-1 0 V
-1 0 V
1 0 V
-1 -1 V
1 0 V
-1 0 V
1 0 V
-1 0 V
1 -1 V
-1 0 V
-1 -1 V
1 0 V
-1 0 V
1 0 V
-1 0 V
0 -1 V
1 0 V
-1 0 V
1 0 V
-1 0 V
0 -1 V
-1 0 V
1 0 V
-1 -1 V
1 0 V
-1 0 V
1 0 V
-1 0 V
0 -1 V
1 0 V
-1 0 V
0 -1 V
-1 0 V
0 -1 V
1 0 V
-1 0 V
1 -1 V
-1 0 V
1 0 V
-1 -1 V
-1 -1 V
1 0 V
-1 0 V
1 0 V
0 -1 V
-1 0 V
1 0 V
-1 -1 V
1 0 V
0 -1 V
-1 0 V
1 0 V
-1 0 V
0 -1 V
0 -1 V
0 -1 V
0 -1 V
-1 -1 V
1 0 V
0 -1 V
-1 0 V
1 0 V
0 -1 V
-1 0 V
0 -1 V
1 0 V
-1 -1 V
1 0 V
-1 -1 V
0 -1 V
1 0 V
-1 0 V
0 -1 V
1 0 V
-1 0 V
1 -1 V
-1 0 V
1 0 V
0 -1 V
-1 0 V
1 0 V
0 -1 V
-1 0 V
1 0 V
0 -1 V
0 -1 V
-1 -1 V
1 0 V
0 -1 V
0 -1 V
0 -1 V
0 -1 V
0 -1 V
1 0 V
-1 -1 V
1 -1 V
0 -1 V
-1 0 V
1 -1 V
0 -1 V
0 -1 V
0 -1 V
0 -1 V
1 -1 V
0 -1 V
0 -1 V
0 -1 V
0 -1 V
0 -1 V
1 0 V
0 -1 V
0 -1 V
0 -1 V
0 -1 V
1 0 V
-1 -1 V
1 0 V
0 -1 V
-1 0 V
1 -1 V
0 -1 V
0 -1 V
0 -1 V
1 -1 V
0 -1 V
0 -1 V
0 -1 V
1 0 V
0 -1 V
0 -1 V
0 -1 V
1 -1 V
0 -1 V
0 -1 V
0 -1 V
1 0 V
0 -1 V
0 -1 V
0 -1 V
1 0 V
-1 -1 V
1 0 V
0 -1 V
0 -1 V
0 -1 V
1 0 V
0 -1 V
0 -1 V
0 -1 V
1 0 V
0 -1 V
0 -1 V
0 -1 V
1 0 V
0 -1 V
0 -1 V
0 -1 V
1 0 V
0 -1 V
0 -1 V
1 0 V
0 -1 V
0 -1 V
1 -1 V
0 -1 V
0 -1 V
1 -1 V
0 -1 V
1 -1 V
0 -1 V
0 -1 V
1 -1 V
0 -1 V
0 -1 V
1 0 V
0 -1 V
1 -1 V
0 -1 V
1 -1 V
0 -1 V
1 -1 V
0 -1 V
1 -1 V
0 -1 V
1 -1 V
0 -1 V
1 -1 V
0 -1 V
1 -1 V
1 -1 V
0 -1 V
1 -1 V
1 -1 V
0 -1 V
1 -1 V
1 0 V
0 -1 V
1 -1 V
0 -1 V
1 0 V
0 -1 V
1 -1 V
1 -1 V
0 -1 V
1 0 V
0 -1 V
1 -1 V
1 -1 V
1 -1 V
1 -1 V
1 -1 V
1 -1 V
0 -1 V
1 0 V
1 -1 V
1 -1 V
1 -1 V
1 -1 V
1 -1 V
1 -1 V
1 -1 V
1 0 V
1 -1 V
1 -1 V
1 -1 V
1 -1 V
1 -1 V
2 0 V
1 -1 V
1 -1 V
1 -1 V
1 0 V
1 -1 V
2 -1 V
1 -1 V
2 -1 V
1 -1 V
2 0 V
1 -1 V
2 -1 V
2 -1 V
1 -1 V
2 0 V
2 -1 V
1 -1 V
2 -1 V
1 0 V
3 -1 V
2 -1 V
1 0 V
2 -1 V
2 -1 V
2 0 V
2 -1 V
2 -1 V
3 0 V
2 -1 V
2 -1 V
3 0 V
2 -1 V
3 0 V
2 -1 V
3 0 V
3 -1 V
2 0 V
3 -1 V
3 0 V
2 -1 V
3 0 V
3 -1 V
2 0 V
3 0 V
3 -1 V
3 0 V
4 0 V
3 0 V
3 -1 V
4 0 V
currentpoint stroke M
3 0 V
3 0 V
4 0 V
3 0 V
4 0 V
3 -1 V
4 0 V
3 0 V
5 0 V
4 1 V
3 0 V
4 0 V
4 0 V
4 0 V
4 0 V
4 0 V
4 1 V
4 0 V
4 0 V
5 0 V
5 1 V
5 0 V
4 1 V
4 0 V
5 0 V
4 1 V
6 1 V
4 0 V
5 1 V
5 0 V
4 1 V
5 0 V
4 1 V
6 1 V
4 0 V
5 1 V
5 1 V
5 1 V
6 0 V
4 1 V
5 1 V
6 1 V
5 1 V
6 1 V
5 1 V
5 1 V
5 1 V
6 1 V
5 0 V
5 1 V
5 1 V
6 2 V
5 0 V
5 2 V
6 1 V
6 1 V
4 1 V
6 1 V
6 1 V
5 1 V
6 1 V
5 2 V
5 1 V
6 1 V
6 1 V
4 1 V
6 2 V
6 1 V
5 1 V
6 1 V
6 2 V
4 1 V
7 1 V
5 1 V
6 2 V
6 1 V
5 1 V
6 2 V
7 1 V
5 1 V
6 2 V
5 1 V
6 2 V
5 1 V
6 1 V
6 2 V
6 1 V
6 2 V
5 1 V
6 1 V
5 1 V
6 2 V
5 1 V
7 2 V
4 1 V
7 2 V
5 1 V
7 2 V
4 1 V
7 2 V
5 1 V
6 2 V
4 1 V
6 1 V
6 2 V
5 1 V
6 2 V
6 1 V
6 2 V
4 1 V
6 1 V
6 2 V
5 1 V
7 2 V
5 1 V
6 2 V
6 2 V
4 1 V
6 2 V
6 1 V
5 1 V
5 2 V
7 2 V
5 1 V
7 2 V
5 1 V
5 2 V
7 2 V
5 1 V
4 1 V
6 2 V
6 2 V
5 1 V
5 1 V
7 2 V
4 2 V
5 1 V
6 2 V
6 2 V
6 1 V
5 2 V
7 2 V
6 1 V
5 2 V
4 1 V
6 2 V
7 2 V
5 2 V
5 1 V
5 2 V
6 2 V
6 1 V
5 2 V
5 1 V
5 2 V
5 1 V
7 3 V
6 2 V
5 1 V
4 2 V
5 1 V
7 2 V
6 2 V
5 2 V
5 1 V
3 2 V
5 1 V
4 2 V
4 1 V
4 1 V
7 2 V
6 3 V
7 2 V
6 2 V
4 2 V
5 1 V
4 2 V
5 2 V
4 1 V
5 2 V
4 1 V
4 2 V
7 2 V
7 3 V
6 2 V
6 3 V
8 2 V
5 3 V
6 2 V
6 2 V
5 3 V
5 1 V
11 5 V
24 11 V
17 7 V
15 8 V
11 5 V
12 6 V
12 7 V
14 8 V
9 6 V
9 6 V
13 9 V
8 5 V
9 8 V
7 6 V
6 6 V
8 8 V
5 6 V
6 8 V
5 7 V
4 7 V
3 6 V
2 7 V
2 7 V
2 8 V
0 7 V
stroke
grestore
end
showpage
}}%
\put(1597,1398){\rjust{$\tau (t=0) = 0.605$}}%
\put(-69,819){\ljust{$\alpha_1$}}%
\put(1130,50){\cjust{$\tau$}}%
\put(2010,200){\cjust{ 5.5}}%
\put(1850,200){\cjust{ 5}}%
\put(1690,200){\cjust{ 4.5}}%
\put(1530,200){\cjust{ 4}}%
\put(1370,200){\cjust{ 3.5}}%
\put(1210,200){\cjust{ 3}}%
\put(1050,200){\cjust{ 2.5}}%
\put(890,200){\cjust{ 2}}%
\put(730,200){\cjust{ 1.5}}%
\put(570,200){\cjust{ 1}}%
\put(410,200){\cjust{ 0.5}}%
\put(250,200){\cjust{ 0}}%
\put(200,1511){\rjust{ 4}}%
\put(200,1338){\rjust{ 3}}%
\put(200,1165){\rjust{ 2}}%
\put(200,992){\rjust{ 1}}%
\put(200,819){\rjust{ 0}}%
\put(200,646){\rjust{-1}}%
\put(200,473){\rjust{-2}}%
\put(200,300){\rjust{-3}}%
\endGNUPLOTpicture
\endgroup
 

%% file: fig_2therm_Q1_angle.tex
\begingroup
  \catcode`\@=11\relax
  \def\GNUPLOTspecial{%
    \def\do##1{\catcode`##1=12\relax}\dospecials
    \catcode`\{=1\catcode`\}=2\catcode\%=14\relax\special}%
\expandafter\ifx\csname GNUPLOTpicture\endcsname\relax
  \csname newdimen\endcsname\GNUPLOTunit
  \gdef\GNUPLOTpicture(#1,#2){\vbox to#2\GNUPLOTunit\bgroup
    \def\put(##1,##2)##3{\unskip\raise##2\GNUPLOTunit
      \hbox to0pt{\kern##1\GNUPLOTunit ##3\hss}\ignorespaces}%
    \def\ljust##1{\vbox to0pt{\vss\hbox to0pt{##1\hss}\vss}}%
    \def\cjust##1{\vbox to0pt{\vss\hbox to0pt{\hss ##1\hss}\vss}}%
    \def\rjust##1{\vbox to0pt{\vss\hbox to0pt{\hss ##1}\vss}}%
    \def\stack##1{\let\\=\cr\tabskip=0pt\halign{\hfil ####\hfil\cr ##1\crcr}}%
    \def\lstack##1{\hbox to0pt{\vbox to0pt{\vss\stack{##1}}\hss}}%
    \def\cstack##1{\hbox to0pt{\hss\vbox to0pt{\vss\stack{##1}}\hss}}%
    \def\rstack##1{\hbox to0pt{\vbox to0pt{\stack{##1}\vss}\hss}}%
    \vss\hbox to#1\GNUPLOTunit\bgroup\ignorespaces}%
  \gdef\endGNUPLOTpicture{\hss\egroup\egroup}%
\fi
\GNUPLOTunit=0.1bp
{\GNUPLOTspecial{!
/gnudict 256 dict def
gnudict begin
/Color false def
/Solid false def
/gnulinewidth 5.000 def
/userlinewidth gnulinewidth def
/vshift -33 def
/dl {10 mul} def
/hpt_ 31.5 def
/vpt_ 31.5 def
/hpt hpt_ def
/vpt vpt_ def
/M {moveto} bind def
/L {lineto} bind def
/R {rmoveto} bind def
/V {rlineto} bind def
/vpt2 vpt 2 mul def
/hpt2 hpt 2 mul def
/Lshow { currentpoint stroke M
  0 vshift R show } def
/Rshow { currentpoint stroke M
  dup stringwidth pop neg vshift R show } def
/Cshow { currentpoint stroke M
  dup stringwidth pop -2 div vshift R show } def
/UP { dup vpt_ mul /vpt exch def hpt_ mul /hpt exch def
  /hpt2 hpt 2 mul def /vpt2 vpt 2 mul def } def
/DL { Color {setrgbcolor Solid {pop []} if 0 setdash }
 {pop pop pop Solid {pop []} if 0 setdash} ifelse } def
/BL { stroke userlinewidth 2 mul setlinewidth } def
/AL { stroke userlinewidth 2 div setlinewidth } def
/UL { dup gnulinewidth mul /userlinewidth exch def
      dup 1 lt {pop 1} if 10 mul /udl exch def } def
/PL { stroke userlinewidth setlinewidth } def
/LTb { BL [] 0 0 0 DL } def
/LTa { AL [1 udl mul 2 udl mul] 0 setdash 0 0 0 setrgbcolor } def
/LT0 { PL [] 1 0 0 DL } def
/LT1 { PL [4 dl 2 dl] 0 1 0 DL } def
/LT2 { PL [2 dl 3 dl] 0 0 1 DL } def
/LT3 { PL [1 dl 1.5 dl] 1 0 1 DL } def
/LT4 { PL [5 dl 2 dl 1 dl 2 dl] 0 1 1 DL } def
/LT5 { PL [4 dl 3 dl 1 dl 3 dl] 1 1 0 DL } def
/LT6 { PL [2 dl 2 dl 2 dl 4 dl] 0 0 0 DL } def
/LT7 { PL [2 dl 2 dl 2 dl 2 dl 2 dl 4 dl] 1 0.3 0 DL } def
/LT8 { PL [2 dl 2 dl 2 dl 2 dl 2 dl 2 dl 2 dl 4 dl] 0.5 0.5 0.5 DL } def
/Pnt { stroke [] 0 setdash
   gsave 1 setlinecap M 0 0 V stroke grestore } def
/Dia { stroke [] 0 setdash 2 copy vpt add M
  hpt neg vpt neg V hpt vpt neg V
  hpt vpt V hpt neg vpt V closepath stroke
  Pnt } def
/Pls { stroke [] 0 setdash vpt sub M 0 vpt2 V
  currentpoint stroke M
  hpt neg vpt neg R hpt2 0 V stroke
  } def
/Box { stroke [] 0 setdash 2 copy exch hpt sub exch vpt add M
  0 vpt2 neg V hpt2 0 V 0 vpt2 V
  hpt2 neg 0 V closepath stroke
  Pnt } def
/Crs { stroke [] 0 setdash exch hpt sub exch vpt add M
  hpt2 vpt2 neg V currentpoint stroke M
  hpt2 neg 0 R hpt2 vpt2 V stroke } def
/TriU { stroke [] 0 setdash 2 copy vpt 1.12 mul add M
  hpt neg vpt -1.62 mul V
  hpt 2 mul 0 V
  hpt neg vpt 1.62 mul V closepath stroke
  Pnt  } def
/Star { 2 copy Pls Crs } def
/BoxF { stroke [] 0 setdash exch hpt sub exch vpt add M
  0 vpt2 neg V  hpt2 0 V  0 vpt2 V
  hpt2 neg 0 V  closepath fill } def
/TriUF { stroke [] 0 setdash vpt 1.12 mul add M
  hpt neg vpt -1.62 mul V
  hpt 2 mul 0 V
  hpt neg vpt 1.62 mul V closepath fill } def
/TriD { stroke [] 0 setdash 2 copy vpt 1.12 mul sub M
  hpt neg vpt 1.62 mul V
  hpt 2 mul 0 V
  hpt neg vpt -1.62 mul V closepath stroke
  Pnt  } def
/TriDF { stroke [] 0 setdash vpt 1.12 mul sub M
  hpt neg vpt 1.62 mul V
  hpt 2 mul 0 V
  hpt neg vpt -1.62 mul V closepath fill} def
/DiaF { stroke [] 0 setdash vpt add M
  hpt neg vpt neg V hpt vpt neg V
  hpt vpt V hpt neg vpt V closepath fill } def
/Pent { stroke [] 0 setdash 2 copy gsave
  translate 0 hpt M 4 {72 rotate 0 hpt L} repeat
  closepath stroke grestore Pnt } def
/PentF { stroke [] 0 setdash gsave
  translate 0 hpt M 4 {72 rotate 0 hpt L} repeat
  closepath fill grestore } def
/Circle { stroke [] 0 setdash 2 copy
  hpt 0 360 arc stroke Pnt } def
/CircleF { stroke [] 0 setdash hpt 0 360 arc fill } def
/C0 { BL [] 0 setdash 2 copy moveto vpt 90 450  arc } bind def
/C1 { BL [] 0 setdash 2 copy        moveto
       2 copy  vpt 0 90 arc closepath fill
               vpt 0 360 arc closepath } bind def
/C2 { BL [] 0 setdash 2 copy moveto
       2 copy  vpt 90 180 arc closepath fill
               vpt 0 360 arc closepath } bind def
/C3 { BL [] 0 setdash 2 copy moveto
       2 copy  vpt 0 180 arc closepath fill
               vpt 0 360 arc closepath } bind def
/C4 { BL [] 0 setdash 2 copy moveto
       2 copy  vpt 180 270 arc closepath fill
               vpt 0 360 arc closepath } bind def
/C5 { BL [] 0 setdash 2 copy moveto
       2 copy  vpt 0 90 arc
       2 copy moveto
       2 copy  vpt 180 270 arc closepath fill
               vpt 0 360 arc } bind def
/C6 { BL [] 0 setdash 2 copy moveto
      2 copy  vpt 90 270 arc closepath fill
              vpt 0 360 arc closepath } bind def
/C7 { BL [] 0 setdash 2 copy moveto
      2 copy  vpt 0 270 arc closepath fill
              vpt 0 360 arc closepath } bind def
/C8 { BL [] 0 setdash 2 copy moveto
      2 copy vpt 270 360 arc closepath fill
              vpt 0 360 arc closepath } bind def
/C9 { BL [] 0 setdash 2 copy moveto
      2 copy  vpt 270 450 arc closepath fill
              vpt 0 360 arc closepath } bind def
/C10 { BL [] 0 setdash 2 copy 2 copy moveto vpt 270 360 arc closepath fill
       2 copy moveto
       2 copy vpt 90 180 arc closepath fill
               vpt 0 360 arc closepath } bind def
/C11 { BL [] 0 setdash 2 copy moveto
       2 copy  vpt 0 180 arc closepath fill
       2 copy moveto
       2 copy  vpt 270 360 arc closepath fill
               vpt 0 360 arc closepath } bind def
/C12 { BL [] 0 setdash 2 copy moveto
       2 copy  vpt 180 360 arc closepath fill
               vpt 0 360 arc closepath } bind def
/C13 { BL [] 0 setdash  2 copy moveto
       2 copy  vpt 0 90 arc closepath fill
       2 copy moveto
       2 copy  vpt 180 360 arc closepath fill
               vpt 0 360 arc closepath } bind def
/C14 { BL [] 0 setdash 2 copy moveto
       2 copy  vpt 90 360 arc closepath fill
               vpt 0 360 arc } bind def
/C15 { BL [] 0 setdash 2 copy vpt 0 360 arc closepath fill
               vpt 0 360 arc closepath } bind def
/Rec   { newpath 4 2 roll moveto 1 index 0 rlineto 0 exch rlineto
       neg 0 rlineto closepath } bind def
/Square { dup Rec } bind def
/Bsquare { vpt sub exch vpt sub exch vpt2 Square } bind def
/S0 { BL [] 0 setdash 2 copy moveto 0 vpt rlineto BL Bsquare } bind def
/S1 { BL [] 0 setdash 2 copy vpt Square fill Bsquare } bind def
/S2 { BL [] 0 setdash 2 copy exch vpt sub exch vpt Square fill Bsquare } bind def
/S3 { BL [] 0 setdash 2 copy exch vpt sub exch vpt2 vpt Rec fill Bsquare } bind def
/S4 { BL [] 0 setdash 2 copy exch vpt sub exch vpt sub vpt Square fill Bsquare } bind def
/S5 { BL [] 0 setdash 2 copy 2 copy vpt Square fill
       exch vpt sub exch vpt sub vpt Square fill Bsquare } bind def
/S6 { BL [] 0 setdash 2 copy exch vpt sub exch vpt sub vpt vpt2 Rec fill Bsquare } bind def
/S7 { BL [] 0 setdash 2 copy exch vpt sub exch vpt sub vpt vpt2 Rec fill
       2 copy vpt Square fill
       Bsquare } bind def
/S8 { BL [] 0 setdash 2 copy vpt sub vpt Square fill Bsquare } bind def
/S9 { BL [] 0 setdash 2 copy vpt sub vpt vpt2 Rec fill Bsquare } bind def
/S10 { BL [] 0 setdash 2 copy vpt sub vpt Square fill 2 copy exch vpt sub exch vpt Square fill
       Bsquare } bind def
/S11 { BL [] 0 setdash 2 copy vpt sub vpt Square fill 2 copy exch vpt sub exch vpt2 vpt Rec fill
       Bsquare } bind def
/S12 { BL [] 0 setdash 2 copy exch vpt sub exch vpt sub vpt2 vpt Rec fill Bsquare } bind def
/S13 { BL [] 0 setdash 2 copy exch vpt sub exch vpt sub vpt2 vpt Rec fill
       2 copy vpt Square fill Bsquare } bind def
/S14 { BL [] 0 setdash 2 copy exch vpt sub exch vpt sub vpt2 vpt Rec fill
       2 copy exch vpt sub exch vpt Square fill Bsquare } bind def
/S15 { BL [] 0 setdash 2 copy Bsquare fill Bsquare } bind def
/D0 { gsave translate 45 rotate 0 0 S0 stroke grestore } bind def
/D1 { gsave translate 45 rotate 0 0 S1 stroke grestore } bind def
/D2 { gsave translate 45 rotate 0 0 S2 stroke grestore } bind def
/D3 { gsave translate 45 rotate 0 0 S3 stroke grestore } bind def
/D4 { gsave translate 45 rotate 0 0 S4 stroke grestore } bind def
/D5 { gsave translate 45 rotate 0 0 S5 stroke grestore } bind def
/D6 { gsave translate 45 rotate 0 0 S6 stroke grestore } bind def
/D7 { gsave translate 45 rotate 0 0 S7 stroke grestore } bind def
/D8 { gsave translate 45 rotate 0 0 S8 stroke grestore } bind def
/D9 { gsave translate 45 rotate 0 0 S9 stroke grestore } bind def
/D10 { gsave translate 45 rotate 0 0 S10 stroke grestore } bind def
/D11 { gsave translate 45 rotate 0 0 S11 stroke grestore } bind def
/D12 { gsave translate 45 rotate 0 0 S12 stroke grestore } bind def
/D13 { gsave translate 45 rotate 0 0 S13 stroke grestore } bind def
/D14 { gsave translate 45 rotate 0 0 S14 stroke grestore } bind def
/D15 { gsave translate 45 rotate 0 0 S15 stroke grestore } bind def
/DiaE { stroke [] 0 setdash vpt add M
  hpt neg vpt neg V hpt vpt neg V
  hpt vpt V hpt neg vpt V closepath stroke } def
/BoxE { stroke [] 0 setdash exch hpt sub exch vpt add M
  0 vpt2 neg V hpt2 0 V 0 vpt2 V
  hpt2 neg 0 V closepath stroke } def
/TriUE { stroke [] 0 setdash vpt 1.12 mul add M
  hpt neg vpt -1.62 mul V
  hpt 2 mul 0 V
  hpt neg vpt 1.62 mul V closepath stroke } def
/TriDE { stroke [] 0 setdash vpt 1.12 mul sub M
  hpt neg vpt 1.62 mul V
  hpt 2 mul 0 V
  hpt neg vpt -1.62 mul V closepath stroke } def
/PentE { stroke [] 0 setdash gsave
  translate 0 hpt M 4 {72 rotate 0 hpt L} repeat
  closepath stroke grestore } def
/CircE { stroke [] 0 setdash 
  hpt 0 360 arc stroke } def
/Opaque { gsave closepath 1 setgray fill grestore 0 setgray closepath } def
/DiaW { stroke [] 0 setdash vpt add M
  hpt neg vpt neg V hpt vpt neg V
  hpt vpt V hpt neg vpt V Opaque stroke } def
/BoxW { stroke [] 0 setdash exch hpt sub exch vpt add M
  0 vpt2 neg V hpt2 0 V 0 vpt2 V
  hpt2 neg 0 V Opaque stroke } def
/TriUW { stroke [] 0 setdash vpt 1.12 mul add M
  hpt neg vpt -1.62 mul V
  hpt 2 mul 0 V
  hpt neg vpt 1.62 mul V Opaque stroke } def
/TriDW { stroke [] 0 setdash vpt 1.12 mul sub M
  hpt neg vpt 1.62 mul V
  hpt 2 mul 0 V
  hpt neg vpt -1.62 mul V Opaque stroke } def
/PentW { stroke [] 0 setdash gsave
  translate 0 hpt M 4 {72 rotate 0 hpt L} repeat
  Opaque stroke grestore } def
/CircW { stroke [] 0 setdash 
  hpt 0 360 arc Opaque stroke } def
/BoxFill { gsave Rec 1 setgray fill grestore } def
/Symbol-Oblique /Symbol findfont [1 0 .167 1 0 0] makefont
dup length dict begin {1 index /FID eq {pop pop} {def} ifelse} forall
currentdict end definefont
end
}}%
\GNUPLOTpicture(2160,1511)
{\GNUPLOTspecial{"
gnudict begin
gsave
0 0 translate
0.100 0.100 scale
0 setgray
newpath
1.000 UL
LTb
500 300 M
63 0 V
1447 0 R
-63 0 V
500 563 M
63 0 V
1447 0 R
-63 0 V
500 827 M
63 0 V
1447 0 R
-63 0 V
500 1090 M
63 0 V
1447 0 R
-63 0 V
500 1353 M
63 0 V
1447 0 R
-63 0 V
500 300 M
0 63 V
0 1148 R
0 -63 V
740 300 M
0 63 V
0 1148 R
0 -63 V
981 300 M
0 63 V
0 1148 R
0 -63 V
1221 300 M
0 63 V
0 1148 R
0 -63 V
1461 300 M
0 63 V
0 1148 R
0 -63 V
1702 300 M
0 63 V
0 1148 R
0 -63 V
1942 300 M
0 63 V
0 1148 R
0 -63 V
1.000 UL
LTb
500 300 M
1510 0 V
0 1211 V
-1510 0 V
500 300 L
1.000 UL
LT0
1647 1398 M
263 0 V
508 930 M
15 -2 V
15 -20 V
15 -15 V
15 -20 V
15 -18 V
15 -36 V
15 15 V
628 694 L
643 548 L
16 171 V
15 -50 V
15 -99 V
15 157 V
15 -73 V
15 63 V
15 80 V
15 38 V
15 23 V
15 11 V
16 7 V
15 -8 V
15 16 V
15 17 V
15 17 V
15 30 V
15 27 V
15 12 V
15 21 V
15 9 V
16 -27 V
15 -9 V
15 -42 V
15 -71 V
15 84 V
15 -26 V
15 32 V
15 75 V
15 20 V
15 69 V
16 -154 V
15 5 V
15 -38 V
15 -24 V
15 7 V
15 -9 V
15 -1 V
15 8 V
15 4 V
15 14 V
16 5 V
15 3 V
15 7 V
15 9 V
15 3 V
15 1 V
15 24 V
15 14 V
15 -198 V
15 98 V
16 167 V
15 120 V
15 -77 V
15 -38 V
15 -11 V
15 -9 V
15 100 V
15 -95 V
15 -22 V
15 -34 V
16 21 V
15 -14 V
15 -27 V
15 -22 V
15 -15 V
15 -19 V
15 -9 V
15 -14 V
15 -21 V
15 -31 V
16 -15 V
15 -61 V
15 24 V
15 105 V
15 -108 V
15 -40 V
15 1 V
15 -145 V
15 66 V
15 114 V
16 83 V
15 102 V
15 29 V
15 -88 V
15 -15 V
15 25 V
15 9 V
15 -4 V
15 6 V
15 0 V
1.000 UL
LT1
1647 1298 M
263 0 V
508 908 M
15 0 V
15 0 V
15 0 V
15 0 V
15 0 V
15 0 V
15 0 V
15 0 V
15 0 V
16 0 V
15 0 V
15 0 V
15 0 V
15 0 V
15 0 V
15 0 V
15 0 V
15 0 V
15 0 V
16 0 V
15 0 V
15 0 V
15 0 V
15 0 V
15 0 V
15 0 V
15 0 V
15 0 V
15 0 V
16 0 V
15 0 V
15 0 V
15 0 V
15 0 V
15 0 V
15 0 V
15 0 V
15 0 V
15 0 V
16 0 V
15 0 V
15 0 V
15 0 V
15 0 V
15 0 V
15 0 V
15 0 V
15 0 V
15 0 V
16 0 V
15 0 V
15 0 V
15 0 V
15 0 V
15 0 V
15 0 V
15 0 V
15 0 V
15 0 V
16 0 V
15 0 V
15 0 V
15 0 V
15 0 V
15 0 V
15 0 V
15 0 V
15 0 V
15 0 V
16 0 V
15 0 V
15 0 V
15 0 V
15 0 V
15 0 V
15 0 V
15 0 V
15 0 V
15 0 V
16 0 V
15 0 V
15 0 V
15 0 V
15 0 V
15 0 V
15 0 V
15 0 V
15 0 V
15 0 V
16 0 V
15 0 V
15 0 V
15 0 V
15 0 V
15 0 V
15 0 V
15 0 V
15 0 V
15 0 V
stroke
grestore
end
showpage
}}%
\put(1597,1298){\rjust{theoretical}}%
\put(1597,1398){\rjust{numerical}}%
\put(1255,50){\cjust{$\theta$}}%
\put(1942,200){\cjust{ 6}}%
\put(1702,200){\cjust{ 5}}%
\put(1461,200){\cjust{ 4}}%
\put(1221,200){\cjust{ 3}}%
\put(981,200){\cjust{ 2}}%
\put(740,200){\cjust{ 1}}%
\put(500,200){\cjust{ 0}}%
\put(450,1353){\rjust{ 0.16}}%
\put(450,1090){\rjust{ 0.1595}}%
\put(450,827){\rjust{ 0.159}}%
\put(450,563){\rjust{ 0.1585}}%
\put(450,300){\rjust{ 0.158}}%
\endGNUPLOTpicture
\endgroup
 

%% file: fig_2therm_Q1_amplitude.tex
\begingroup
  \catcode`\@=11\relax
  \def\GNUPLOTspecial{%
    \def\do##1{\catcode`##1=12\relax}\dospecials
    \catcode`\{=1\catcode`\}=2\catcode\%=14\relax\special}%
\expandafter\ifx\csname GNUPLOTpicture\endcsname\relax
  \csname newdimen\endcsname\GNUPLOTunit
  \gdef\GNUPLOTpicture(#1,#2){\vbox to#2\GNUPLOTunit\bgroup
    \def\put(##1,##2)##3{\unskip\raise##2\GNUPLOTunit
      \hbox to0pt{\kern##1\GNUPLOTunit ##3\hss}\ignorespaces}%
    \def\ljust##1{\vbox to0pt{\vss\hbox to0pt{##1\hss}\vss}}%
    \def\cjust##1{\vbox to0pt{\vss\hbox to0pt{\hss ##1\hss}\vss}}%
    \def\rjust##1{\vbox to0pt{\vss\hbox to0pt{\hss ##1}\vss}}%
    \def\stack##1{\let\\=\cr\tabskip=0pt\halign{\hfil ####\hfil\cr ##1\crcr}}%
    \def\lstack##1{\hbox to0pt{\vbox to0pt{\vss\stack{##1}}\hss}}%
    \def\cstack##1{\hbox to0pt{\hss\vbox to0pt{\vss\stack{##1}}\hss}}%
    \def\rstack##1{\hbox to0pt{\vbox to0pt{\stack{##1}\vss}\hss}}%
    \vss\hbox to#1\GNUPLOTunit\bgroup\ignorespaces}%
  \gdef\endGNUPLOTpicture{\hss\egroup\egroup}%
\fi
\GNUPLOTunit=0.1bp
{\GNUPLOTspecial{!
/gnudict 256 dict def
gnudict begin
/Color false def
/Solid false def
/gnulinewidth 5.000 def
/userlinewidth gnulinewidth def
/vshift -33 def
/dl {10 mul} def
/hpt_ 31.5 def
/vpt_ 31.5 def
/hpt hpt_ def
/vpt vpt_ def
/M {moveto} bind def
/L {lineto} bind def
/R {rmoveto} bind def
/V {rlineto} bind def
/vpt2 vpt 2 mul def
/hpt2 hpt 2 mul def
/Lshow { currentpoint stroke M
  0 vshift R show } def
/Rshow { currentpoint stroke M
  dup stringwidth pop neg vshift R show } def
/Cshow { currentpoint stroke M
  dup stringwidth pop -2 div vshift R show } def
/UP { dup vpt_ mul /vpt exch def hpt_ mul /hpt exch def
  /hpt2 hpt 2 mul def /vpt2 vpt 2 mul def } def
/DL { Color {setrgbcolor Solid {pop []} if 0 setdash }
 {pop pop pop Solid {pop []} if 0 setdash} ifelse } def
/BL { stroke userlinewidth 2 mul setlinewidth } def
/AL { stroke userlinewidth 2 div setlinewidth } def
/UL { dup gnulinewidth mul /userlinewidth exch def
      dup 1 lt {pop 1} if 10 mul /udl exch def } def
/PL { stroke userlinewidth setlinewidth } def
/LTb { BL [] 0 0 0 DL } def
/LTa { AL [1 udl mul 2 udl mul] 0 setdash 0 0 0 setrgbcolor } def
/LT0 { PL [] 1 0 0 DL } def
/LT1 { PL [4 dl 2 dl] 0 1 0 DL } def
/LT2 { PL [2 dl 3 dl] 0 0 1 DL } def
/LT3 { PL [1 dl 1.5 dl] 1 0 1 DL } def
/LT4 { PL [5 dl 2 dl 1 dl 2 dl] 0 1 1 DL } def
/LT5 { PL [4 dl 3 dl 1 dl 3 dl] 1 1 0 DL } def
/LT6 { PL [2 dl 2 dl 2 dl 4 dl] 0 0 0 DL } def
/LT7 { PL [2 dl 2 dl 2 dl 2 dl 2 dl 4 dl] 1 0.3 0 DL } def
/LT8 { PL [2 dl 2 dl 2 dl 2 dl 2 dl 2 dl 2 dl 4 dl] 0.5 0.5 0.5 DL } def
/Pnt { stroke [] 0 setdash
   gsave 1 setlinecap M 0 0 V stroke grestore } def
/Dia { stroke [] 0 setdash 2 copy vpt add M
  hpt neg vpt neg V hpt vpt neg V
  hpt vpt V hpt neg vpt V closepath stroke
  Pnt } def
/Pls { stroke [] 0 setdash vpt sub M 0 vpt2 V
  currentpoint stroke M
  hpt neg vpt neg R hpt2 0 V stroke
  } def
/Box { stroke [] 0 setdash 2 copy exch hpt sub exch vpt add M
  0 vpt2 neg V hpt2 0 V 0 vpt2 V
  hpt2 neg 0 V closepath stroke
  Pnt } def
/Crs { stroke [] 0 setdash exch hpt sub exch vpt add M
  hpt2 vpt2 neg V currentpoint stroke M
  hpt2 neg 0 R hpt2 vpt2 V stroke } def
/TriU { stroke [] 0 setdash 2 copy vpt 1.12 mul add M
  hpt neg vpt -1.62 mul V
  hpt 2 mul 0 V
  hpt neg vpt 1.62 mul V closepath stroke
  Pnt  } def
/Star { 2 copy Pls Crs } def
/BoxF { stroke [] 0 setdash exch hpt sub exch vpt add M
  0 vpt2 neg V  hpt2 0 V  0 vpt2 V
  hpt2 neg 0 V  closepath fill } def
/TriUF { stroke [] 0 setdash vpt 1.12 mul add M
  hpt neg vpt -1.62 mul V
  hpt 2 mul 0 V
  hpt neg vpt 1.62 mul V closepath fill } def
/TriD { stroke [] 0 setdash 2 copy vpt 1.12 mul sub M
  hpt neg vpt 1.62 mul V
  hpt 2 mul 0 V
  hpt neg vpt -1.62 mul V closepath stroke
  Pnt  } def
/TriDF { stroke [] 0 setdash vpt 1.12 mul sub M
  hpt neg vpt 1.62 mul V
  hpt 2 mul 0 V
  hpt neg vpt -1.62 mul V closepath fill} def
/DiaF { stroke [] 0 setdash vpt add M
  hpt neg vpt neg V hpt vpt neg V
  hpt vpt V hpt neg vpt V closepath fill } def
/Pent { stroke [] 0 setdash 2 copy gsave
  translate 0 hpt M 4 {72 rotate 0 hpt L} repeat
  closepath stroke grestore Pnt } def
/PentF { stroke [] 0 setdash gsave
  translate 0 hpt M 4 {72 rotate 0 hpt L} repeat
  closepath fill grestore } def
/Circle { stroke [] 0 setdash 2 copy
  hpt 0 360 arc stroke Pnt } def
/CircleF { stroke [] 0 setdash hpt 0 360 arc fill } def
/C0 { BL [] 0 setdash 2 copy moveto vpt 90 450  arc } bind def
/C1 { BL [] 0 setdash 2 copy        moveto
       2 copy  vpt 0 90 arc closepath fill
               vpt 0 360 arc closepath } bind def
/C2 { BL [] 0 setdash 2 copy moveto
       2 copy  vpt 90 180 arc closepath fill
               vpt 0 360 arc closepath } bind def
/C3 { BL [] 0 setdash 2 copy moveto
       2 copy  vpt 0 180 arc closepath fill
               vpt 0 360 arc closepath } bind def
/C4 { BL [] 0 setdash 2 copy moveto
       2 copy  vpt 180 270 arc closepath fill
               vpt 0 360 arc closepath } bind def
/C5 { BL [] 0 setdash 2 copy moveto
       2 copy  vpt 0 90 arc
       2 copy moveto
       2 copy  vpt 180 270 arc closepath fill
               vpt 0 360 arc } bind def
/C6 { BL [] 0 setdash 2 copy moveto
      2 copy  vpt 90 270 arc closepath fill
              vpt 0 360 arc closepath } bind def
/C7 { BL [] 0 setdash 2 copy moveto
      2 copy  vpt 0 270 arc closepath fill
              vpt 0 360 arc closepath } bind def
/C8 { BL [] 0 setdash 2 copy moveto
      2 copy vpt 270 360 arc closepath fill
              vpt 0 360 arc closepath } bind def
/C9 { BL [] 0 setdash 2 copy moveto
      2 copy  vpt 270 450 arc closepath fill
              vpt 0 360 arc closepath } bind def
/C10 { BL [] 0 setdash 2 copy 2 copy moveto vpt 270 360 arc closepath fill
       2 copy moveto
       2 copy vpt 90 180 arc closepath fill
               vpt 0 360 arc closepath } bind def
/C11 { BL [] 0 setdash 2 copy moveto
       2 copy  vpt 0 180 arc closepath fill
       2 copy moveto
       2 copy  vpt 270 360 arc closepath fill
               vpt 0 360 arc closepath } bind def
/C12 { BL [] 0 setdash 2 copy moveto
       2 copy  vpt 180 360 arc closepath fill
               vpt 0 360 arc closepath } bind def
/C13 { BL [] 0 setdash  2 copy moveto
       2 copy  vpt 0 90 arc closepath fill
       2 copy moveto
       2 copy  vpt 180 360 arc closepath fill
               vpt 0 360 arc closepath } bind def
/C14 { BL [] 0 setdash 2 copy moveto
       2 copy  vpt 90 360 arc closepath fill
               vpt 0 360 arc } bind def
/C15 { BL [] 0 setdash 2 copy vpt 0 360 arc closepath fill
               vpt 0 360 arc closepath } bind def
/Rec   { newpath 4 2 roll moveto 1 index 0 rlineto 0 exch rlineto
       neg 0 rlineto closepath } bind def
/Square { dup Rec } bind def
/Bsquare { vpt sub exch vpt sub exch vpt2 Square } bind def
/S0 { BL [] 0 setdash 2 copy moveto 0 vpt rlineto BL Bsquare } bind def
/S1 { BL [] 0 setdash 2 copy vpt Square fill Bsquare } bind def
/S2 { BL [] 0 setdash 2 copy exch vpt sub exch vpt Square fill Bsquare } bind def
/S3 { BL [] 0 setdash 2 copy exch vpt sub exch vpt2 vpt Rec fill Bsquare } bind def
/S4 { BL [] 0 setdash 2 copy exch vpt sub exch vpt sub vpt Square fill Bsquare } bind def
/S5 { BL [] 0 setdash 2 copy 2 copy vpt Square fill
       exch vpt sub exch vpt sub vpt Square fill Bsquare } bind def
/S6 { BL [] 0 setdash 2 copy exch vpt sub exch vpt sub vpt vpt2 Rec fill Bsquare } bind def
/S7 { BL [] 0 setdash 2 copy exch vpt sub exch vpt sub vpt vpt2 Rec fill
       2 copy vpt Square fill
       Bsquare } bind def
/S8 { BL [] 0 setdash 2 copy vpt sub vpt Square fill Bsquare } bind def
/S9 { BL [] 0 setdash 2 copy vpt sub vpt vpt2 Rec fill Bsquare } bind def
/S10 { BL [] 0 setdash 2 copy vpt sub vpt Square fill 2 copy exch vpt sub exch vpt Square fill
       Bsquare } bind def
/S11 { BL [] 0 setdash 2 copy vpt sub vpt Square fill 2 copy exch vpt sub exch vpt2 vpt Rec fill
       Bsquare } bind def
/S12 { BL [] 0 setdash 2 copy exch vpt sub exch vpt sub vpt2 vpt Rec fill Bsquare } bind def
/S13 { BL [] 0 setdash 2 copy exch vpt sub exch vpt sub vpt2 vpt Rec fill
       2 copy vpt Square fill Bsquare } bind def
/S14 { BL [] 0 setdash 2 copy exch vpt sub exch vpt sub vpt2 vpt Rec fill
       2 copy exch vpt sub exch vpt Square fill Bsquare } bind def
/S15 { BL [] 0 setdash 2 copy Bsquare fill Bsquare } bind def
/D0 { gsave translate 45 rotate 0 0 S0 stroke grestore } bind def
/D1 { gsave translate 45 rotate 0 0 S1 stroke grestore } bind def
/D2 { gsave translate 45 rotate 0 0 S2 stroke grestore } bind def
/D3 { gsave translate 45 rotate 0 0 S3 stroke grestore } bind def
/D4 { gsave translate 45 rotate 0 0 S4 stroke grestore } bind def
/D5 { gsave translate 45 rotate 0 0 S5 stroke grestore } bind def
/D6 { gsave translate 45 rotate 0 0 S6 stroke grestore } bind def
/D7 { gsave translate 45 rotate 0 0 S7 stroke grestore } bind def
/D8 { gsave translate 45 rotate 0 0 S8 stroke grestore } bind def
/D9 { gsave translate 45 rotate 0 0 S9 stroke grestore } bind def
/D10 { gsave translate 45 rotate 0 0 S10 stroke grestore } bind def
/D11 { gsave translate 45 rotate 0 0 S11 stroke grestore } bind def
/D12 { gsave translate 45 rotate 0 0 S12 stroke grestore } bind def
/D13 { gsave translate 45 rotate 0 0 S13 stroke grestore } bind def
/D14 { gsave translate 45 rotate 0 0 S14 stroke grestore } bind def
/D15 { gsave translate 45 rotate 0 0 S15 stroke grestore } bind def
/DiaE { stroke [] 0 setdash vpt add M
  hpt neg vpt neg V hpt vpt neg V
  hpt vpt V hpt neg vpt V closepath stroke } def
/BoxE { stroke [] 0 setdash exch hpt sub exch vpt add M
  0 vpt2 neg V hpt2 0 V 0 vpt2 V
  hpt2 neg 0 V closepath stroke } def
/TriUE { stroke [] 0 setdash vpt 1.12 mul add M
  hpt neg vpt -1.62 mul V
  hpt 2 mul 0 V
  hpt neg vpt 1.62 mul V closepath stroke } def
/TriDE { stroke [] 0 setdash vpt 1.12 mul sub M
  hpt neg vpt 1.62 mul V
  hpt 2 mul 0 V
  hpt neg vpt -1.62 mul V closepath stroke } def
/PentE { stroke [] 0 setdash gsave
  translate 0 hpt M 4 {72 rotate 0 hpt L} repeat
  closepath stroke grestore } def
/CircE { stroke [] 0 setdash 
  hpt 0 360 arc stroke } def
/Opaque { gsave closepath 1 setgray fill grestore 0 setgray closepath } def
/DiaW { stroke [] 0 setdash vpt add M
  hpt neg vpt neg V hpt vpt neg V
  hpt vpt V hpt neg vpt V Opaque stroke } def
/BoxW { stroke [] 0 setdash exch hpt sub exch vpt add M
  0 vpt2 neg V hpt2 0 V 0 vpt2 V
  hpt2 neg 0 V Opaque stroke } def
/TriUW { stroke [] 0 setdash vpt 1.12 mul add M
  hpt neg vpt -1.62 mul V
  hpt 2 mul 0 V
  hpt neg vpt 1.62 mul V Opaque stroke } def
/TriDW { stroke [] 0 setdash vpt 1.12 mul sub M
  hpt neg vpt 1.62 mul V
  hpt 2 mul 0 V
  hpt neg vpt -1.62 mul V Opaque stroke } def
/PentW { stroke [] 0 setdash gsave
  translate 0 hpt M 4 {72 rotate 0 hpt L} repeat
  Opaque stroke grestore } def
/CircW { stroke [] 0 setdash 
  hpt 0 360 arc Opaque stroke } def
/BoxFill { gsave Rec 1 setgray fill grestore } def
/Symbol-Oblique /Symbol findfont [1 0 .167 1 0 0] makefont
dup length dict begin {1 index /FID eq {pop pop} {def} ifelse} forall
currentdict end definefont
end
}}%
\GNUPLOTpicture(2160,1511)
{\GNUPLOTspecial{"
gnudict begin
gsave
0 0 translate
0.100 0.100 scale
0 setgray
newpath
1.000 UL
LTb
500 300 M
63 0 V
1447 0 R
-63 0 V
500 473 M
63 0 V
1447 0 R
-63 0 V
500 646 M
63 0 V
1447 0 R
-63 0 V
500 819 M
63 0 V
1447 0 R
-63 0 V
500 992 M
63 0 V
1447 0 R
-63 0 V
500 1165 M
63 0 V
1447 0 R
-63 0 V
500 1338 M
63 0 V
1447 0 R
-63 0 V
500 1511 M
63 0 V
1447 0 R
-63 0 V
500 300 M
0 63 V
0 1148 R
0 -63 V
689 300 M
0 63 V
0 1148 R
0 -63 V
878 300 M
0 63 V
0 1148 R
0 -63 V
1066 300 M
0 63 V
0 1148 R
0 -63 V
1255 300 M
0 63 V
0 1148 R
0 -63 V
1444 300 M
0 63 V
0 1148 R
0 -63 V
1633 300 M
0 63 V
0 1148 R
0 -63 V
1821 300 M
0 63 V
0 1148 R
0 -63 V
2010 300 M
0 63 V
0 1148 R
0 -63 V
1.000 UL
LTb
500 300 M
1510 0 V
0 1211 V
-1510 0 V
500 300 L
1.000 UL
LT0
1647 1398 M
263 0 V
515 453 M
15 486 V
15 186 V
15 -3 V
575 309 L
15 476 V
15 528 V
620 848 L
636 536 L
15 409 V
15 29 V
681 613 L
15 -89 V
15 167 V
15 -74 V
15 209 V
15 330 V
771 301 L
16 20 V
15 543 V
817 573 L
15 91 V
847 415 L
15 294 V
15 540 V
892 424 L
15 86 V
15 -62 V
16 99 V
15 -52 V
15 42 V
983 313 L
15 152 V
15 -151 V
15 86 V
15 101 V
15 246 V
15 -102 V
16 -193 V
15 -11 V
15 203 V
15 -88 V
15 -102 V
15 69 V
15 -8 V
15 -176 V
15 22 V
15 -28 V
16 273 V
15 160 V
15 -406 V
15 -29 V
15 199 V
15 61 V
15 -236 V
15 100 V
15 215 V
15 39 V
16 -14 V
15 -233 V
15 90 V
15 -141 V
15 57 V
15 108 V
15 9 V
15 -217 V
15 -33 V
15 21 V
16 37 V
15 42 V
15 -76 V
15 75 V
15 -43 V
15 -26 V
15 -20 V
15 18 V
15 19 V
15 -12 V
16 -60 V
15 -4 V
15 68 V
15 39 V
15 11 V
15 -84 V
15 55 V
15 -28 V
15 -20 V
15 -19 V
16 10 V
15 52 V
15 -54 V
15 -13 V
15 -14 V
15 24 V
15 -22 V
15 -1 V
15 7 V
15 23 V
16 -13 V
15 -26 V
stroke
grestore
end
showpage
}}%
\put(1597,1398){\rjust{$ \left| f_{\rm theo}^{\rm amp}(r) - f_{\rm num}^{\rm amp}(r) \right| $}}%
\put(1255,50){\cjust{$r=\sqrt{q^2+p^2}$}}%
\put(2010,200){\cjust{ 4}}%
\put(1821,200){\cjust{ 3.5}}%
\put(1633,200){\cjust{ 3}}%
\put(1444,200){\cjust{ 2.5}}%
\put(1255,200){\cjust{ 2}}%
\put(1066,200){\cjust{ 1.5}}%
\put(878,200){\cjust{ 1}}%
\put(689,200){\cjust{ 0.5}}%
\put(500,200){\cjust{ 0}}%
\put(450,1511){\rjust{ 0.0035}}%
\put(450,1338){\rjust{ 0.003}}%
\put(450,1165){\rjust{ 0.0025}}%
\put(450,992){\rjust{ 0.002}}%
\put(450,819){\rjust{ 0.0015}}%
\put(450,646){\rjust{ 0.001}}%
\put(450,473){\rjust{ 0.0005}}%
\put(450,300){\rjust{ 0}}%
\endGNUPLOTpicture
\endgroup
 

%% file: fig_2therm_Q1_disc.tex
\begingroup
  \catcode`\@=11\relax
  \def\GNUPLOTspecial{%
    \def\do##1{\catcode`##1=12\relax}\dospecials
    \catcode`\{=1\catcode`\}=2\catcode\%=14\relax\special}%
\expandafter\ifx\csname GNUPLOTpicture\endcsname\relax
  \csname newdimen\endcsname\GNUPLOTunit
  \gdef\GNUPLOTpicture(#1,#2){\vbox to#2\GNUPLOTunit\bgroup
    \def\put(##1,##2)##3{\unskip\raise##2\GNUPLOTunit
      \hbox to0pt{\kern##1\GNUPLOTunit ##3\hss}\ignorespaces}%
    \def\ljust##1{\vbox to0pt{\vss\hbox to0pt{##1\hss}\vss}}%
    \def\cjust##1{\vbox to0pt{\vss\hbox to0pt{\hss ##1\hss}\vss}}%
    \def\rjust##1{\vbox to0pt{\vss\hbox to0pt{\hss ##1}\vss}}%
    \def\stack##1{\let\\=\cr\tabskip=0pt\halign{\hfil ####\hfil\cr ##1\crcr}}%
    \def\lstack##1{\hbox to0pt{\vbox to0pt{\vss\stack{##1}}\hss}}%
    \def\cstack##1{\hbox to0pt{\hss\vbox to0pt{\vss\stack{##1}}\hss}}%
    \def\rstack##1{\hbox to0pt{\vbox to0pt{\stack{##1}\vss}\hss}}%
    \vss\hbox to#1\GNUPLOTunit\bgroup\ignorespaces}%
  \gdef\endGNUPLOTpicture{\hss\egroup\egroup}%
\fi
\GNUPLOTunit=0.1bp
{\GNUPLOTspecial{!
/gnudict 256 dict def
gnudict begin
/Color false def
/Solid false def
/gnulinewidth 5.000 def
/userlinewidth gnulinewidth def
/vshift -33 def
/dl {10 mul} def
/hpt_ 31.5 def
/vpt_ 31.5 def
/hpt hpt_ def
/vpt vpt_ def
/M {moveto} bind def
/L {lineto} bind def
/R {rmoveto} bind def
/V {rlineto} bind def
/vpt2 vpt 2 mul def
/hpt2 hpt 2 mul def
/Lshow { currentpoint stroke M
  0 vshift R show } def
/Rshow { currentpoint stroke M
  dup stringwidth pop neg vshift R show } def
/Cshow { currentpoint stroke M
  dup stringwidth pop -2 div vshift R show } def
/UP { dup vpt_ mul /vpt exch def hpt_ mul /hpt exch def
  /hpt2 hpt 2 mul def /vpt2 vpt 2 mul def } def
/DL { Color {setrgbcolor Solid {pop []} if 0 setdash }
 {pop pop pop Solid {pop []} if 0 setdash} ifelse } def
/BL { stroke userlinewidth 2 mul setlinewidth } def
/AL { stroke userlinewidth 2 div setlinewidth } def
/UL { dup gnulinewidth mul /userlinewidth exch def
      dup 1 lt {pop 1} if 10 mul /udl exch def } def
/PL { stroke userlinewidth setlinewidth } def
/LTb { BL [] 0 0 0 DL } def
/LTa { AL [1 udl mul 2 udl mul] 0 setdash 0 0 0 setrgbcolor } def
/LT0 { PL [] 1 0 0 DL } def
/LT1 { PL [4 dl 2 dl] 0 1 0 DL } def
/LT2 { PL [2 dl 3 dl] 0 0 1 DL } def
/LT3 { PL [1 dl 1.5 dl] 1 0 1 DL } def
/LT4 { PL [5 dl 2 dl 1 dl 2 dl] 0 1 1 DL } def
/LT5 { PL [4 dl 3 dl 1 dl 3 dl] 1 1 0 DL } def
/LT6 { PL [2 dl 2 dl 2 dl 4 dl] 0 0 0 DL } def
/LT7 { PL [2 dl 2 dl 2 dl 2 dl 2 dl 4 dl] 1 0.3 0 DL } def
/LT8 { PL [2 dl 2 dl 2 dl 2 dl 2 dl 2 dl 2 dl 4 dl] 0.5 0.5 0.5 DL } def
/Pnt { stroke [] 0 setdash
   gsave 1 setlinecap M 0 0 V stroke grestore } def
/Dia { stroke [] 0 setdash 2 copy vpt add M
  hpt neg vpt neg V hpt vpt neg V
  hpt vpt V hpt neg vpt V closepath stroke
  Pnt } def
/Pls { stroke [] 0 setdash vpt sub M 0 vpt2 V
  currentpoint stroke M
  hpt neg vpt neg R hpt2 0 V stroke
  } def
/Box { stroke [] 0 setdash 2 copy exch hpt sub exch vpt add M
  0 vpt2 neg V hpt2 0 V 0 vpt2 V
  hpt2 neg 0 V closepath stroke
  Pnt } def
/Crs { stroke [] 0 setdash exch hpt sub exch vpt add M
  hpt2 vpt2 neg V currentpoint stroke M
  hpt2 neg 0 R hpt2 vpt2 V stroke } def
/TriU { stroke [] 0 setdash 2 copy vpt 1.12 mul add M
  hpt neg vpt -1.62 mul V
  hpt 2 mul 0 V
  hpt neg vpt 1.62 mul V closepath stroke
  Pnt  } def
/Star { 2 copy Pls Crs } def
/BoxF { stroke [] 0 setdash exch hpt sub exch vpt add M
  0 vpt2 neg V  hpt2 0 V  0 vpt2 V
  hpt2 neg 0 V  closepath fill } def
/TriUF { stroke [] 0 setdash vpt 1.12 mul add M
  hpt neg vpt -1.62 mul V
  hpt 2 mul 0 V
  hpt neg vpt 1.62 mul V closepath fill } def
/TriD { stroke [] 0 setdash 2 copy vpt 1.12 mul sub M
  hpt neg vpt 1.62 mul V
  hpt 2 mul 0 V
  hpt neg vpt -1.62 mul V closepath stroke
  Pnt  } def
/TriDF { stroke [] 0 setdash vpt 1.12 mul sub M
  hpt neg vpt 1.62 mul V
  hpt 2 mul 0 V
  hpt neg vpt -1.62 mul V closepath fill} def
/DiaF { stroke [] 0 setdash vpt add M
  hpt neg vpt neg V hpt vpt neg V
  hpt vpt V hpt neg vpt V closepath fill } def
/Pent { stroke [] 0 setdash 2 copy gsave
  translate 0 hpt M 4 {72 rotate 0 hpt L} repeat
  closepath stroke grestore Pnt } def
/PentF { stroke [] 0 setdash gsave
  translate 0 hpt M 4 {72 rotate 0 hpt L} repeat
  closepath fill grestore } def
/Circle { stroke [] 0 setdash 2 copy
  hpt 0 360 arc stroke Pnt } def
/CircleF { stroke [] 0 setdash hpt 0 360 arc fill } def
/C0 { BL [] 0 setdash 2 copy moveto vpt 90 450  arc } bind def
/C1 { BL [] 0 setdash 2 copy        moveto
       2 copy  vpt 0 90 arc closepath fill
               vpt 0 360 arc closepath } bind def
/C2 { BL [] 0 setdash 2 copy moveto
       2 copy  vpt 90 180 arc closepath fill
               vpt 0 360 arc closepath } bind def
/C3 { BL [] 0 setdash 2 copy moveto
       2 copy  vpt 0 180 arc closepath fill
               vpt 0 360 arc closepath } bind def
/C4 { BL [] 0 setdash 2 copy moveto
       2 copy  vpt 180 270 arc closepath fill
               vpt 0 360 arc closepath } bind def
/C5 { BL [] 0 setdash 2 copy moveto
       2 copy  vpt 0 90 arc
       2 copy moveto
       2 copy  vpt 180 270 arc closepath fill
               vpt 0 360 arc } bind def
/C6 { BL [] 0 setdash 2 copy moveto
      2 copy  vpt 90 270 arc closepath fill
              vpt 0 360 arc closepath } bind def
/C7 { BL [] 0 setdash 2 copy moveto
      2 copy  vpt 0 270 arc closepath fill
              vpt 0 360 arc closepath } bind def
/C8 { BL [] 0 setdash 2 copy moveto
      2 copy vpt 270 360 arc closepath fill
              vpt 0 360 arc closepath } bind def
/C9 { BL [] 0 setdash 2 copy moveto
      2 copy  vpt 270 450 arc closepath fill
              vpt 0 360 arc closepath } bind def
/C10 { BL [] 0 setdash 2 copy 2 copy moveto vpt 270 360 arc closepath fill
       2 copy moveto
       2 copy vpt 90 180 arc closepath fill
               vpt 0 360 arc closepath } bind def
/C11 { BL [] 0 setdash 2 copy moveto
       2 copy  vpt 0 180 arc closepath fill
       2 copy moveto
       2 copy  vpt 270 360 arc closepath fill
               vpt 0 360 arc closepath } bind def
/C12 { BL [] 0 setdash 2 copy moveto
       2 copy  vpt 180 360 arc closepath fill
               vpt 0 360 arc closepath } bind def
/C13 { BL [] 0 setdash  2 copy moveto
       2 copy  vpt 0 90 arc closepath fill
       2 copy moveto
       2 copy  vpt 180 360 arc closepath fill
               vpt 0 360 arc closepath } bind def
/C14 { BL [] 0 setdash 2 copy moveto
       2 copy  vpt 90 360 arc closepath fill
               vpt 0 360 arc } bind def
/C15 { BL [] 0 setdash 2 copy vpt 0 360 arc closepath fill
               vpt 0 360 arc closepath } bind def
/Rec   { newpath 4 2 roll moveto 1 index 0 rlineto 0 exch rlineto
       neg 0 rlineto closepath } bind def
/Square { dup Rec } bind def
/Bsquare { vpt sub exch vpt sub exch vpt2 Square } bind def
/S0 { BL [] 0 setdash 2 copy moveto 0 vpt rlineto BL Bsquare } bind def
/S1 { BL [] 0 setdash 2 copy vpt Square fill Bsquare } bind def
/S2 { BL [] 0 setdash 2 copy exch vpt sub exch vpt Square fill Bsquare } bind def
/S3 { BL [] 0 setdash 2 copy exch vpt sub exch vpt2 vpt Rec fill Bsquare } bind def
/S4 { BL [] 0 setdash 2 copy exch vpt sub exch vpt sub vpt Square fill Bsquare } bind def
/S5 { BL [] 0 setdash 2 copy 2 copy vpt Square fill
       exch vpt sub exch vpt sub vpt Square fill Bsquare } bind def
/S6 { BL [] 0 setdash 2 copy exch vpt sub exch vpt sub vpt vpt2 Rec fill Bsquare } bind def
/S7 { BL [] 0 setdash 2 copy exch vpt sub exch vpt sub vpt vpt2 Rec fill
       2 copy vpt Square fill
       Bsquare } bind def
/S8 { BL [] 0 setdash 2 copy vpt sub vpt Square fill Bsquare } bind def
/S9 { BL [] 0 setdash 2 copy vpt sub vpt vpt2 Rec fill Bsquare } bind def
/S10 { BL [] 0 setdash 2 copy vpt sub vpt Square fill 2 copy exch vpt sub exch vpt Square fill
       Bsquare } bind def
/S11 { BL [] 0 setdash 2 copy vpt sub vpt Square fill 2 copy exch vpt sub exch vpt2 vpt Rec fill
       Bsquare } bind def
/S12 { BL [] 0 setdash 2 copy exch vpt sub exch vpt sub vpt2 vpt Rec fill Bsquare } bind def
/S13 { BL [] 0 setdash 2 copy exch vpt sub exch vpt sub vpt2 vpt Rec fill
       2 copy vpt Square fill Bsquare } bind def
/S14 { BL [] 0 setdash 2 copy exch vpt sub exch vpt sub vpt2 vpt Rec fill
       2 copy exch vpt sub exch vpt Square fill Bsquare } bind def
/S15 { BL [] 0 setdash 2 copy Bsquare fill Bsquare } bind def
/D0 { gsave translate 45 rotate 0 0 S0 stroke grestore } bind def
/D1 { gsave translate 45 rotate 0 0 S1 stroke grestore } bind def
/D2 { gsave translate 45 rotate 0 0 S2 stroke grestore } bind def
/D3 { gsave translate 45 rotate 0 0 S3 stroke grestore } bind def
/D4 { gsave translate 45 rotate 0 0 S4 stroke grestore } bind def
/D5 { gsave translate 45 rotate 0 0 S5 stroke grestore } bind def
/D6 { gsave translate 45 rotate 0 0 S6 stroke grestore } bind def
/D7 { gsave translate 45 rotate 0 0 S7 stroke grestore } bind def
/D8 { gsave translate 45 rotate 0 0 S8 stroke grestore } bind def
/D9 { gsave translate 45 rotate 0 0 S9 stroke grestore } bind def
/D10 { gsave translate 45 rotate 0 0 S10 stroke grestore } bind def
/D11 { gsave translate 45 rotate 0 0 S11 stroke grestore } bind def
/D12 { gsave translate 45 rotate 0 0 S12 stroke grestore } bind def
/D13 { gsave translate 45 rotate 0 0 S13 stroke grestore } bind def
/D14 { gsave translate 45 rotate 0 0 S14 stroke grestore } bind def
/D15 { gsave translate 45 rotate 0 0 S15 stroke grestore } bind def
/DiaE { stroke [] 0 setdash vpt add M
  hpt neg vpt neg V hpt vpt neg V
  hpt vpt V hpt neg vpt V closepath stroke } def
/BoxE { stroke [] 0 setdash exch hpt sub exch vpt add M
  0 vpt2 neg V hpt2 0 V 0 vpt2 V
  hpt2 neg 0 V closepath stroke } def
/TriUE { stroke [] 0 setdash vpt 1.12 mul add M
  hpt neg vpt -1.62 mul V
  hpt 2 mul 0 V
  hpt neg vpt 1.62 mul V closepath stroke } def
/TriDE { stroke [] 0 setdash vpt 1.12 mul sub M
  hpt neg vpt 1.62 mul V
  hpt 2 mul 0 V
  hpt neg vpt -1.62 mul V closepath stroke } def
/PentE { stroke [] 0 setdash gsave
  translate 0 hpt M 4 {72 rotate 0 hpt L} repeat
  closepath stroke grestore } def
/CircE { stroke [] 0 setdash 
  hpt 0 360 arc stroke } def
/Opaque { gsave closepath 1 setgray fill grestore 0 setgray closepath } def
/DiaW { stroke [] 0 setdash vpt add M
  hpt neg vpt neg V hpt vpt neg V
  hpt vpt V hpt neg vpt V Opaque stroke } def
/BoxW { stroke [] 0 setdash exch hpt sub exch vpt add M
  0 vpt2 neg V hpt2 0 V 0 vpt2 V
  hpt2 neg 0 V Opaque stroke } def
/TriUW { stroke [] 0 setdash vpt 1.12 mul add M
  hpt neg vpt -1.62 mul V
  hpt 2 mul 0 V
  hpt neg vpt 1.62 mul V Opaque stroke } def
/TriDW { stroke [] 0 setdash vpt 1.12 mul sub M
  hpt neg vpt 1.62 mul V
  hpt 2 mul 0 V
  hpt neg vpt -1.62 mul V Opaque stroke } def
/PentW { stroke [] 0 setdash gsave
  translate 0 hpt M 4 {72 rotate 0 hpt L} repeat
  Opaque stroke grestore } def
/CircW { stroke [] 0 setdash 
  hpt 0 360 arc Opaque stroke } def
/BoxFill { gsave Rec 1 setgray fill grestore } def
/Symbol-Oblique /Symbol findfont [1 0 .167 1 0 0] makefont
dup length dict begin {1 index /FID eq {pop pop} {def} ifelse} forall
currentdict end definefont
end
}}%
\GNUPLOTpicture(2519,1728)
{\GNUPLOTspecial{"
gnudict begin
gsave
0 0 translate
0.100 0.100 scale
0 setgray
newpath
1.000 UL
LTb
450 300 M
31 0 V
1988 0 R
-31 0 V
450 447 M
31 0 V
1988 0 R
-31 0 V
450 562 M
31 0 V
1988 0 R
-31 0 V
450 655 M
31 0 V
1988 0 R
-31 0 V
450 734 M
31 0 V
1988 0 R
-31 0 V
450 803 M
31 0 V
1988 0 R
-31 0 V
450 863 M
31 0 V
1988 0 R
-31 0 V
450 917 M
63 0 V
1956 0 R
-63 0 V
450 1273 M
31 0 V
1988 0 R
-31 0 V
450 1481 M
31 0 V
1988 0 R
-31 0 V
450 1628 M
31 0 V
1988 0 R
-31 0 V
450 300 M
0 63 V
0 1265 R
0 -63 V
695 300 M
0 31 V
0 1297 R
0 -31 V
839 300 M
0 31 V
0 1297 R
0 -31 V
941 300 M
0 31 V
0 1297 R
0 -31 V
1020 300 M
0 31 V
0 1297 R
0 -31 V
1084 300 M
0 31 V
0 1297 R
0 -31 V
1139 300 M
0 31 V
0 1297 R
0 -31 V
1186 300 M
0 31 V
0 1297 R
0 -31 V
1228 300 M
0 31 V
0 1297 R
0 -31 V
1265 300 M
0 63 V
0 1265 R
0 -63 V
1510 300 M
0 31 V
0 1297 R
0 -31 V
1654 300 M
0 31 V
0 1297 R
0 -31 V
1756 300 M
0 31 V
0 1297 R
0 -31 V
1835 300 M
0 31 V
0 1297 R
0 -31 V
1899 300 M
0 31 V
0 1297 R
0 -31 V
1954 300 M
0 31 V
0 1297 R
0 -31 V
2001 300 M
0 31 V
0 1297 R
0 -31 V
2043 300 M
0 31 V
0 1297 R
0 -31 V
2080 300 M
0 63 V
0 1265 R
0 -63 V
2325 300 M
0 31 V
0 1297 R
0 -31 V
2469 300 M
0 31 V
0 1297 R
0 -31 V
1.000 UL
LTb
450 300 M
2019 0 V
0 1328 V
-2019 0 V
450 300 L
1.000 UL
LT0
2106 1515 M
263 0 V
695 1439 M
246 -28 V
143 -194 V
102 14 V
79 1 V
65 -32 V
54 -107 V
47 -72 V
42 -47 V
37 -32 V
34 -35 V
31 -15 V
28 -15 V
27 -5 V
24 -16 V
23 -58 V
21 -53 V
20 -24 V
20 -21 V
18 10 V
17 39 V
17 -9 V
15 2 V
15 -25 V
15 -56 V
14 22 V
13 16 V
13 -22 V
12 38 V
12 -36 V
12 -16 V
11 19 V
11 19 V
11 20 V
10 -6 V
10 -24 V
10 -10 V
9 11 V
9 -2 V
9 -40 V
9 5 V
8 -5 V
9 18 V
8 -9 V
8 -30 V
8 -18 V
7 6 V
8 4 V
7 -12 V
7 -7 V
7 -3 V
7 8 V
7 -9 V
6 -12 V
7 -15 V
6 11 V
6 -16 V
7 -11 V
6 25 V
6 -17 V
6 6 V
5 -19 V
6 -4 V
6 -19 V
5 -9 V
5 -12 V
6 1 V
5 -27 V
5 -18 V
5 -14 V
5 -1 V
5 -31 V
5 3 V
5 -1 V
5 -21 V
4 6 V
5 9 V
5 6 V
4 5 V
4 2 V
5 -17 V
4 -10 V
5 -8 V
4 8 V
4 -5 V
4 4 V
4 -11 V
4 -7 V
4 24 V
4 2 V
4 -9 V
4 5 V
4 4 V
4 2 V
3 -10 V
4 6 V
4 8 V
3 -1 V
4 -1 V
3 -11 V
1.000 UL
LT1
2106 1415 M
263 0 V
695 1532 M
941 1361 L
143 -101 V
102 -71 V
79 -55 V
65 -45 V
54 -38 V
47 -33 V
42 -30 V
37 -26 V
34 -23 V
31 -22 V
28 -20 V
27 -18 V
24 -17 V
23 -16 V
21 -15 V
20 -14 V
20 -13 V
18 -13 V
17 -12 V
17 -12 V
15 -11 V
15 -10 V
15 -10 V
14 -10 V
13 -9 V
13 -9 V
12 -9 V
12 -9 V
12 -8 V
11 -7 V
11 -8 V
11 -8 V
10 -7 V
10 -7 V
10 -6 V
9 -7 V
9 -6 V
9 -7 V
9 -6 V
8 -6 V
9 -6 V
8 -5 V
8 -6 V
8 -5 V
7 -6 V
8 -5 V
7 -5 V
7 -5 V
7 -5 V
7 -5 V
7 -4 V
6 -5 V
7 -5 V
6 -4 V
6 -4 V
7 -5 V
6 -4 V
6 -4 V
6 -4 V
5 -4 V
6 -4 V
6 -4 V
5 -4 V
5 -4 V
6 -3 V
5 -4 V
5 -4 V
5 -3 V
5 -4 V
5 -3 V
5 -4 V
5 -3 V
5 -3 V
4 -4 V
5 -3 V
5 -3 V
4 -3 V
4 -3 V
5 -3 V
4 -3 V
5 -3 V
4 -3 V
4 -3 V
4 -3 V
4 -3 V
4 -3 V
4 -3 V
4 -2 V
4 -3 V
4 -3 V
4 -3 V
4 -2 V
3 -3 V
4 -2 V
4 -3 V
3 -3 V
4 -2 V
3 -3 V
stroke
grestore
end
showpage
}}%
\put(2056,1415){\rjust{LMS approximation}}%
\put(2056,1515){\rjust{$D_N$}}%
\put(1459,50){\cjust{$N$}}%
\put(2080,200){\cjust{ 1e+09}}%
\put(1265,200){\cjust{ 1e+08}}%
\put(450,200){\cjust{ 1e+07}}%
\put(400,917){\rjust{ 0.001}}%
\endGNUPLOTpicture
\endgroup
 